\documentclass[11pt,leqno]{amsart}
\usepackage{amsthm,amsfonts,amssymb,amsmath,oldgerm,upgreek,xcolor}
\numberwithin{equation}{section}
\usepackage{fullpage}
\usepackage{amsmath}
\usepackage{amsfonts}
\usepackage{amssymb}
\usepackage{setspace}
\usepackage{stmaryrd}
\usepackage{mathrsfs}
\usepackage{fancyhdr}
\usepackage{graphicx}
\usepackage[colorlinks = true, linkcolor=blue,
filecolor=darkgreen,
urlcolor=cadmiumgreen,
citecolor=red,pdftex]{hyperref}
\usepackage{psfrag}
\usepackage{listings} 
\usepackage{bbold,rotating,enumitem}
\usepackage[T1]{fontenc}
\usepackage[latin1]{inputenc}
\usepackage[french]{babel}

\definecolor{airforceblue}{rgb}{0.36, 0.54, 0.66}
\definecolor{darkgreen}{rgb}{0.0, 0.2, 0.13}
\definecolor{darkslategray}{rgb}{0.18, 0.31, 0.31}
\definecolor{mediumjunglegreen}{rgb}{0.11, 0.21, 0.18}
\definecolor{prussianblue}{rgb}{0.0, 0.19, 0.33}
\definecolor{warmblack}{rgb}{0.0, 0.26, 0.26}
\definecolor{cadmiumgreen}{rgb}{0.0, 0.42, 0.24}
\definecolor{mauvetaupe}{rgb}{0.57, 0.37, 0.43}
\definecolor{maroon(x11)}{rgb}{0.69, 0.19, 0.38}

\usepackage[top=30mm,bottom=30mm,left=25mm,right=25mm,a4paper]{geometry}
\newenvironment{prof}{\noindent {\sl{Proof.}}}{
\vskip -2mm
\hfill $\Box$ 
}

\def\eps{\varepsilon }

\newcommand\br{\begin{remark}}
\newcommand\er{\end{remark}}
\newcommand\bp{\begin{pmatrix}}
\newcommand\ep{\end{pmatrix}}
\newcommand{\be}{\begin{equation}}
\newcommand{\ee}{\end{equation}}
\newcommand{\ba}[1]{\begin{array}{#1}}
\newcommand{\ea}{\end{array}}

\newcommand{\beg}{\begin{example}}
\newcommand{\eeg}{\end{exaplem}}
\newcommand{\bpr}{\begin{proposition}}
\newcommand{\epr}{\end{proposition}}
\newcommand{\bt}{\begin{theorem}}
\newcommand{\et}{\end{theorem}}
\newcommand{\bc}{\begin{corollary}}
\newcommand{\ec}{\end{corollary}}
\newcommand{\bl}{\begin{lemma}}
\newcommand{\el}{\end{lemma}}
\newcommand{\bd}{\begin{definition}}
\newcommand{\ed}{\end{definition}}
\newcommand{\brs}{\begin{remarks}}
\newcommand{\ers}{\end{remarks}}

\newtheorem{theo}{Theorem}
\newtheorem{prop}[theo]{Proposition}
\newtheorem{cor}[theo]{Corollary}
\newtheorem{lem}[theo]{Lemma}
\newtheorem{defi}[theo]{Definition}
\newtheorem{ass}[theo]{Assumption}

\numberwithin{equation}{section}

\def\eps {\varepsilon}

\def\part{\partial}

\newcommand{\CC}{{\mathbb C}}

\newcommand{\EE }{{\mathbb E}}

\newcommand{\NN}{{\mathbb N}}

\newcommand{\RR}{{\mathbb R}}
\newcommand{\TT}{{\mathbb T}}

\newcommand\cB{{\mathcal B}}

\newcommand\cF{{\mathcal F}}

\newcommand\cL{{\mathcal L}}
\newcommand\cM{{\mathcal M}}
\newcommand\cN{{\mathcal N}}

\newcommand\cP{{\mathcal P}}

\newcommand\cR{{\mathcal R}}

\newcommand\cT{{\mathcal T}}

\newcommand\cX{{\mathcal X}}

\title{\large The Relativistic Vlasov-Maxwell system \\
{\small Local smooth solvability for a weak topology}}
\author[C. Cheverry]{Christophe Cheverry \& Slim Ibrahim}
\address[Christophe Cheverry]{Univ Rennes, CNRS, IRMAR - UMR 6625, F-35000 Rennes, France}
\address[Slim Ibrahim]{Department of Mathematics and Statistics, University of Victoria, British Columbia, Canada}

\begin{document}

\maketitle

\noindent \textbf{\small Abstract}. \small{This article is devoted to the Relativistic Vlasov-Maxwell system in 
space dimension three. We prove the local existence and uniqueness of solutions for initial data $ ({\rm f}_0,
{\rm E}_0,{\rm B}_0) \in L^\infty \times H^1 \times H^1 $, with $ {\rm f}_0 $ compactly supported in momentum.
This result is the consequence of the local smooth solvability for the ``weak" topology associated with 
$ L^\infty \times H^1 \times H^1 $. It is derived from a representation formula decoding how the momentum 
spreads and revealing that the domain of influence in momentum is controlled by mild information (involving 
only conserved quantities). We do so by developing a Radon Fourier analysis on the RVM system, leading to 
the study of a class of singular weighted integrals. In parallel, we implement our method to construct smooth 
solutions to the RVM system in the regime of dense, hot and strongly magnetized plasmas.}

\bigskip


\bigskip

{\small \parskip=1pt
\tableofcontents
}


\section{Main results} \label{introduction} Let  $ C^n_c (\RR^m) $ be the class of compactly supported 
$ n $ times continuously differentiable functions on $ \RR^m $, endowed with the norm
\[ \parallel {\rm U} \parallel_n := \sup \ \bigl \lbrace \parallel \partial^\alpha {\rm U} \parallel_{L^\infty 
(\RR^m)} \, ; \, \vert \alpha \vert \leq n \bigr \rbrace \, , \qquad \parallel \cdot \parallel_0 \equiv \parallel 
\cdot \parallel_{L^\infty} \, , \qquad (n,m) \in \NN^2 . \]
Consider a Cauchy problem which is locally well-posed in $ C^n_c $ with $ n \in \NN^* $. Given some 
initial data $ {\rm U}_0 \in C^n_c $, denote by $ T ({\rm U}_0) \in \RR_+^* \cup \{ + \infty \} $ the lifespan 
of the associated smooth (say $ C^1 $) solution. Select a subspace $ \mathfrak N \subset C^n_c $ which 
is equipped with a \og small \fg{} norm $ \cN $, weaker than $ \parallel \cdot \parallel_n $ meaning that $ \cN \lesssim 
\parallel \cdot \parallel_n $ on $ \mathfrak N $. From there, for all $ {\rm S}_0 \in \RR_+ $, we can introduce 
the lower bound
\begin{equation}\label{lower bound fo T} 
 T(\mathfrak N, \cN;{\rm S}_0) := \inf \ \bigl \lbrace T ({\rm U}_0)  \, ; \, {\rm U}_0 \in \mathfrak N \text{ and } 
\cN({\rm U}_0) \leq {\rm S}_0 \bigr \rbrace  \in \RR_+ \cup \{+\infty \} . 
\end {equation}

\begin{defi}[Smooth solvability of initial value problems] \label{initialvalueproblem} We say that a Cauchy 
problem is locally smoothly solvable (LSS) for $ (\mathfrak N, \cN) $ when $ T(\mathfrak N, \cN;{\rm S}_0) $ 
is a positive number for all $ {\rm S}_0 \in \RR_+ $. It is globally smoothly solvable (GSS) for $ (\mathfrak N, 
\cN) $ when $ T(\mathfrak N, \cN;{\rm S}_0) = + \infty $ for all $ {\rm S}_0 \in \RR_+ $.
\end{defi}

\noindent The notion of smooth solvability for $ (\mathfrak N, \cN) $ focuses on the persistence of regularity. 
It must be clearly distinguished from the classical concept of \og well-posedness \fg{}, which is usually understood 
in the sense of Hadamard that is with existence, uniqueness and also stability related to the normed space 
associated with $ \cN $. It is less demanding in the sense that it does not necessarily require the continuous 
dependence with respect to the \og small \fg{} norm $ \cN $. Such an alleviation may have an important implication 
to analyze systems in turbulent regimes.

\noindent For a large class of first-order quasilinear symmetric hyperbolic systems (QSH systems) whose 
prototypes are Burgers' equation and the 3D compressible Euler equations (away from the vacuum), the Cauchy 
problem  is known to be well posed in $ (C^n_c, \parallel \cdot \parallel_n ) $ as long as $ n$ is large enough. 
But, due to the finite-time singularity formation, it is certainly not globally smoothly solvable for such 
$ (C^n_c, \parallel \cdot \parallel_n ) $. A simple scaling argument indicates that it cannnot be LSS for 
$ (C^n_c, \parallel \cdot \parallel_0 ) $.

\noindent For another whole range of (three-dimensional) nonlinear equations including the incompressible 
Euler equations, the Navier-Stokes equations and the Relativistic Vlasov-Maxwell system (RVM system in 
abbreviated form), we have the local well-posedness (at least for adequate \og large \fg{} norms), but 
the global well-posedness is still open. Standard results involving norms $ \parallel \cdot \parallel_n $ with 
$ n $ large provide with a lower bound for $ T ({\rm U}_0) $ which is typically of the form $ \parallel {\rm U}_0 
\parallel_n^{-1} $. Starting from there, the local smooth solvability for weaker norms $ \cN $ (with $ \cN 
\lesssim \parallel \cdot \parallel_n $) may not be granted. It is recognized that the challenge behind this 
issue (and behind GSS) lies in the development of turbulence.

\noindent This text is devoted to the three-dimensional RVM system. The local Hadamard well-posedness for 
the norm $ C^1 \times C^2 \times C^2 $ is known \cite{MR816621,MR1039487,zbMATH03951266}. The same 
applies to the global weak existence \cite{MR1003433,Rhein} (without uniqueness and stability)  for data in 
$ L^\infty \times L^2 \times L^2 $ (with finite energy). Our goal here is to make a significant progress at the 
interface. It is to show the local smooth solvability for the norm $ L^\infty \times H^1 \times H^1 $. As a 
consequence, we can keep the local existence with the regularity $ L^\infty \times H^1 \times H^1 $ while 
conserving the uniqueness (Corollary \ref{strong-weak}).

\noindent The present section summarizes our outcomes. In Subsection \ref{momentum spread by m}, we 
specify our result (Theorem \ref{inimaintheo}) on LSS. In Subsection \ref{dhsmplasmas}, the focus is on dense, 
hot, collisionless and strongly  magnetized plasmas. As a by-product of our analysis, we are able to construct a 
framework within which the local uniform existence of smooth solutions may be achieved while a large magnetic 
field is applied (Theorem \ref{mainmaintheo}). This furnishes mathematical tools for a better understanding of 
plasmas during their confinement (like in fusion devices), including the several complicated behaviors that can occur.


\subsection{On the Cauchy problem for the RVM system in a strong-weak setting}\label{momentum spread by m}
The unknowns are the distribution function $ {\rm f} $ and the electromagnetic field $ ({\rm E} , {\rm B} ) $ 
depending on the time $ t \in \RR $, on the spatial position $x\in\RR^3$ and on the momentum $ \xi \in \RR^3 $ as 
indicated below 
\[ {\rm f} : \RR_t \times \RR^3_x \times \RR^3_\xi \rightarrow \RR_+ , \qquad {\rm E} : \RR_t \times \RR^3_x 
\rightarrow \RR^3 , \qquad {\rm B} : \RR_t \times \RR^3_x \rightarrow \RR^3 . \]
The speed of light is normalized to one. We deal with the relativistic velocity 
\[ \nu (\xi) := \langle \xi \rangle^{-1} \, \xi , \qquad \langle \xi \rangle := \sqrt {1+ \vert \xi \vert^2} , \qquad \vert 
\nu (\xi) \vert < 1 . \]
The RVM system is composed of the Vlasov equation
\begin{equation}\label{VlasoveqRVMD} 
\quad \partial_t {\rm f} + \nu (\xi) \cdot \nabla_x {\rm f} + {\rm F} (t,x,\xi) \cdot \nabla_\xi {\rm f} = 0 \, , \qquad {\rm F} 
(t,x,\xi) :=  {\rm E} (t,x) + \nu(\xi) \times {\rm B} (t,x)
\end{equation}
coupled with Maxwell's equations
 \begin{subequations}\label{maxwelleqRVMD} 
 \begin{align}
 & \displaystyle \partial_t {\rm E} - \nabla_x \times {\rm B} = {\rm J} := - \int_{\RR^3} \nu(\xi) \ {\rm f} (t,x,\xi) \, d \xi , & \nabla_x 
 \cdot {\rm E} = \int_{\RR^3} {\rm f} (t,x,\xi) \, d \xi - \rho ,\label{maxwelleqRVMD1} \\
 & \displaystyle \partial_t {\rm B} + \nabla_x \times {\rm E} = 0 , &  \nabla_x \cdot {\rm B} = 0 . \qquad \quad \qquad 
 \qquad \ \label{maxwelleqRVMD2} 
 \end{align}
 \end{subequations}  
 
 \noindent The system (\ref{VlasoveqRVMD})-(\ref{maxwelleqRVMD}) on $ {\rm U} := ({\rm f},{\rm E},{\rm B}) $ is
 completed with  initial data
\begin{equation}\label{stationarystate} 
{\rm U}_{\mid t=0} = {\rm U}_0 = ({\rm f}_0,{\rm E}_0,{\rm B}_0) ,
\end{equation}  
satisfying, for some given function $ \rho(x) \in C^1_c (\RR^3;\RR_+) $, the compatibility conditions
\begin{equation}\label{compatibility conditionsrm} 
\nabla_x \cdot {\rm E}_0 = \int_{\RR^3} {\rm f}_0 (x,\xi) \, d \xi - \rho \, , \qquad \nabla_x \cdot {\rm B}_0 = 0 . 
\end{equation} 

 \begin{theo} \label{inimaintheo} Fix $ {\rm P}_0 \in \RR_+ $, and denote by $ B(0,{\rm P}_0] $ the closed ball of $ \RR^3 $
 with radius $ {\rm P}_0 $. Select any Sobolev exponent $ \bar p \in ]3/2,2] $. The RVM is locally smoothly solvable for 
 the subspace
 \begin{equation}\label{mathfrak N} 
\quad \ \mathfrak N := \bigl \lbrace {\rm U}_0 \in C^1_c (\RR^3 \times \RR^3) \times C^2_c (\RR^3) \times C^2_c (\RR^3) \, ; 
 \, \text{\rm supp} \, {\rm f}_0 \subset \RR^3 \times B(0,{\rm P}_0] \bigr \rbrace 
 \end{equation} 
 equipped with the product norm
 \begin{equation}\label{mathcal N}  
 \cN ({\rm U}_0) := \parallel {\rm f}_0 \parallel_{L^\infty(\RR^3 \times \RR^3)} + \parallel {\rm E}_0 \parallel_{W^{1,\bar p} 
 (\RR^3)} + \parallel {\rm B}_0 \parallel_{W^{1,\bar p}  (\RR^3)} . \qquad \qquad
 \end{equation} 
   \end{theo} 

\noindent At the level of (\ref{mathfrak N}) and (\ref{mathcal N}), the topologies related to $ {\rm f}_0 $ and $ ({\rm E}_0,{\rm B}_0) $ 
are managed separately. The choice of $ C^1 $ for $ {\rm f}_0 $ and $ C^2 $ for $ ({\rm E}_0,{\rm B}_0) $ inside (\ref{mathfrak N}) is 
classical \cite {MR816621}. Passing to (\ref{mathcal N}), each of these indices go down one level. To put into perspective Theorem 
\ref{inimaintheo}, as a first step, we can work in a simplified context. With this in mind, we handle initial data which 
are \underline{r}estricted to
\[ \mathfrak N_r := \bigl \lbrace {\rm U}_0 \, ;  \, {\rm f}_0 \in C^1_c (\RR^3 \times \RR^3) \, , \, \text{supp} \ {\rm f}_0 
\subset \RR^3 \times B(0,{\rm P}_0] \, , \, {\rm E}_0 \equiv 0 \, , \, {\rm B}_0 \equiv 0 \bigr \rbrace \subset  
\mathfrak N . \]
Observe that in such a case, a non trivial electromagnetic field $ ({\rm E} , {\rm B} ) $ is generated by interaction. Now, 
tested on $ \mathfrak N_r  $, 
the norm $ \cN $ reduces to the sup-norm (on the sole component $ {\rm f}_0 $). As a direct consequence of Theorem 
\ref{inimaintheo}, the RVM system is locally solvable for $ (\mathfrak N_r, \parallel \cdot \parallel_0 ) $. Since the RVM 
system shares many commonalities with QSH systems, it is instructive to compare this information with what is 
obtained for hyperbolic systems. Starting from a smooth solution leading to blow-up in finite time, a scaling argument 
indicates that general QSH systems (with constant coefficients) are not locally solvable for $ (C^n_c, \parallel \cdot 
\parallel_0 ) $. Why is there such discrepancy ? 

\noindent This is because there are also notable differences separating the RVM system from common QSH systems. 
First, the RVM system is not scaling invariant (due to the Lorentz factor). Secondly, the nonlinearity occurs only on a transport 
part which preserves the sup-norm. There is however a coupling (through the integral term $ {\rm J} $) which produces a
self-consistent electromagnetic field. The feedback of this $ ({\rm E} , {\rm B} ) $ on $ {\rm f} $ may destroy the 
$ C^1 $-regularity. But this impact is lessened due to a number of specificities, among which relativistic effects and transfers 
from time (in $ t $) and space (in $ x $) derivatives to kinetic (in $ \xi $)  derivatives, which can be neutralized at the 
level of the electric current $ {\rm J } $ through arguments from \cite{Bouchut,MR2124491,MR1877669}. 

\smallskip

\noindent For data in $ \mathfrak N_r $, we simply work with 
\begin{equation}\label{canassumethat} 
{\rm U}_{\mid t=0} = {\rm U}_0 = ({\rm f}_0,{\rm E}_0,{\rm B}_0) = ( f_0 , 0, 0 ) ,
\end{equation} 
where $ f_0 $ satisfies (for some $ {\rm P}_0 \in \RR_+ $ and $ {\rm S}_0 \in \RR_+ $)
 \begin{subequations}\label{cestpourf0} 
 \begin{align}
 & \displaystyle \ f_0 \in C^1_c (\RR^3 \times \RR^3 ; \RR)  , \label{cestpourf0r} \\
 & \displaystyle \ \text{supp} \ f_0 \subset \RR^3  \times B(0,{\rm P}_0] , \label{cestpourf0s} \\
 & \displaystyle \parallel f_0 \parallel_0 := \! \! \sup_{(x,\xi) \in \RR^3 \times \RR^3} \vert f_0 (x,\xi) \vert \leq {\rm S}_0 , 
 \label{cestpourf0infty} 
 \end{align}
 \end{subequations}  
together with the compatibility condition inherited from (\ref{compatibility conditionsrm}) and (\ref{canassumethat}), that is
\begin{equation}\label{compatibility conditionsrmf0} 
0 = \int_{\RR^3} f_0 (x,\xi) \, d \xi - \rho . 
\end{equation} 
The initial data $ f_0 $, while remaining smooth, can undergo large fluctuations allowing to trigger filamentation or coherent 
structures \cite{MR3582249,gyrokineticCF}, whose counterparts in the fluid description are shearing and cascade of phases. 
The point is that the regularity persists locally in time, uniformly with respect to the small norm $ \cN $ (instead of $ \parallel 
\cdot \parallel_1 $). In other words, despite arbitrarily large derivatives of $ f_0 $, a control remains within reach. Basically, 
this relates to the minimal radius $ {\rm P} (t) $ of the balls  containing the momentum support of $ {\rm f}(t,\cdot) $, i.e.:
  \begin{equation}\label{minimal radius of the balls} 
\qquad {\rm P} (t) := \inf \ \bigl \lbrace R \in \RR_+ \, ; \, {\rm f} (t,x,\xi) = 0 \ \text{for all} \  x \in \RR^3 \ \text{and for all} \ \xi 
\in \RR^3 \ \text{with} \ R \leq \vert \xi \vert \bigr \rbrace .
\end{equation}
The role of $ {\rm P} $ is essential in the forthcoming discussion. From the definition (\ref{mathfrak N}) of $ \mathfrak N $, 
we can assert that $ {\rm P} (0) \leq {\rm P}_0 $. In what follows, we will exhibit (Proposition \ref{theorempms}) a positive 
increasing continuous function $ \cF $ depending only on $ {\rm S}_0 $, allowing to control $ {\rm P} $ as indicated in 
(\ref{conofcpbis}). Now, from Glassey-Strauss' conditional theorem \cite{MR816621}, singularities may develop only if 
$ {\rm P} (t) $ explodes. Looking at (\ref{deefdeQFF}), this may occur only if $ \parallel {\rm E}(t,\cdot) \parallel_{L^\infty} $ 
goes to $ + \infty $. This explosion mechanism does not necessarily apply to  $ \parallel {\rm f} (t,\cdot) \parallel_{W^{1,\infty}} $. 
Large values of  $ \parallel {\rm f} (t,\cdot) \parallel_{W^{1,\infty}} $ are not known to be a triggering factor of the rapid 
breakdown of smooth solutions. In line with this, Theorem \ref{inimaintheo} shows that large Lipschitz norms of $ f_0 $ 
are compatible with a uniform life-span. Another way to put Theorem \ref{inimaintheo} in context is to derive a priori 
sup-norm estimates on all fluid quantities (Corollary \ref{corrempms}), without the need for looking at derivatives of $ f_0 $.

  \noindent In fact, our approach brings into play strong and  weak characteristics: {\it strong} in view of the (technical) regularity  
 assumption (\ref{cestpourf0r}), {\it weak} in the sense of the relaxed conditions (\ref{cestpourf0s}) and (\ref{cestpourf0infty}). 
 Now, by compactness arguments, we can exhibit a notion of local well-posedness (without stability) in the framework of 
 $ L^\infty \times H^1 \times H^1$. Such \og strong-weak \fg{} solutions satisfy (\ref{conofcpbis}) and they are unique, see 
 Corollary \ref{strong-weak}. Note that the existence results \cite{MR1003433,Rhein} of weak solutions do not provide information 
 about the persistence of regularity nor about the property (\ref{conofcpbis}). On the other hand, little is known about the 
 uniqueness of these weak solutions, see \cite{MR3356994} for an overview and a result (in the absence of coupling).

 \noindent Theorem \ref{inimaintheo} is derived from Sections \ref{Repformulamomentum} and \ref{RadonFourier}. It is a compilation 
 of two intermediate important stages. First, in Section \ref{Repformulamomentum}, we bring to the fore a representation  formula 
 for the momentum increment (Definition \ref{momentum spread}) with some original proof (resorting to a Radon transform 
 \cite{MR1723736} with respect to $ x \in \RR^3 $ and oscillatory integrals). Second, in Section \ref{RadonFourier}, we control the 
 momentum spread. To this end, we have to study a class of singular weighted integrals. A few ideas leading to (\ref{conofcpbis}) 
 are new in comparison to the preceding approaches (presented in Subsection \ref{precedingstrategies}). They are exposed in 
 Subsection \ref{Plan of the text}.


\subsection{Application: the regime of dense, hot and strongly magnetized plasmas}\label{dhsmplasmas} In kinetic theory, the time evolution 
of a single species of charged particles (typically electrons in a background of ions) is described (in dimensionless units) by the 
system (\ref{VlasoveqRVMD})-(\ref{maxwelleqRVMD}). The study of (\ref{VlasoveqRVMD})-(\ref{maxwelleqRVMD}) is of fundamental 
interest because it reveals the interactions between the matter (driven by $ {\rm f} $) and the fields (represented by $ {\rm E} $ and 
$ {\rm B} $). In (\ref{maxwelleqRVMD1}), the function $ \rho \in C^1_c (\RR^3;\RR_+) $ stands for the density of charges associated 
with protons. 

\noindent Property (\ref{conofcpbis}), which prevails for strong and (adequate) weak solutions, is likely to furnish a wide 
range of applications. It is related to precise quantitative information on the lifespan, and therefore it is well suited to the study 
of the RVM system at specific scales. Particular emphasis is placed here on situations involving (after nondimensionalization) 
large magnetic fields. As a matter of fact, such models constitute a real challenge in plasma physics (Paragraph \ref{physical model}) 
and provide means (Paragraph \ref{Perturbation theory}) to generate solutions undergoing rapid oscillations.

\subsubsection{The physical model} \label{physical model} In many important applications such as planetary magnetospheres 
(like the Van Allen Belts in the Earth situation) or fusion devices (like tokamaks), a large inhomogeneous external magnetic field 
$ \eps^{-1} \, {\rm B}_e(x) $ is applied. Here, the number $ \eps $ is the inverse of the electron gyrofrequency; it is a small 
dimensionless parameter which is typically in the order of $ \eps \approx 10^{-5} $. By extension and to highlight the smallness 
of $ \eps $, we will work with $ \eps \in ]0,1] $. Since the large weight $ \eps^{-1} $ is in factor of $ {\rm B}_e $, the plasma is 
{\it strongly magnetized} as soon as $ {\rm B}_e \not \equiv 0 $. The function $ {\rm B}_e(\cdot)$ often takes the form of a smooth bounded solenoidal 
and irrotational vector field:
\begin{equation}\label{conditionforthebe} 
{\rm B}_e \in C^1_b ( \RR^3 ; \RR^3 ) \, , \qquad \nabla_x \cdot {\rm B}_e \equiv 0 \, , \qquad \nabla_x \times 
{\rm B}_e \equiv 0 . 
\end{equation}
By this way, the given external magnetic field $ {\rm B}_e $ does fit in with the second condition inside 
(\ref{compatibility conditionsrm}), which means that it induces zero current. The variations of $ {\rm B}_e $ 
are quite important since they account for the spatial inhomogeneities which are usually issued from the 
underlying physical geometries (such as toroidal shapes). Physical plasmas are generally comprised of 
a dominant part which stays  at the thermodynamic equilibrium. In the relativistic framework, the reference 
model at rest  is the Maxwell-J\"uttner distribution. We consider here that most of charged particles are in 
a steady state. This may be represented by a distribution profile for $ {\rm f} (\cdot) $ of the form $ {\rm M} 
(\eps ,\langle \xi \rangle) $, which may depend smoothly on $ \eps \in [0,1] $. More precisely 
\begin{equation}\label{conditionforthedenM} 
 {\rm M} \equiv {\rm M}^\eps  \equiv {\rm M} (\eps , \langle \xi \rangle) \, , \qquad {\rm M} (\eps,r) \in C^1_c
 ( [0,1] \times [1,+\infty [ ;\RR_+) . 
\end{equation}
The plasma is called {\it dense} when there exists a constant $ c_1 \in \RR_+^* $ such that
\begin{equation}\label{denseconditionfornM} 
 0 < c_1 \leq \parallel {\rm M} (\eps , \cdot) \parallel_\infty \, , \qquad \forall \eps \in [0,1] .  
 \end{equation}
It is called {\it hot} when there exists a constant $ c_2 \in \RR_+^* $ such that
\begin{equation}\label{hotconditionfornM} 
 \text{supp} \ {\rm M} (\eps , \cdot) \cap [1+c_2,+\infty[ \not = \emptyset \, , \qquad \forall \eps \in [0,1] .  
\end{equation}
\noindent To fit with (\ref{compatibility conditionsrm}), we must assume that $ \rho \equiv \rho^\eps $ is a 
constant adjusted in such a way that
\begin{equation}\label{basiccompatibility conditionsrm} 
\rho^\eps = \int_{\RR^3} {\rm M}(\eps , \langle \xi \rangle) \, d \xi \, , \qquad \forall \eps \in [0,1] . 
\end{equation} 
Then, for all parameter $ \eps \in ]0,1] $, the expression
\begin{equation}\label{prototype}  
\tilde {\rm U}_a^\eps (t,x,\xi) := \bigl( {\rm M} (\eps,\langle \xi \rangle), 0, \eps^{-1} \, {\rm B}_e(x) \bigr) \equiv 
\tilde {\rm U}_a^\eps (0,x,\xi) 
\end{equation}  
is clearly a stationary solution to the RVM system. When $ {\rm M} \not \equiv 0 $ does not depend on $ \eps $ and 
$ {\rm B}_e \not \equiv 0 $, it is the prototype of a solution belonging to the regime of dense, hot and strongly 
magnetized plasmas. It can serve as a working example. Retain however that we deal with a more general class 
of \underline{\it a}pp{\it roximate solutions} denoted by $ {\rm U}_a^\eps = ( {\rm f}_a^\eps , 
{\rm E}_a^\eps , {\rm B}_a^\eps ) $, adjusted as follows.

\begin{ass}\label{assumptiononoUepsa}
The family $ \{ {\rm U}_a^\eps \}_\eps $ is \underline {well}-p\underline {re}p\underline {ared} in the sense of 
Definition \ref{Local well-prepared approximate solution}. 
\end{ass}

\noindent The reader is referred to Paragraph \ref{Well-prepared approximate solutions} for a precise description and a 
discussion about the content of $ {\rm U}_a^\eps $. The expression $ {\rm U}_a^\eps $ does not necessarily have to be 
an exact solution to (\ref{VlasoveqRVMD})-(\ref{maxwelleqRVMD}). It may produce a non-zero remainder $ {\rm R}_a^\eps $. 
It may also reveal a full range of plasma distinctive features which are not detected by $ \tilde {\rm U}_a^\eps $, see for instance 
\cite{MR3582249,gyrokineticCF} for a preview of what can happen. From now on, we fix a family 
$ \{ {\rm U}_a^\eps \}_\eps $ of approximate solutions and, at the initial time $ t = 0 $, we modify $ {\rm U}_a^\eps (0,\cdot) $ 
according to 
\begin{equation}\label{stationarystateperturb} 
{\rm U}_0^\eps = {\rm U}_a^\eps (0,\cdot) + U_0^\eps \, , \qquad U_0^\eps = (f_0^\eps,0,0) , \qquad \eps \in ]0,1] .
\end{equation}  
In (\ref{stationarystateperturb}), there is no initial electromagnetic perturbation $ (E_0^\eps,B_0^\eps) $ interfering 
with $ {\rm U}_a^\eps $. This is because, under adequate assumptions that will be investigated in Paragraph 
\ref{Well-prepared approximate solutions}, the influence of such $ (E_0^\eps,B_0^\eps) $ can be incorporated inside 
$ ({\rm E}_a^\eps , {\rm B}_a^\eps) (0,\cdot) $.

\subsubsection{Perturbation theory} \label{Perturbation theory} The realistic plasmas always include a larger or smaller 
part of matter which is out of equilibrium and which may have destabilizing effects (like electron beams). Such aspects 
can be taken into account by introducing some $ f_0^\eps \not \equiv 0 $ inside (\ref{stationarystateperturb}). 

\begin{ass}\label{assumptiononofoepsprep} 
The family $ \{ f_0^\eps \}_\eps $ is \underline {com}p\underline {atible} in the sense of Definition \ref{assumptiononofoeps}.
\end{ass}

\noindent Basically, Assumption \ref{assumptiononofoepsprep} means that we impose (\ref{cestpourf0})-(\ref{compatibility conditionsrmf0})
uniformly on $ \{ f_0^\eps \}_\eps $, and moreover that we add the $ L^2 $-smallness condition (\ref{absorptioforsm}) which is 
calibrated by $ {\rm N}_0 $ given in Definition \ref{assumptiononofoeps}. Below, we fix the values of $ {\rm P}_0 \in \RR_+^* $, 
$ {\rm S}_0 \in \RR_+^* $ and  $ {\rm N}_0 \in \RR_+^* $ occurring in Definition \ref{assumptiononofoeps}, and we select $ \{ f_0^\eps \}_\eps  $ 
accordingly. Then, we seek the solution to 
(\ref{VlasoveqRVMD})-(\ref{maxwelleqRVMD}) in the form
\begin{equation}\label{typeofsolution} 
{\rm U}^\eps (t,x,\xi) = {\rm U}_a^\eps (t,x,\xi) + U^\eps (t,x,\xi) , \qquad U^\eps = (f^\eps,E^\eps,B^\eps) .
\end{equation}  

 \begin{theo}[Uniform life-span for classical solutions] \label{mainmaintheo} Fix a well-prepared family $ \{ {\rm U}_a^\eps \}_\eps $ 
 of approximate solutions. Select a compatible family $ \{ f_0^\eps \}_\eps $ of initial data controlled by $ {\rm P}_0 $, $ {\rm S}_0 $ 
 and  $ {\rm N}_0 $. Then, there exists a positive time $ \cT $ depending 
 only on $ ({\rm P}_0,{\rm S}_0,{\rm N}_0) $, say $ \cT \equiv \cT ({\rm P}_0, {\rm S}_0,{\rm N}_0) \in \RR_+^* $, such that, for all 
 $ \eps \in ]0,1] $, the Cauchy problem (\ref{VlasoveqRVMD})-(\ref{maxwelleqRVMD})-(\ref{stationarystate}) with $ {\rm U}_0^\eps $ 
 as in (\ref{stationarystateperturb}) has a smooth $ C^1_c $-solution on the domain $ [0,\cT] \times \RR^3 
 \times \RR^3 $. Moreover, there exists a positive continuous function $ \cF $ depending only on $ ({\rm P}_0,{\rm S}_0,{\rm N}_0) $, 
 say $ \cF \equiv \cF ({\rm P}_0,{\rm S}_0,{\rm N}_0) : [0,\cT] \rightarrow  \RR_+^* $, such that
\begin{equation}\label{conofcp} 
\text{\rm supp} \ f^\eps (t,\cdot) \subset \RR^3 \times B \bigl(0,\cF({\rm P}_0,{\rm S}_0,{\rm N}_0) (t) \bigr) , \qquad 
\forall \, t \in [0,\cT({\rm P}_0,{\rm S}_0,{\rm N}_0)] . 
\end{equation}
 \end{theo}

 \noindent Resorting to approximate solutions $ {\rm U}^\eps_a $ allows to break (\ref{canassumethat}) since 
 $ ({\rm E}_a^\eps ,{\rm B}_a^\eps ) (0,\cdot) $ is aimed at being non-trivial. As can be guessed by looking at (\ref{prototype}), 
 Definition \ref{Local well-prepared approximate solution} gives access to large magnetic fields $ {\rm B}_a^\eps $. 
 That is, however, not the only aspect. Implementing $ {\rm U}^\eps_a $ or $ f_0^\eps $ is also a way to generate rapid fluctuations, 
 and then to  measure their quantitative impact. In Theorem \ref{mainmaintheo}, the absorption of large fields and oscillations 
 requires the  $L^2 $-smallness conditions (\ref{remainderupesa2}) on the remainder $ {\rm R}_a^\eps $ and 
 (\ref{absorptioforsm}) on the initial data $ f_0^\eps $. This  is really not demanding in comparison with the common  
 assumptions (like coherence, $ H^s_\eps $-estimates and so on) in nonlinear geometric optics \cite{MR2562165,Rauch}, 
 in comparison with the already improved contexts of \cite{MR4084146} or \cite{MR4357273}, and in comparison with the 
strong conditions imposed in \cite{MR1086751} and recently in \cite{MR4244270}. 
  
\noindent The regime of dense, hot and strongly magnetized plasmas is of practical interest. For this reason, it has been 
intensively studied in physics. Up to now, the mathematical advances \cite{MR2368978,MR1003433} are based on weak 
solutions which do not confer the ability to describe and control the solutions (due to lack of uniqueness). Theorem \ref{mainmaintheo} 
achieves that goal with smooth solutions on a pertinent observation time ($ t \sim 1 $). It can be interpreted as a kind 
of stability statement about the regularity of solutions near WKB expansions (such as $ {\rm U}^\eps_a $), when some 
large $ C^1 $-perturbations are applied. Theorem \ref{mainmaintheo} is proved in Section \ref{Application}. For more 
details on Theorems \ref{inimaintheo} and \ref{mainmaintheo}, the reader  is recommended to look at next Section \ref{detailintroduction}.


\section{Detailed introduction} \label{detailintroduction} 

The construction of solutions to the RVM system has a long history going back to the pioneering works of Wollman \cite{zbMATH03951266},
Glassey-Strauss \cite{MR816621,MR1039487} and DiPerna-Lions \cite{MR1003433} in the 1980s. It also includes the alternative methods
proposed by Klainerman-Staffilani \cite{MR1877669} and Bouchut-Golse-Pallard \cite{Bouchut,MR2124491} in the early 2000s. This also 
relates to the recent developments \cite{MR4084146,MR4357273,MR1086751,MR4244270} in connection with plasmas involving 
large fields. In this section, we aim to explain in depth how our results and techniques stand in relation to the preceding contributions, 
and also how they differ from them. In Subsection \ref{Historical background}, we recall the historical backdrop. In Subsection \ref{impactofoscillations}, 
we give some insight into localized and oscillating solutions. In Subsection \ref{precedingstrategies}, we list the preceding approaches. In Subsection 
 \ref{Plan of the text}, we compare them with our strategy; we also detail the plan and the content of the text, while giving a sketch 
 of ideas and proofs. 


\subsection{Historical background} \label{Historical background} We give here a quick overview of previous contributions.

\subsubsection{Local and global smooth well-posedness} \label{smooth well-posedness}
The existence of local stable classical solutions dates back to \cite{zbMATH03951266}. The issue of global existence was 
then raised in \cite{MR816621,MR1039487}. Then, much has been done to make progress on the global smooth solvability. 
The advances concern small initial data \cite{MR969207} (much smaller than in Theorem \ref{mainmaintheo}), reduced 
dimensions \cite{MR1463042,MR1620506} (less than the six dimensions of the actual phase space), symmetry conditions 
\cite{Wang} (broken by the inhomogeneities of the external field $ {\rm B}_e $), stability properties \cite{MR1086751} 
(under stronger topologies), and so on. For further information, we can refer to the survey article \cite{MR1379589}. 
The most elaborate result in relation with our theme is perhaps the latest article \cite{MR4244270} of Wei-Yang which 
(among other things) can include a large electromagnetic field but which also requires (in compensation) a very small 
density distribution. In view of \cite{MR4244270}-Proposition 4.8, interpreted in our framework, this density should be of size 
less that $ \eps^8 \, \vert \ln \eps \vert^{-11} $. This is far from the regime under consideration in Theorem \ref{mainmaintheo}. 
Now, for general data, many continuation criteria have been established, see \cite{MR3437855,MR3291372,MR3721415} 
and the references therein. However, the global existence of classical solutions remains open. This is a longstanding and still
active problem. On the other hand, failing to reach global existence, in connection with  concrete applications, we can seek 
for more precise quantitative information on the life-span of smooth solutions. That is the position of the present text.


\subsubsection{About nonlinear geometric optics} \label{geometric optics} The construction of oscillating solutions for
quasilinear systems of conservation laws is the core subject of nonlinear geometric optics. However, the general results
\cite{MR2562165,Rauch} do not furnish optimal information when applied to the RVM system. There are several reasons 
for this. The transport part (the Vlasov equation) and the integral source term (the electric courant) make things easier.
Also, importantly, the variable coefficients with $ \eps^{-1} $ in factor (generated by $ {\rm B}_e $) can be handled through
different arguments than the usual quite restrictive conditions.

\noindent With the goal of better exploiting the specificities of (\ref{VlasoveqRVMD})-(\ref{maxwelleqRVMD}) and also 
motivated by the applications, there has been some progress to remedy this situation. The {\it cold} configuration 
for which $ \vert \xi \vert \lesssim \eps $, or $ {\rm M}^\eps \equiv M(\vert \xi \vert /\eps) $ with $ M \in C^1_c 
(\RR_+;\RR_+) $, is examined in \cite{MR4084146}. The hot but {\it dilute} situation, for which $ {\rm M}^\eps \equiv \eps \, 
M ( \xi ) $ with $ M \in C^1_c (\RR^3;\RR_+)$ and $ (f_0^\eps,  E_0^\eps , B_0^\eps) $ is {\it small in Lipschitz 
norm} (of size $ \eps $), is investigated in \cite{MR4357273}. In \cite{MR4357273}, large amplitude profiles 
$ {\rm M}^\eps \approx 1 $ as well as weak norms of $ U^\eps_0 $ are excluded for deeper reasons related to the 
method used.

\smallskip

\noindent To our knowledge, Theorem \ref{mainmaintheo} is completely new. Implemented in the dense, hot and strongly 
magnetized framework, that is when $ {\rm M} \approx 1 $, $ {\rm P}_0 \approx 1 $ and $ |{\rm B}_e| \approx 1 $, the criteria
yielding global existence are clearly not applicable. On the other hand, all preceding approaches seem to furnish a life-span 
$ \cT^\eps \in \RR_+^* $ that shrinks very rapidly to $ 0 $ when $ \eps $ goes to $ 0 $. 

\subsubsection{Global weak existence}  \label{Global weak well-posedness} As is well-known, weak global solutions 
are available. This has been shown by DiPerna and Lions in the seminal contribution \cite{MR1003433}, just after 
\cite{MR816621,MR1039487} in the late 1980s. But such solutions are obtained by compactness methods, and 
much information has not yet been provided. In the context of (\ref{VlasoveqRVMD})-(\ref{maxwelleqRVMD}), 
little can currently be said about the uniqueness, the stability, the regularity, or the (oscillating) form of these solutions. 
Note however that uniqueness and stability can be addressed from the perspective of the renormalization property
\cite{MR1022305}. We refer to \cite{MR3356994}-Paragraph 1.3 for recent developments (in a decoupled situation)
and a nice presentation of results in this direction. Recall also that weak solutions have been considered in the 
framework (similar to Subsection \ref{dhsmplasmas}) of large magnetic fields in order to give partial information, 
develop asymptotic models or enrich gyrokinetic theory (see for instance \cite{MR2368978}). In view of 
Theorem \ref{mainmaintheo}, some strong control (related to the smoothness) on the solutions does persist. But 
how is this possible, and why$ \, $? The first step (Subsection \ref{impactofoscillations}) is to identify the underlying 
difficulties; the second step (Subsection \ref{precedingstrategies}) is to complement our presentation by recalling the 
previous strategies; the third step (Subsection \ref{Plan of the text}) is to present our plan and to explain our approach.


\subsection{The impact of localizations and oscillations} \label{impactofoscillations} Plasmas out of equilibrium can exhibit 
a wide variety of complicated behaviors, which are manifested by oscillating coherent structures (in subdomains of the phase 
space \cite{MR3582249,gyrokineticCF}) or by chaotic motions (shearing due to the sensitivity of characteristics under 
changes of initial data, especially near separatrices). As a consequence, the Lipschitz norm of $ f^\eps $ is in practice often very 
large. Our aim in this subsection is not to inventory all the phenomena that can occur. Instead, we want to indicate why the 
distribution function $ f^\eps $ should undergo large rapid variations (when $ {\rm B}_e \not \equiv 0 $). For the sake of simplicity, 
we do not (always) mention the dependence on the parameter $ \eps $ at the level of 
$ {\rm U} $, $ U $, $ f $, $ E $ and $ B $. 

\noindent In the simplified framework (\ref{prototype}), emphasis may be placed on the role 
of $ f_0^\eps $, and its impact on $ f $ and on the self-consistent electromagnetic field $ (E,B) $. The presence of some 
$ f_0^\eps \not \equiv 0 $ has many consequences, among which:
 \begin{itemize}
\item [i)] the onset of {\it anisotropic features} (in momentum variable $ \xi $). A spatial localization of $ f_0^\eps $ 
necessarily implies that $ \nabla_x f \not \equiv 0 $. Then, the part $ \nu (\xi) \cdot \nabla_x f \not \equiv 0 $ inside 
(\ref{VlasoveqRVMD}) is switched on, and it cannot be only a function of $ \vert \xi \vert $. The same applies for $ f $. 
\item [ii)] the emergence of {\it fast oscillations}. The Lorentz force $ {\rm F} $ can be decomposed into
\begin{equation}\label{defdeF} 
{\rm F} = F + \eps^{-1} \, \nu(\xi) \times {\rm B_e} (x) \, , \qquad F (t,x,\xi) :=  E(t,x) + \nu(\xi) \times B(t,x) .
\end{equation}
Due to the anisotropy of $ f $, when computing the contribution $ (\xi \times {\rm B_e}) \cdot \nabla_\xi {\rm f} $, the above fast 
rotating term is certainly activated. This means that, as allowed by (\ref{cestpourf0infty}), $ \parallel f \parallel_{W^{1,\infty}} $ 
is designed to be not uniformly bounded with respect to $ \eps \in ]0,1] $. The condition (\ref{cestpourf0infty}) does not preclude the 
choice of initial data $ f_0^\eps $ having arbitrarily large Lipschitz norms (as compared to $ {\rm P}_0 $ and $ {\rm S}_0 $). 
\end{itemize}

\noindent The above considerations deserve to be highlighted by a concrete example (more elaborate models can be found in 
\cite{MR3582249,gyrokineticCF,MR2139205}). Let $ \phi : \RR^3 \rightarrow \RR^3 $ be a global diffeomorphism. Select
\[ f_0^\eps (x,\xi) := \cP \bigl( x,  \eps^{-1} \, \phi(x) , \xi \bigr) \, , \qquad \cP \in C^1_c (\RR^3_x \times \RR^3 
\times \RR^3_\xi ;\RR) . \]
The profile $ \cP $ may be periodic or (partially) compactly supported with respect to its second variable (to represent a localized 
input of electrons for the initialization of an electron beam). For adequate choices of $ \cP $ and $ \phi $, we have clearly access 
to (\ref{absorptioforsm}), and Theorem \ref{mainmaintheo} does apply. Now, we can get a preview of what happens at least when 
neglecting most of the terms inside (\ref{VlasoveqRVMD}). With this in mind, we can consider the elementary transport equation
\begin{equation}\label{simplifiedfmodel} 
\partial_t f + \eps^{-1} \ \bigl( \nu (\xi) \times {\rm B}_e (x) \bigr) \cdot \nabla_\xi f = 0 , \qquad {\rm B}_e (x) = 
b_e (x_3) \ {}^t (0,0,1) .
\end{equation}
In (\ref{simplifiedfmodel}), the field $ {\rm B}_e $ has a fixed (vertical) direction and a varying amplitude $ b_e $ (depending only on $ x_3 $ to satisfy 
$\nabla\cdot B_e=0,\;=\nabla\times B_e=0$) . In the cylindrical coordinate system for $ \xi $, with 
\[ \xi = (r \, \cos \theta,r \, \sin \theta , \xi_3) , \quad \ r^2 = \xi_1^2 + \xi_2^2  , \quad \ \theta \in \TT , \quad \ \tilde \cP (\cdot, r,\theta, 
\xi_3) := \cP (\cdot, r \, \cos \theta,r \, \sin \theta, \xi_3) , \] 
the solution is simply given by
\begin{equation}\label{rapidoscillations} 
\tilde f (x, r,\theta, \xi_3) :=  f(t,x,r \, \cos \theta,r \, \sin \theta, \xi_3) = \tilde \cP \Bigl( x, \frac{\phi(x)}{\eps} ,r,\theta + \frac{t \, b_e
(x_3)}{\eps \, \sqrt{1+r^2}}, \xi_3 \Bigr) . 
 \end{equation}
 The transport part $ \nu (\xi) \cdot \nabla_x f $ of (\ref{VlasoveqRVMD}) has been removed at the level of (\ref{simplifiedfmodel}). 
 This suppresses the first effect i). But still we are faced with ii). Indeed, when $ \partial_\theta \tilde \cP \not \equiv 0 $, the first order 
 derivatives in almost all directions (time, space and momentum) of the expression $ \tilde f $ given by (\ref{rapidoscillations}) are of 
 large size $ \eps^{-1} $. As a consequence, uniform Lipschitz estimates (with respect to $ \eps $) are certainly not available. The same is 
 sure to apply to the solutions provided by Theorem \ref{mainmaintheo}.

\noindent In the dilute situation \cite{MR4357273}, the expression $ {\rm U}^\eps_a $ and the perturbation $ f^\eps_a $
are of small amplitude $ \eps $. This weight $ \eps $ which can be put in factor of the data can be used to absorb the oscillations 
at the frequency $ \eps^{-1} $. This boils down to a sort of weakly nonlinear regime. Then, uniform Lipschitz estimates may become 
available and, by this way, oscillating approximate solutions to (\ref{VlasoveqRVMD}) can be constructed and justified by standard 
arguments. There is nothing like this in the framework of Theorem \ref{mainmaintheo}.

\noindent In fact, the main reason why it is complicated to achieve $ \cT^\eps \sim 1 $ is the following. So far, the methods which 
have been implemented require, at a moment or another, to compute derivatives of $ U $. But this is proving to be very costly in 
terms of negative powers of $ \eps $ because these derivatives are - in all reasonable norms - at least of size $ \eps^{-1} $, 
compromising the existence of a uniform life-span. By contrast, we are able here to avoid this problem. 


\subsection{Previous approaches for smooth solutions} \label{precedingstrategies} 
\noindent Estimates on $ {\rm P} $ (or on similar quantities) have been a central part in the study of the RVM system because 
singularities do not develop as long as the momentum support of $  {\rm f} (t,\cdot) $ remains bounded. This is Glassey-Strauss' 
conditional theorem \cite{MR816621} which has inspired many works and which has been revisited in \cite{Bouchut}. Such bounds 
have been achieved by pursuing three noteworthy lines of research.

\subsubsection{The historical procedure} \label{The historical procedure} As already mentioned, the standard method has been 
initiated in \cite{MR969207,MR816621,MR1039487}. It is based on a representation formula for $ ({\rm E},{\rm B}) $, see 
\cite{MR816621}-Theorem 3 and \cite{MR1620506} for its simplified two-dimensional version. In this line, the stability under 
(very) small perturbations in Lipschitz norm of smooth solutions (such as $ {\rm U}^\eps_a $) has been investigated in \cite{MR1086751}.

\subsubsection{An alternative path} \label{An alternative path} The second way is to  proceed similarly to what has been 
done in Klainerman-Staffilani \cite{MR1877669} or even simpler in Bouchut-Golse-Pallard \cite{Bouchut,MR2124491}. These 
authors have obtained with an economy of means Lipschitz estimates on the field $ ({\rm E},{\rm B}) $. To this end, they have 
implemented three principal arguments: 

\noindent {\it a)} {\it commuting vector field techniques} \cite{MR1877669} for the wave equation related to
Maxwell's equations; 

\noindent {\it b)} a {\it non-resonant smoothing property} \cite{MR2124491} stating that $ 1 \pm \nu(\xi) \cdot \omega $ 
with $ \omega \in \mathbb S^2 $ remains away from zero. This property of ellipticity deteriorates when $ \vert \xi \vert $ 
grows up, and this is one of the difficulties;

\noindent {\it c)} a {\it division lemma}, see \cite{Bouchut}-Lemma 3.1, whose aim is to convert transversal derivatives
to the light cone into the derivative $ \part_t + \nu(\xi) \cdot \nabla_x $ (appearing in the transport part of RVM).

\subsubsection{Global existence through sharp decay estimates} \label{through sharp decay estimates} The third method 
is to implement (by way of a fixed-point iteration) sup-norm controls and decay estimates (related to the dispersion) on 
$ ({\rm E},{\rm B}) $, as was achieved by Wei-Yang in \cite{MR4244270}. But this requires controls in $ C^1 \times 
C^2 \times C^2 $ which are at the origin of the strong smallness condition ($ \leq \eps^8 \, \vert \ln \eps \vert^{-11} $)
 on $ f^\eps_0 $ noted above.

\smallskip

\noindent The three approaches \ref{The historical procedure}, \ref{An alternative path} and \ref{through sharp decay estimates} 
are based on different assumptions; they resort to distinct techniques; they produce complementary information; and they are 
good starting points. But, to achieve our goal, they need to be supplemented by novel procedures. We will adapt these 
methods by following at some points a different logic, yielding other consequences. 


\subsection{Plan of the text, and main ideas} \label{Plan of the text} This article is organized around three sections
(\ref{Repformulamomentum}, \ref{RadonFourier} and \ref{Application}) whose contents are detailed in below.

\subsubsection{Content of Section \ref{Repformulamomentum}}
In Section \ref{Repformulamomentum}, we deal with smooth compactly supported solutions to the RVM system 
(\ref{VlasoveqRVMD})-(\ref{maxwelleqRVMD}). We derive a representation formula (Proposition \ref{theorem1})
which is reminiscent of \cite{MR816621}. In contrast to \cite{Bouchut,MR816621}, we do not look at $ ({\rm E},{\rm B}) $. 
Instead, we focus on the 
{\it momentum increment} $ {\rm D} (t,y,\eta) $, see Definition \ref{momentum spread}. The construction of 
$ {\rm D} (t,y,\eta) $ is related to the characteristic $ (X,\Xi) (t,y,\eta) $ emanating from a point $ (y,\eta) $ 
and associated with the Vlasov equation (\ref{VlasoveqRVMD}). Equivalently, $ {\rm D} (t,y,\eta) $ is the difference 
$ \langle \Xi (t,y,\eta) \rangle - \langle \eta \rangle $ which is built from the solution $ (X,\Xi) $ to the ordinary differential 
equation (\ref{diffeoom}). The quantity $ {\rm D} $ depends on $ t $. It is a dynamical quantity which allows to control 
the size $ {\rm P} (t) $ of the momentum support. Obviously, when $ {\rm E} \equiv 0 $, we find that $ {\rm D} $ is invariant, 
that is $ {\rm D}(t) =  {\rm D}(0) $ for all time $ t $. But certainly this becomes not true in the presence of interactions. This 
raises an interesting question: 

\smallskip

\noindent {\bf Question 1.} {\it Can we identify the factors that alter $ {\rm D} \, $? And evaluate their quantitative effects$ \, $?}

\smallskip

\noindent The determination of $ {\rm D} $ is based on $ \vert \Xi \vert $, not on $ \Xi $. This remark is important because, in the 
regime under consideration, the directions of $ \Xi $ may be strongly oscillating with respect to the parameter $ \eps \in ]0,1] $, 
see for instance (\ref{rapidoscillations}), 
and therefore they would be impractical to control. In fact, resorting on $ {\rm D} $ is a way to avoid rapid fluctuations of 
$ \Xi /  \vert \Xi \vert $, and to dispense with having to implement the (possibly large) momentum speed of propagation. 
The quantity $ {\rm D} $  is recovered after some integration (with underlying cancellation effects), while still informing 
about the expanse of the momentum support.

\noindent To evaluate $ {\rm D} $, a first option is to consider (\ref{deefdeQFF}), and to directly extract the sup-norm of $ E $. 
But this would imply pointwise estimates (like the estimates of von Wahl \cite{MR280885} or the Strichartz inequalities 
\cite{MR1646048,MR1240537}) at the level of a three-dimensional wave equation with source term $ \part_t {\rm J} $. 
This would take us back to the control of derivatives with again important losses in terms of negative powers of $ \eps $.
This is why we adopt an alternative strategy which is partly inspired by \cite{Bouchut,MR3291372}. As in \cite{Bouchut}, 
in Paragraph \ref{Lienard-Wiechertpotentials}, we interpret the RVM system in terms of the Lienard-Wiechert potential 
$ {\rm u} (t,x,\xi) $. As in \cite{MR3291372}, we aim at substituting for (\ref{deefdeQFF}) integral expressions which are 
not based on $ dt $ but instead which integrate with respect to the time, the space and the velocity. But, to this end, we 
exploit other arguments. Our main innovation to prove (\ref{cD}) is:

\vskip 0.4mm

- to perform a Radon transform (with respect to the position variable $ x \in \RR^3 $). The Radon transform (see the book 
\cite{MR1723736} and Paragraph \ref{basictools}) allows to convert $ {\rm f} (t,x,\xi) $ into $ {\rm g} (t,\omega,p,\xi) $ and 
$ {\rm u} (t,x,\xi) $ into $ {\rm v} (t,\omega,p,\xi) $. The equation (\ref{vlasoveqonura}) for $ {\rm g} (\cdot,\omega,\cdot,\xi) $ 
becomes (with respect to $ t $ and $ p $) a one-dimensional transport equation; the equation (\ref{waveeqonura}) for 
$ {\rm v} (\cdot,\omega,\cdot,\xi)$ becomes (with respect to $ t $ and $ p $) a one-dimensional wave equation. These 
one-dimensional features entail many simplifications. Indeed, we recover a simple one dimensional mean field equation 
(comparable to the model studied in \cite{MR2646061}) with all the three dimensional geometry encoded in the sphere 
(with the angle $ \omega \in \mathbb S^2 $ serving as a parameter). 

\vskip 0.2mm

- to interpret $ {\rm D} $ as an oscillatory integral. By this way,  the geometry of propagation is driven by the phase $ \pm 
\vert X (s) - x \vert + s -r $ which reveals very well the {\it joint} properties of the transport equation (\ref{VlasoveqRVMD})
and of the wave equation inside (\ref{maxwelleqRVMD}).

\smallskip

\noindent It follows that:

\vskip 0.2mm

{\it a)} the commutating vector fields considerations of \cite{MR1877669} are bypassed; 

\vskip 0.2mm

{\it b)} the non-resonant smoothing property of \cite{MR2124491} reduces to the manipulation of a non-zero 
coefficient (vanishing to $ 0 $ when $ \vert \xi \vert $ tends to $ + \infty $) occurring in a non-stationary phase; 

\vskip 0.2mm

{\it c)} the division lemma of \cite{Bouchut} is replaced by a straightforward argument (there is no more three 
spatial directions but instead only one direction $ p $ with a parameter $ \omega \in \mathbb S^2 $).

\noindent On the other hand, the return through the inverse Radon transform to the original spatial configuration 
(in $ x $) involves oscillatory integrals which, thanks to Lemmas \ref{splitD1} and \ref{splitDhgdzajh}, simplify and 
lead to the explicit formula (\ref{cD}), which is easy to use. 

\noindent The content of (\ref{cD}) does not imply any derivative of $ {\rm f}  $, 
$ {\rm E} $ or $ {\rm B} $. It does neither require any pointwise estimate on $ {\rm f}  $, $ {\rm E} $ or $ {\rm B} $. Instead, 
it relies on integral computations with respect to the volume element $ ds \, d \omega \, d \xi $. Modulo the (singular) Jacobian 
(\ref{jacobiadit}) which is issued from a pushforward to the original phase space, this amounts to work with the Lebesgue 
measure $ dx \, d \xi $. The great advantage is that this Liouville measure $ dx \, d \xi $ is preserved by the flow (\ref{diffeoom}). 
Henceforth, we obtain a connection between the value of $ {\rm D}(t,y,\eta) $, regarding one particular characteristic 
issued from $ (y,\eta)$, and the computation of the total energy
 \begin{equation}\label{consertotenRVMgen} 
 \text{\small $ \pmb{\mathscr{E}} $} (t) := \int_{\RR^3} \int_{\RR^3} \langle \xi \rangle \ {\rm f}(t,x,\xi) \ dx \, d\xi + \frac12 
\int_{\RR^3} |{\rm E}(t,x)|^2 \ dx + \frac12 \int_{\RR^3} |{\rm B}(t,x)|^2 \ dx ,
\end{equation}
which is for smooth solutions a conserved quantity \cite{Rhein} while, for weak solutions, this issue has not yet been completely 
resolved (see \cite{MR4077461} for a discussion about Onsager's conjecture). The fact remains that, in the actual context of 
smooth solutions, we have
 \begin{equation}\label{consertotenRVM} 
\text{\small $ \pmb{\mathscr{E}} $} (t) = \text{\small $ \pmb{\mathscr{E}} $}_{\! 0} := \text{\small $ \pmb{\mathscr{E}} $} (0) \, , 
\qquad t \geq 0 . 
\end{equation}
Thus, the search for pointwise estimates on $ {\rm D} $ is driven by the invariant energy estimate (\ref{consertotenRVM}). 
In that sense, we can again assert that our approach is at the junction between results about weak solutions (which rely 
on the sole preservation of $ \text{\small $ \pmb{\mathscr{E}} $} $ and which do not bring into play derivatives) and results
yielding a kind of strong information (namely sup-norm estimates for adequate quantities) without resorting to the 
computation of derivatives.

\noindent Amongst other things, the formula (\ref{cD}) has the benefit of revealing and delineating the role of a series of weights 
$ {\rm W}_\star $ defined on $ \RR_+ \times \mathbb S^2_\omega \times \RR^3_\xi $. Depending on the respective positions of 
the directions $ \Xi(s) $, $ \omega $  and $ \xi $, these weights may be larger or smaller. This remark gives access to quantitative 
estimates which differentiate between the localizations of $ \bigl(\Xi(s),\omega,\xi \bigr) $. This process is precisely what motivates 
Section \ref{RadonFourier}.
 
 
\subsubsection{Content of Section \ref{RadonFourier}} Section \ref{RadonFourier} answers the question 1. To this end,  in 
Subsection \ref{Preparatory work}, we start by a preparatory work: we tie together the controls of $ {\rm D} $ and $ {\rm P} $ 
(Paragraph \ref{Momentum spread}); we express the pushforward of the measure $ ds \, d \omega \, d \xi $ as a density with 
respect to $ dx \, d \xi $ (Paragraph \ref{Comparison of measures}); and we introduce useful tools concerning the microlocal 
weights $ {\rm W}_\star $ (Paragraph \ref{Three useful tools}); this allows (in Subsection \ref{Radon microlocal analysis of 
the weight functions}) to evaluate the sizes of the $ {\rm W}_\star $. This leads (in Subsection \ref{Weighted integrals}) to the 
study and to the control of weighted integrals. In the end (in Subsection \ref{Proof of Propositiontheorempms}), this provides 
key inputs for proving Proposition \ref{theorempms} (which also investigates large time issues).


\subsubsection{Content of Section \ref{Application}} Section \ref{Application} is devoted to Theorem \ref{mainmaintheo}. In Subsection
\ref{Preliminary material}, we specify the underlying framework: we introduce the notion of well-prepared approximate solutions 
$ {\rm U}^\eps_a $ (Paragraph \ref{Well-prepared approximate solutions}); we write the equations for the perturbation $ U $ 
(Paragraph \ref{equations for the perturbation}); and we explain what is meant by compatible initial data $ f^\eps_0 $ (Paragraph 
\ref{compilitdazh}). In Subsection \ref{discussenergy}, we come back to the question of energy estimates. It turns out that 
(\ref{consertotenRVM}) is not very useful. Instead, we have to propagate the $ L^2 $-norm of $ U(t) $ at the level of the linearized 
equations along $ {\rm U}^\eps_a $. That is the only way to ensure that the smallness $ L^2 $-conditions (\ref{remainderupesa2}) and 
(\ref{absorptioforsm}) is passed on to $ U (t) $. In Subsection \ref{prove Theorem mainmain}, we prove Theorem \ref{mainmaintheo}. 
This is done by estimating separately the approximate momentum increment $ {\rm D}^\eps_a $, bilinear terms obtained by freezing
$ {\rm D}_n $ along $ {\rm U}^\eps_a $, and a full nonlinear contribution $ D_n $. In this process, it is essential to know that $ U (t) $
is of size $ \eps $ in $ L^2 $. Indeed, this emerges as an indispensable prerequisite to compensate the presence inside (\ref{Vlasoveq}) 
of the large factor $ \eps^{-1} $. Again, it is noteworthy that the preceding approaches would furnish, if any, much less information than in 
Theorem \ref{mainmaintheo}.


\section{The proof of a representation formula} \label{Repformulamomentum} In this section, we consider a solution $ ({\rm f},{\rm E},{\rm B}) $ 
to the Cauchy problem (\ref{VlasoveqRVMD})-(\ref{maxwelleqRVMD})-(\ref{stationarystate}), which is assumed to be smooth, 
which is compactly supported with respect to $ \xi $, and which is defined on $ [0,T[ $ for some $ T \in \RR_+^* \cup \{+\infty \} $. 
In what follows, the time $ t $ is chosen in $ [0,T[ $. 

\noindent The Vlasov equation (\ref{VlasoveqRVMD}) is linked to a dynamical system on the phase space 
$ \RR^3_x \times \RR^3_\xi $. Recall that $ {\rm F} $ is the Lorentz force given by (\ref{VlasoveqRVMD}), and consider the flow $ (X,\Xi) $ 
obtained by solving the ordinary differential equation:
\renewcommand\arraystretch{2}
\begin{equation} \label{diffeoom}
\left \lbrace \begin{array} {ll}
\displaystyle \frac{d X}{dt} (t,y,\eta) = \nu (\Xi) , \quad & X(0,y,\eta) = y , \\
\displaystyle \frac{d \Xi}{dt} (t,y,\eta) = {\rm F} (t,X,\Xi) , \quad & \Xi(0,y,\eta) = \eta .
\end{array} \right. 
\end{equation}
The solution $ {\rm f} $ may be recovered by integrating along the characteristics, in the sense that
\begin{equation}\label{characteristicsmeaning} 
{\rm f}(t,x,\xi) = {\rm f}_0 \bigl( X(-t,x,\xi) , \Xi(-t,x,\xi) \bigr) .
 \end{equation}

\begin{defi}[Momentum increment] \label{momentum spread} The momentum increment at the time $ t $ associated with 
an initial phase point $ (y,\eta) $ is the difference $ {\rm D} (t,y,\eta)  := \langle \Xi (t,y,\eta) \rangle - \langle \eta \rangle $.
\end{defi}

\noindent In the absence of an electric field, that is when $ {\rm E} \equiv 0 $ as it is the case concerning the stationary 
solution $ \tilde {\rm U}^\eps_a $ of (\ref{prototype}), the kinetic energy $ \vert \Xi \vert^2/2 $ is just constant, and the same 
applies to $ \langle \Xi \rangle $ so that $ {\rm D} \equiv 0 $. But in general, we have $ {\rm D} \not \equiv 0 $ for 
two main reasons:
\begin{itemize}
\item [-] Impact of a non-zero initial data  $ ({\rm E}_0,{\rm B}_0) \not \equiv 0 $, with $ \nabla_x \times {\rm B}_0 \not \equiv 0 $
when $ {\rm E}_0 \equiv 0 $. Then, a non-trivial electric field $ {\rm E} \not \equiv 0 $ persists or is created (at least for small times 
$ t $). It can be approximated by solving the homogeneous version of Maxwell's equations: 
\begin{equation} \label{homogeneous version}
\part^2_{tt} {\rm E}_h - \Delta_x {\rm E}_h = 0 \, , \qquad {{\rm E}_h}_{\mid t=0} = {\rm E}_0 \, , \qquad \part_t {{\rm E}_h}_{\mid t=0} 
= \nabla_x \times {\rm B}_0 ,
\end{equation}
where the subscript $ h $ is for \underline{h}omogeneous. Retain that the access to $ {\rm E}_h $ is determined only by 
$ ({\rm E}_0 , {\rm B}_0 ) $; it is obtained by solving a linear wave equation; and it is completely decoupled from the Vlasov equation.
\item [-] Effect of nonlinear interactions (the quadratic terms in the Vlasov equation) together with the coupling (electric current in Maxwell's 
equations).
\end{itemize}

\noindent The momentum increment $ {\rm D} (t)$ can be computed from $ {\rm E}_h $, $ {\rm f} $ and $ {\rm F} $ as indicated below.

\begin{prop}[Representation formula for $ {\rm D} $] \label{theorem1} We have $ {\rm D} = {\rm D}_0 + {\rm D}_h + {\rm D}_l +  {\rm D}_n $ 
with 
 \begin{subequations}\label{cD} 
 \begin{align}
 & \displaystyle {\rm D}_0 :=  \int_0^t \int_{\mathbb S^2} \int_{\RR^3} {\rm W}_0 (s,\omega,\xi) \ {\rm f}_0 \bigl( X (s) + 
 s \omega,\xi \bigr) \ ds \, d \omega \, d \xi , \label{cD0} \\
 & \displaystyle {\rm D}_h := \int_0^t \nu \circ \Xi (s) \cdot {\rm E}_h \bigl( s,X(s)\bigr) \ ds , \label{cDl} \\
& \displaystyle {\rm D}_l := \int_0^t \Bigl( \int_r^t \int_{\mathbb S^2} \int_{\RR^3} {\rm W}_l (s,\omega,\xi) \  
{\rm f} \bigl(r,X(s) + (s-r) \omega,\xi \bigr) \ ds \, d \omega \, d \xi \Bigr) \, dr , \label{cDi1} \\
& \displaystyle \! \! \! \begin{array}{rl} 
\displaystyle {\rm D}_n := \int_0^t \Bigl( \int_r^t \int_{\mathbb S^2} \int_{\RR^3} \! \! \! & {\rm W}_n (r,s,\omega,\xi) \cdot 
{\rm F} \bigl(r,X(s) + (s-r) \omega,\xi \bigr) \\
& \times \, {\rm f} \bigl(r,X(s) + (s-r) \omega,\xi \bigr) \ ds \, d \omega \, d \xi \Bigr) \, dr ,
 \end{array} \label{cDi2} 
 \end{align}
 \end{subequations}  
 
 \noindent where the microlocal weights $ {\rm W}_\star $ are given by
  \begin{subequations}\label{the weights } 
 \begin{align}
 & \displaystyle {\rm W}_0 (s,\omega,\xi) := - \frac{s}{4 \pi} \ \frac{\nu \circ \Xi (s) \cdot  \bigl(\nu(\xi) + \omega \bigr)}{1
 + \omega \cdot \nu(\xi)} , \label{weights0} \\
& \displaystyle {\rm W}_l (s,\omega,\xi) := - \frac{1}{4 \pi} \ \frac{1}{\langle \xi \rangle^2} \ \frac{\nu \circ \Xi (s) \cdot 
\bigl( \omega + \nu(\xi) \bigr)}{\bigl( 1+ \omega \cdot \nu(\xi) \bigr)^2}  , \label{weightsi} \\
& \displaystyle {\rm W}_n (r,s,\omega,\xi) := - \frac{s-r}{4 \pi} \ \nabla_\xi \Bigl( \frac{\nu \circ \Xi (s) \cdot  \bigl(\nu(\xi) + 
\omega \bigr)}{1 + \omega \cdot \nu(\xi)} \Bigr) . \label{weightsn} 
 \end{align}
 \end{subequations}  
 \end{prop}
\renewcommand\arraystretch{1} 

 \noindent In (\ref{cD}), the subscripts $ 0 $, $ h $, $ l $ and $ n $ stand respectively for $ t = \underline{0} $ (initial time), 
 \underline{h}omogeneous, \underline{l}inear with respect to $ {\rm f} $, and \underline{n}onlinear in terms of $ {\rm f} $ 
 and $ {\rm F} $. We can separate inside $ {\rm D}_n $ the bilinear interactions involving $ {\rm E} $ and $ {\rm B} $. We 
 have $ {\rm D}_n = {\rm D}_{ne} + {\rm D}_{nb} $. The expression $ {\rm D}_{ne} $ and $ {\rm D}_{nb} $ are defined as 
 $ {\rm D}_n $ with $ ({\rm W}_n , {\rm F} ) $ replaced respectively by $ ({\rm W}_{ne} , {\rm E} )  $ and $ ({\rm W}_{nb} , 
 {\rm B} ) $ where 
 \begin{equation} \label{cDndecomponbinftyssb}
\qquad {\rm W}_{ne} (r,s,\omega,\xi) := {\rm W}_n (r,s,\omega,\xi) \, , \qquad  {\rm W}_{nb} (r,s,\omega,\xi) := {\rm W}_n 
(r,s,\omega,\xi) \times \nu (\xi) .
\end{equation}

\noindent In the pioneering works of Glassey and Strauss \cite{MR816621,MR1039487},  the electromagnetic field 
$ ({\rm E},{\rm B}) $ was represented in terms of $ {\rm E}_h$, $ {\rm f} $ and $ {\rm F} $. The pointwise estimates on $ {\rm E} $ 
and $ {\rm B} $ thus obtained were exploited (in \cite{MR816621,MR1039487} and subsequent works) to extract information 
on $ {\rm D} $. By contrast, we focus here directly on $ {\rm D} $. As a consequence, we will be able to control $ {\rm D} $ 
without resorting to (costly) sup-norm estimates on $ {\rm E} $ and $ {\rm B} $ but only through integrals involving $ {\rm f} $, 
$ {\rm E} $ and $ {\rm B} $. Seen in this light, our approach is more in line with \cite{Bouchut,MR1877669}, and especially 
\cite{Pallard}. It can be interpreted as an alternative to these contributions. Still, it differs from those both in its formulation, 
conception (some arguments we use to prove Proposition \ref{theorem1} are original), and consequences.
  
 \noindent This section is devoted to the decomposition of $ {\rm D} $ into (\ref{cD}). In Subsection \ref{prerequisites}, 
 we introduce basic tools. In Subsection \ref{Twooscillatoryintegrals}, we compute two oscillatory integrals. Then, in 
 Subsection \ref{proofoftheorem1}, we show Proposition \ref{theorem1}.

\noindent For the sake of simplicity, we will prove Proposition \ref{theorem1} for smooth solutions which are compactly 
supported with respect to both variables $ x $ and $ \xi $. The finite speed of propagation in $ x $ allows ultimately to 
relax this condition on the spatial support. Note also that Proposition \ref{theorem1} should remain true for less regular 
(weak) solutions under adequate integrability conditions in $ \xi $. But this aspect will not be investigated here.


\subsection{Prerequisites} \label{prerequisites} In Paragraph \ref{Lienard-Wiechertpotentials}, we start by adopting the 
intrinsic viewpoint of \cite{Bouchut,MR2124491,Pallard} to express $ {\rm D} $ in terms of the microscopic electromagnetic potential. 
In Paragraph \ref{basictools}, we recall basic facts about the Radon transform \cite{MR1723736}. 


\subsubsection{Lienard-Wiechert potentials} \label{Lienard-Wiechertpotentials} Choose a vector field $ {\rm A}_i : \RR^3_x 
\rightarrow \RR^3  $, where $ i $ stands for \underline{i}nitial, such that $ \nabla_x \cdot {\rm A}_i = 0 $ and $ \nabla_x \times 
{\rm A}_i = {\rm B}_0 $. Then, solve the wave equation 
\begin{equation} \label{homogeneous versionforpot}
\part^2_{tt} {\rm A}_h - \Delta_x {\rm A}_h = 0 \, , \qquad {{\rm A}_h}_{\mid t=0} = {\rm A}_i \, , \qquad \part_t {{\rm A}_h}_{\mid t=0} 
= - {\rm E}_0 .
\end{equation}
This allows to recover $ {\rm E}_h  $ through the relation $ {\rm E}_h = - \part_t {\rm A}_h $. As first noted in \cite{MR2124491} 
and exploited for instance in \cite{Bouchut,MR4084146}, the RVM system can be recast as a coupling between a wave equation 
and a Vlasov equation. To this end, it suffices to define the microscopic electromagnetic potential $ {\rm u} (t,x,\xi) $ which solves 
the initial value problem:
\begin{equation}\label{waveeqonu} 
(\partial^2_{tt} - \Delta_x) {\rm u} = {\rm f} , \qquad {\rm u}_{\mid t=0} = 0 , \qquad \partial_t {\rm u}_{\mid t=0} = 0.
\end{equation}
Then, the electromagnetic field $ ({\rm E},{\rm B}) $ can be computed from $u$ by the two identities
\begin{subequations}\label{defdeEetB} 
\begin{align}
& \displaystyle {\rm E}(t,x) = - \part_t {{\rm A}_h} (t,x) - \int_{\RR^3} \bigl \lbrack \nu(\xi) \partial_t {\rm u} + \nabla_x {\rm u} \bigr 
\rbrack (t,x,\xi) \, d \xi , \label{defdeEetBpourE} \\
& \displaystyle {\rm B}(t,x) = \nabla_x \times {\rm A}_h (t,x) + \int_{\RR^3} \nabla_x \times \bigl \lbrack {\rm u} \, \nu(\xi) \bigr 
\rbrack (t,x,\xi) \, d \xi . \label{defdeEetBpourB} 
 \end{align}
 \end{subequations}  
System \eqref{VlasoveqRVMD}-\eqref{waveeqonu}-\eqref{defdeEetB} is self-contained. It is equivalent to 
\eqref{VlasoveqRVMD}-\eqref{maxwelleqRVMD}. This is why, in what follows, it will also be referred to as the "RVM system". 
From (\ref{diffeoom}) and due to the special structure inside (\ref{VlasoveqRVMD}) of $ {\rm F} $, we have 
\begin{equation}\label{deefdeQFF}
{\rm D}(t,y,\eta) = \int_0^t \frac{d}{ds} \langle \Xi(s) \rangle \ ds = \int_0^t \nu \circ \Xi (s) \cdot {\rm E} \bigl( s,X(s)\bigr) \ ds .
\end{equation}
We can plug (\ref{defdeEetBpourE}) into (\ref{deefdeQFF}). The part $ - \part_t {{\rm A}_h} $ inside (\ref{defdeEetBpourE}) 
leads to $ {\rm D}_h $. We have $ {\rm D} = {\rm D}_h + \breve {\rm D} $ with
\begin{equation}\label{deefdeQ}
\qquad \breve {\rm D}(t,y,\eta) = - \int_0^t \int_{\RR^3} \Bigl \lbrace \nu \circ \Xi (s) \cdot \nu (\xi) \, \partial_s {\rm u} \bigl( s,X(s) , 
\xi \bigr) + \nu \circ \Xi (s) \cdot \nabla_x {\rm u} \bigl( s,X (s) , \xi \bigr) \Bigr \rbrace \, ds \, d\xi .
\end{equation}
This explains the origin of the term $ {\rm D}_h $ inside (\ref{cD}). Observe that the influence of $ {\rm A}_h $ is not limited to 
$ {\rm D}_h $. It does also impact $ \breve {\rm D} $ (through $ {\rm u} $). Indeed, equations \eqref{VlasoveqRVMD}
and \eqref{waveeqonu} are coupled with $ {\rm A}_h $ appearing inside \eqref{VlasoveqRVMD} because the Lorentz 
force $ {\rm F} = {\rm E} + \nu (\xi) \times {\rm B} $ must be computed with $ {\rm E} $ and $ {\rm B} $ given by (\ref{defdeEetB}). 

\noindent On the one hand, the formula (\ref{deefdeQFF}) seems to indicate that the control of $ {\rm D} $ should require 
a sup-norm estimate on $ {\rm E} $. On the other hand, the second identity (\ref{deefdeQ}) suggests that the access to 
$ \breve {\rm D} $ should imply sup-norm estimates on $ \partial_s {\rm u} $ and $ \nabla_x {\rm u} $. As is well-known
\cite{Bouchut,MR2124491,MR1877669,Pallard}, this gives false impressions. As will be seen, more can be done.


\subsubsection{Reminders on the Radon transform} \label{basictools} The Radon transform \cite{MR1723736} is the map defined by
\[ \begin{array}{rcl} 
R : C^0_c ( \RR^3 ; \CC ) & \longrightarrow & C^0_c ( \mathbb S^2 \times \RR ; \CC ) \\
{\rm f} (x) & \longmapsto & R {\rm f} (\omega ,p) ,
\end{array} \]
where 
\begin{equation}\label{defradontransform}
{\rm g} (\omega,p) := R {\rm f} (\omega,p) := \int_{ H_{\omega,p} } {\rm f}(x) \,  dm(x) = {\rm g} (- \omega,-p) , 
\end{equation}
and $ dm $ is the Euclidean measure on the hyperplane $ H_{\omega,p} := \{ x \in \RR^3 \, ; \, \omega \cdot x = p \} $.
Recall that $ R {\rm f} $ is linked to the Fourier transform $ \hat {\rm f} $ through the relation
\[ \hat {\rm f} (\rho \, \omega) = \int_{-\infty}^{+\infty} e^{-i \rho p } \, R {\rm f}(\omega, p) \, dp, \qquad \rho \in \RR_+ . \]
It follows that  
\begin{equation}\label{laplaceradon} 
R (\partial_{x_i} {\rm f}) (\omega,p) = \omega_i \ \part_{p} {\rm g} (\omega,p) \, , \qquad 
R (\Delta {\rm f}) (\omega,p) = \part^2_{pp} {\rm g} (\omega,p) . 
\end{equation}
As explained for instance in line (38)-p. 19 of \cite{MR1723736}, with $ \eta = \rho \, \omega \in \RR^3 $, we 
obtain that
\begin{equation}\label{inverseformr} 
\begin{array} {rl}
{\rm f} (x) \! \! \! & \displaystyle = \frac{1}{(2 \pi)^3} \int_{\RR^3} e^{ix\cdot \eta} \, \hat {\rm f} (\eta) \, d \eta = \frac{1}{(2 \pi)^3} 
\int_{\mathbb S^2} \! \int_{0}^{+\infty} e^{i \rho (\omega\cdot x)} \, \hat {\rm f} (\rho \, \omega) \, \rho^2 \, d \rho \, d \omega \\
& \displaystyle = \frac{1}{2} \, \frac{1}{(2 \pi)^3} \int_{\mathbb S^2} \Bigl \lbrace \int_{-\infty}^{+\infty} \rho^2 \, e^{i \rho (\omega\cdot x)} 
\Bigl \lbrack \int_{-\infty}^{+\infty} e^{-i \rho p} \, R {\rm f}(\omega,p) \, dp \Bigr \rbrack \, d \rho \Bigr \rbrace \, d \omega ,
\end{array}  
 \end{equation}
where the change of variables $ (\omega,\rho) $ into $ (-\omega,-\rho) $ has been exploited to pass from the first to the second line.
This yields the (three dimensional) inversion formula
\begin{equation}\label{inverseformrexpl} 
{\rm f} (x) = R^{-1} {\rm g} (x) := - \frac{1}{2} \, \frac{1}{(2 \pi)^2} \int_{\mathbb S^2} \partial^2_{pp} {\rm g} (\omega, \omega\cdot x ) \, d \omega .
 \end{equation}
Also observe that 
\begin{equation}\label{linkrhat} 
\forall n \, \in \NN, \qquad (i \, \rho)^n \ \hat {\rm f} ( \rho \omega) = \cF_p \bigl( (\partial^n_p R {\rm f}) (\omega,\cdot) \bigr) (\rho) . 
\end{equation}
As a consequence, we have 
\begin{equation}\label{retourr} 
\forall n \, \in \NN, \qquad (\partial^n_p R {\rm f}) (\omega,p) = \frac{1}{2 \pi} \, \int_\RR \int_{\RR^3} e^{i \rho (p-\omega \cdot x)} \, 
(i \, \rho)^n \,  {\rm f}(x) \, d \rho \, dx .
\end{equation}
The relation (\ref{laplaceradon}) indicates that a derivative in $ x $ is on the Radon side a derivative in $ p $. But 
the Radon transform is clearly associated with a smoothing effect (due to the integration). And thereby, its inverse can consume 
derivatives of $ {\rm g} $, as demonstrated on (\ref{inverseformrexpl}). This means that gains of derivatives (in $ p $) are hidden 
behind the integration with respect to $ \omega $. This phenomenon is detected in Subsection \ref{Twooscillatoryintegrals} and exploited 
in Subsection \ref{proofoftheorem1}.


\subsection{Two oscillatory integrals} \label{Twooscillatoryintegrals} As in the case of the Fourier transform, the study of partial 
differential equations with the Radon transform may reveal the role of oscillatory integrals. In the present context, these integrals 
are of two types: the first (Paragraph \ref{Twooscillatoryintegrals1}) involves a compact surface (the sphere $ \mathbb S^2 $); the 
second (Paragraph \ref{Twooscillatoryintegrals2}) implies the whole space $ \RR \times \RR^3 $. In both cases, the aim is to take 
advantage of cancellation properties induced by oscillations. 

\subsubsection{Oscillatory integrals on the sphere} \label{Twooscillatoryintegrals1} In view of (\ref{linkrhat}), a derivative $ \partial_p $ 
on the Radon side costs a multiplication by $ \rho $. On the other hand, as can be inferred from (\ref{retourr}), this loss may be 
associated with oscillatory integrals implying phases looking like $ \rho \, (\omega \cdot {\rm X} + \tau) $. It is compensated at the 
level of (\ref{inverseformrexpl}) by an integration with respect to $ \omega $. Thus, there is some underlying spatial averaging effect 
(on $ \mathbb S^2 $). This is highlighted below.

\begin{lem} [Gain of one derivative] \label{splitD1} Fix $ {\rm X} \in \RR^3 \setminus \{ 0 \} $ and $ (\rho,\tau) \in \RR^2 $. We have
\begin{equation} \label{hefkzalqq}
\int_{\mathbb S^2}  e^{i \rho \, (\omega \cdot {\rm X} + \tau) } \ (i \rho) \ d \omega = \frac{2 \pi}{\vert {\rm X} \vert} \, \sum_\pm \pm 
e^{i \rho \, ( \pm \vert {\rm X} \vert + \tau )} = 4 \pi \, i \, e^{i \rho \, \tau} \, \frac{\sin ( \rho \, \vert {\rm X} \vert)}{\vert {\rm X} \vert} .
 \end{equation}
\end{lem}

\begin{prof} Let $ \cR $ be a rotation such that $ \cR \, ( {\rm X} / \vert {\rm X} \vert ) = {}^t (0,0,1) $. We can change $ \omega $ into 
$ \upomega := \cR \omega $ and then work in spherical coordinates, that is with
\begin{equation} \label{sphericalcoordinates}
\quad \upomega = \left( \begin{array}{c}
\cos \varphi \, \sin \vartheta \\
\sin \varphi \, \sin \vartheta \\
\cos \vartheta 
\end{array} \right) , \qquad \varphi \in [0,2\pi] , \qquad \vartheta \in [0,\pi] , \qquad d \upomega = \sin \vartheta \, d \vartheta
\, d \varphi . 
\end{equation}
Since $ \cR^{-1} \upomega \cdot {\rm X} = \upomega \cdot \cR {\rm X} $, this furnishes 
\[ \begin{array}{rl}
\displaystyle \int_{\mathbb S^2} \! \! \!  & \displaystyle e^{i \rho \, (\omega \cdot {\rm X} + \tau)} \ (i \rho) \ d \omega 
= \displaystyle \int_{\mathbb S^2}  e^{i \rho \, ( \cR^{-1} \upomega \cdot {\rm X} + \tau)} \ (i \rho) \ d \upomega =
 \int_{\mathbb S^2}  e^{i \rho \, \bigl( \vert X \vert \, \upomega \cdot \cR ({\rm X} / \vert {\rm X} \vert)  + \tau \bigr)} \ (i \rho) \ 
 d \upomega \\
\! \! \!  & \displaystyle = \int_0^\pi \int_0^{2\pi} e^{i \rho \, (\vert {\rm X} \vert \cos \vartheta + \tau)} \, (i \rho) \ \sin 
\vartheta \, d \vartheta \, d \varphi  = - \frac{2 \pi}{\vert {\rm X} \vert } \int_0^\pi \partial_\vartheta \Bigl \lbrace 
e^{i \rho \, [\vert {\rm X} \vert \cos \vartheta +\tau)} \Bigr \rbrace \ d \vartheta \\
\! \! \!  & \displaystyle = \displaystyle \frac{2 \pi}{\vert {\rm X} \vert } \, \Bigl \lbrace e^{i \rho [ + \vert {\rm X} \vert + \tau]} 
- e^{i \rho [ - \vert {\rm X} \vert + \tau]} \Bigr \rbrace ,
 \end{array} \]
 which is exactly (\ref{hefkzalqq}). \end{prof} 
 
 \noindent Let $ \nu \in \RR^3 $. Applying the differential operator $ \nu \cdot \nabla_{\rm X} $ to the identity (\ref{hefkzalqq}),
 we end up with
 \begin{equation} \label{end up}
\int_{\mathbb S^2} \nu \cdot \omega \ e^{i \rho \, (\omega \cdot {\rm X} + \tau) } \ \rho^2 \ d \omega = 2 \pi \ \frac{\nu \cdot 
{\rm X}}{\vert {\rm X} \vert^3} \, \sum_\pm ( \pm 1 - i \, \rho \, \vert {\rm X} \vert) \ e^{i \rho \, (\pm \vert {\rm X} \vert + \tau )} .
 \end{equation}
 
 
 \subsubsection{Oscillatory integrals on the whole space} \label{Twooscillatoryintegrals2} Replace $ {\rm X} $ by $ {\rm X} - x $. 
 Looking at the right hand side of (\ref{hefkzalqq}), we see that the integration with respect to $ d \omega $ may produce singular weights (like $ 1/ 
 \vert {\rm X} - x \vert $ near $ {\rm X} = x $) in factor of oscillations. Now, the integration with respect to $ d \rho \, dx $ of such 
 expressions multiplied by $ {\rm f} $ can produce the integral of the trace of $ {\rm f} $ on spheres.

\begin{lem} [Passage from singular weights to traces] \label{splitDhgdzajh} Let $ {\rm f} : \RR^3 \rightarrow \RR $ be 
a compactly supported function of class $ C^1 $. Let $ {\rm K} : \mathbb S^2 \rightarrow \RR $ be a bounded function.
Fix $ {\rm X} \in \RR^3 $, $ \alpha \in \RR $ and $ \tau \in \RR $. Then, for $ \alpha \leq 2 $, we have
\begin{equation} \label{hlv	za}
\begin{array}{rl} 
\displaystyle \int_\RR \int_{\RR^3} \! \! \! & \displaystyle  e^{i \rho [ \pm \vert {\rm X} -x \vert+ \tau ]} \ {\rm K} \Bigl( \frac{
x-{\rm X}}{\vert x - {\rm X} \vert} \Bigr) \ \frac{{\rm f}(x)}{\vert {\rm X} - x \vert^\alpha} \ d \rho \, dx \\
\ & \displaystyle = 2 \pi \int_{\mathbb S^2} \vert \tau \vert^{2-\alpha} \ H(\mp \tau) \ {\rm K} (\omega) \ {\rm f} (X + \vert 
\tau \vert \, \omega) \ d \omega ,
\end{array} 
 \end{equation}
where $ H $ is the Heaviside function (in the half-maximum convention which is important at least when $ \alpha = 2 $ 
and $ \tau = 0 $), namely
\[ H(\tau) := \frac{1}{2} \bigl \lbrack \mathbb 1_{\RR_+} (\tau-) + \mathbb 1_{\RR_+} (\tau+) \bigr \rbrack = \left \lbrace 
\begin{array}{lcl} 
0 \ & \text{ if } \ & \tau < 0 , \\
1/2 \ & \text{ if } \ & \tau = 0 , \\
1 \ & \text{ if } \ & \tau > 0 .\\
\end{array} \right. \]
\end{lem}

\begin{prof} In spherical coordinates $ r \, \omega $ for $ x - {\rm X} $, we have to deal with
\[ \begin{array}{l} 
\displaystyle \int_{\mathbb S^2} \int_\RR {\rm K} (\omega) \, \Bigl( \int_{\RR_+} r^{2 - \alpha} \ {\rm f} ({\rm X} + r \, \omega) \ 
e^{- i \rho \, ( \mp r - \tau)} \, dr \Bigr) \, d \rho \, d \omega \\
\displaystyle \qquad  = \int_{\mathbb S^2}  {\rm K} (\omega) \, \Bigl \lbrace \int_{\RR} \Bigl( \int_{\RR} \psi_\mp 
({\rm r}) \, e^{- i \rho \, {\rm r}}  \ d {\rm r} \Bigr) \, d \rho \Bigr \rbrace \, d \omega = \int_{\mathbb S^2} {\rm K} (\omega) \, 
\Bigl \lbrace \int_{\RR} \cF (\psi_\mp) (\rho) \, d \rho \Bigr \rbrace \, d \omega \, ,
\end{array} \]
where we have changed $ r $ into $ {\rm r}  := \mp r - \tau $, and where we have introduced the function $ \psi_\mp : 
\RR \rightarrow \RR $ (depending on $ \omega $) given by
\[ \psi_\mp ( {\rm r}) := \vert {\rm r} + \tau \vert^{2 - \alpha} \ \mathbb 1_{\RR_+} \bigl( \mp ({\rm r} + \tau) \bigr) \ {\rm f} 
\bigl( {\rm X} + \vert {\rm r} + \tau \vert \, \omega) . \]
The function $ \psi_\mp $ is compactly supported. It is bounded and piecewise $ C^1 $ as long as $\alpha\leq 2$ with only one 
possible discontinuity when $ \alpha= 2 $ (located at $ {\rm r} = - \tau $). Thus, we can apply the Dirichlet condition 
for inversion of Fourier integrals which furnishes
\[ \int_{\RR} \cF (\psi_\mp) (\rho) \, d \rho = 2 \pi \ ( \cF^{-1} \circ \cF) (\psi_\mp) (0) = \pi \, \bigl( \psi_\mp(0+) 
+ \psi_\mp(0-) \bigr) . \]
After substitution, we find (\ref{hlv za}). \end{prof} 


\subsection{Proof of Proposition \ref{theorem1}} \label{proofoftheorem1} The demonstration is done in three stages.
In Paragraph \ref{TransitiontotheRadonside}, we express $ \breve {\rm D} $ in terms of $ \part^2_{pp} {\rm g} $. In Paragraph 
\ref{Returnwithoscillatoryintegrals}, we exploit (\ref{retourr}) to exhibit adequate oscillatory integrals. Then, in Paragraph 
\ref{Outcoming}, we draw the conclusions.


\subsubsection{The Radon side picture} \label{TransitiontotheRadonside} Introduce
\[ {\rm v} (t,\omega,p,\xi) := R \bigl \lbrack {\rm u} (t,\cdot,\xi) \bigr \rbrack ( \omega,p), \qquad {\rm g} (t,\omega,p,\xi) := R \bigl \lbrack 
{\rm f} (t,\cdot,\xi) \bigr \rbrack ( \omega,p) . \]
Under the action of the Radon transform, the three-dimensional wave equation is transformed into a 
one--dimensional wave equation. With $ {\rm v}_1 \equiv 0 $, the Cauchy problem (\ref{waveeqonu}) becomes
\begin{equation}\label{waveeqonura} 
(\partial^2_{tt} - \partial_{pp}^2) {\rm v} = {\rm g} , \qquad {\rm v}_{\mid t=0} = 0 , \qquad \partial_t {\rm v}_{\mid t=0} = {\rm v}_1 .
\end{equation}
On the other hand, from the Vlasov equation, we can deduce that
\begin{equation}\label{vlasoveqonura} 
 \partial_t {\rm g} + \nu (\xi) \cdot \omega \ \partial_p {\rm g} + \text {\rm div}_\xi \, R ( {\rm f} \, {\rm F})  = 0 , \qquad {\rm g}_{\mid t=0} 
 = {\rm g}_0 := R \bigl \lbrack {\rm f}_0 (\cdot,\xi) \bigr \rbrack . 
 \end{equation}
At the level of (\ref{deefdeQ}), the integral $ \breve {\rm D} $ is built with four types of derivative: $ \partial_s {\rm u} $ and $ \part_{x_i} {\rm u} $
with $ i \in \{ 1,2,3 \} $. In Subparagraph a), we show that  $ \partial_s {\rm u} $ can be expressed as a function of $ \partial^2_{pp} 
{\rm g} $. In Subaragraph b), we do the same for $ \nabla_x {\rm u} $. By this way, in Subparagraph c), we can extract a reformulation 
of $ \breve {\rm D} $ in terms of one type of derivative, namely $ \partial^2_{pp} {\rm g} $. By this way, contrary to the division lemma 
\cite{Bouchut,MR1877669}, looking at $ \breve  {\rm D} $ as depending on $ {\rm g} $ instead of $ {\rm u} $ allows to reduce the kinds 
of derivatives which are needed: we can concentrate on $ \part_p $ only.
 
 \smallskip 
 
\noindent $ \bullet $ a) {\it Computation of $ \partial_s {\rm u} $.}  
The solution to the one-dimensional wave equation (\ref{waveeqonura}) can be obtained through Duhamel's formula
\begin{equation}\label{Duhamelforv}  
{\rm v} (t,\omega,p,\xi) = \frac{1}{2} \int_{p-t}^{p+t} {\rm v}_1 (\omega,y,\xi) \, dy + \frac{1}{2} \int \! \! \int_{\triangle(t,p)}  
{\rm g}(s,\omega,y,\xi) \, ds \, dy , 
\end{equation}
where $ \triangle(t,p) $ is the triangle
\[ \triangle(t,p) := \bigl \lbrace \, (s,y) \, ; \, s \in [0,t] \ , \ p-(t-s) \leq y \leq p+(t-s) \, \bigr \rbrace . \]
From (\ref{waveeqonura}) with $ {\rm v}_1 \equiv 0 $, we find also
\[ (\partial^2_{tt} - \partial_{pp}^2) (\partial_t {\rm v}) = \partial_t {\rm g} , \qquad \partial_t {\rm v}_{\mid t=0} = 0 , \]
as well as 
\[ \partial_t (\partial_t {\rm v})_{\mid t=0} = {\rm g}_{\mid t=0} + (\partial_{pp}^2 {\rm v})_{\mid t=0} = {\rm g}_0 . \]
Applying (\ref{Duhamelforv}) to the above equation, we get
\begin{equation}\label{forforv}  
\begin{array}{rl}
\partial_t {\rm v}(t,\omega,p,\xi) \! \! \! & \displaystyle =  \frac{1}{2} \int_{p-t}^{p+t}  {\rm g}_0 (y,\xi) \, dy + 
\frac{1}{2} \int \! \! \int_{\triangle(t,p)}  \partial_s {\rm g} (s,\omega,y,\xi) \, ds \, dy \\
\ & \displaystyle = \frac{1}{2} \int_{p-t}^{p+t}  {\rm g}_0 (y,\xi) \, dy + \frac{1}{2} \Bigl \lbrace \int_{p-t}^p 
\Bigl( \int_0^{y-p+t} \partial_s {\rm g}(s,\omega,y,\xi) \, ds \Bigr) \, dy \\
\ & \displaystyle \quad  +  \int^{p+t}_p \Bigl( \int_0^{-y+p+t} \partial_s {\rm g}(s,\omega,y,\xi) \, ds \Bigr) \, dy \Bigr \rbrace \\
 \ & \displaystyle = \frac{1}{2} \int_0^t \bigl \lbrack {\rm g} (s,\omega,p+t-s,\xi) + {\rm g} (s,\omega,p-t+s,\xi) \bigr \rbrack  \, ds .
\end{array} 
\end{equation}
From (\ref{forforv}), we can infer that
\begin{equation} \label{waveissued}
\begin{array}{rl}
\displaystyle \partial_s {\rm u} ( s,x , \xi ) \! \! \! & \displaystyle = R^{-1} (\partial_s {\rm v}) ( s,x, \xi ) = - \frac{1}{2} \, \frac{1}{(2 \pi)^2} 
\int_{\mathbb S^2} (\partial^2_{pp} \partial_s {\rm v}) (s,\omega, \omega \cdot x ,\xi ) \, d \omega \\
\ & \displaystyle = - \frac{1}{4} \, \frac{1}{(2 \pi)^2} \, \sum_{\pm} \int_{\mathbb S^2}  \Bigl \lbrace \int_0^s \partial^2_{pp} {\rm g} 
\bigl(r,\omega, \omega \cdot x \pm (s-r) ,\xi \bigr) \, dr \Bigr \rbrace \, d \omega .
 \end{array} 
 \end{equation}
 Since $ \partial^2_{pp} {\rm g} (r,\cdot,\xi) $ is in view of (\ref{defradontransform}) an even function, we know that
 \begin{equation} \label{parite}
\partial^2_{pp} {\rm g} \bigl(r,-\omega, - \omega \cdot x - (s-r) ,\xi \bigr) =  \partial^2_{pp} {\rm g} \bigl(r,\omega, \omega \cdot x 
+ (s-r) ,\xi \bigr) , 
 \end{equation}
 and therefore, changing $ \omega $ into $ - \omega $, we can deduce that
  \begin{equation} \label{fh	eihkaz}
 \partial_s {\rm u} ( s,x , \xi ) = - \frac{1}{2} \, \frac{1}{(2 \pi)^2} \, \int_{\mathbb S^2}  \Bigl \lbrace \int_0^s \partial^2_{pp} {\rm g} 
 \bigl(r,\omega, \omega \cdot x + (s-r) ,\xi \bigr) \, dr \Bigr \rbrace \, d \omega . 
  \end{equation}
  
\noindent $ \bullet $ b) {\it Computation of $ \nabla_x {\rm u} $.} From (\ref{waveeqonura}) with $ {\rm v}_1 \equiv 0 $, we can extract
\begin{equation}\label{waveeqonurappp} 
(\partial^2_{tt} - \partial_{pp}^2) (\partial_p {\rm v}) = \partial_p {\rm g} , \qquad (\partial_p {\rm v})_{\mid t=0} = 0 , \qquad \partial_t 
(\partial_p {\rm v})_{\mid t=0} = 0 .
\end{equation}
Applying (\ref{Duhamelforv}), this yields 
\begin{equation}\label{forforpartpv}  
\begin{array}{rl}
\partial_p {\rm v}(t,\omega,p,\xi) \! \! \! & \displaystyle = \frac{1}{2} \int \! \! \int_{\triangle(t,p)}  \partial_p {\rm g} (s,\omega,y,\xi) \, ds \, dy \\
 \ & \displaystyle = \frac{1}{2} \int_0^t \bigl \lbrack {\rm g} (s,\omega,p+t-s,\xi) - {\rm g} (s,\omega,p-t+s,\xi) \bigr \rbrack  \, ds .
\end{array} 
\end{equation}
On the other hand, from (\ref{laplaceradon}), we have
 \[ R \bigl \lbrack  \nu \circ \Xi (s) \cdot \nabla_x {\rm u} \bigr \rbrack ( s,\omega, p, \xi ) = \nu \circ \Xi (s) \cdot  \omega \ \partial_p 
{\rm v} ( s,\omega, p, \xi ) , \] 
 and therefore, using (\ref{parite}) again and (\ref{forforpartpv}), we obtain that
 \begin{equation} \label{waveissuedvid}
\begin{array}{l}
\displaystyle \nu \circ \Xi (s) \cdot \nabla_x {\rm u} (s,x,\xi) = R^{-1} \bigl \lbrack \nu \circ \Xi (s) \cdot  \omega \ \partial_p {\rm v} \bigr 
\rbrack (s,x,\xi) \\
\qquad \displaystyle = - \frac{1}{2} \, \frac{1}{(2 \pi)^2} \int_{\mathbb S^2} \nu \circ \Xi (s) \cdot \omega \ \partial^2_{pp} (\partial_p {\rm v}) 
(s,\omega, \omega \cdot x ,\xi ) \, d \omega \\
\qquad \displaystyle = - \frac{1}{2} \, \frac{1}{(2 \pi)^2} \, \int_{\mathbb S^2} \Bigl \lbrace \int_0^s \nu \circ \Xi (s) \cdot \omega \ 
\partial^2_{pp} {\rm g} (r,\omega, \omega \cdot x + s-r ,\xi ) \, dr \Bigr \rbrace \, d \omega .
 \end{array} 
 \end{equation}
 
 \smallskip
 
\noindent $ \bullet $ c) {\it Summary.} It suffices to plug (\ref{fh	eihkaz}) and (\ref{waveissuedvid}) inside (\ref{deefdeQ}) 
to get a representation formula for $ \breve {\rm D} $ on the Radon side
\[ \breve  {\rm D} = \frac{1}{8 \pi^2} \int_{\RR^3} \int_{\mathbb S^2}  \int_0^t \Bigl \lbrace \int_0^s \nu \circ \Xi (s) 
\cdot \bigl( \omega + \nu (\xi) \bigr) \ \partial^2_{pp} {\rm g} \bigl(r,\omega, \omega \cdot X (s) + s-r ,
\xi \bigr) \, dr \Bigr \rbrace \, ds \, d \omega \, d\xi . \] 
Recall (\ref{laplaceradon}) which says that $ \partial^2_{pp} = R \Delta $. Thus, the computation of $ \breve {\rm D} $ seems
to consume two derivatives of $ {\rm f} $. The aim of the next paragraph is to show that this is not the case. The goal 
is to remove the presence of $ \partial^2_{pp} $ inside $ {\rm D} $.


\subsubsection{Analysis through oscillatory integrals} \label{Returnwithoscillatoryintegrals}
To better understand the content of $ \breve {\rm D} $, we can apply (\ref{retourr}) with $ n = 2 $ to exhibit the following 
oscillatory integral
\renewcommand\arraystretch{1.4} 
\[ \begin{array}{rl}
\displaystyle \breve {\rm D} = - \frac{1}{16 \pi^3} \int_0^t \! \int_0^s \! \int_\RR \! \int_{\RR^3} \! \int_{\RR^3} \Bigl( \int_{\mathbb S^2} 
\nu \circ \Xi (s) \cdot \bigl( \omega + \nu (\xi) \bigr) \! \! \! & \displaystyle  e^{i \rho [ \omega \cdot X (s) - \omega \cdot x + s-r]} \, 
\rho^2 \ d \omega \Bigr) \\
\ & \displaystyle \times \, {\rm f} (r,x,\xi) \ ds \, dr \, d \rho \, dx \, d \xi . 
\end{array} \qquad \quad \]
\renewcommand\arraystretch{1} 

\noindent We can apply Lemma \ref{splitD1} and (\ref{end up}) with $ {\rm X} = X(s) -x $ and $ \tau = s-r $ to get 
\[ \begin{array}{rl}
\displaystyle \breve {\rm D} = \! \! \! & \displaystyle \sum_\pm  \mp \, \frac{1}{8 \pi^2} \int_0^t \! \int_0^s \! \int_\RR \! \int_{\RR^3} \! 
\int_{\RR^3} \frac{\nu \circ \Xi (s) \cdot \bigl( X(s)-x \bigr)}{\vert X (s) -x \vert^3} \ {\rm f} (r,x,\xi) \ e^{i \rho \, [ \pm \vert X (s) 
- x \vert + s -r]} \ ds \, dr \, d \rho \, dx \, d \xi \\
+ \! \! \! & \displaystyle \sum_\pm  \pm \, \frac{1}{8 \pi^2} \int_0^t \! \int_0^s \! \int_\RR \! \int_{\RR^3} \! \int_{\RR^3} 
\Bigl( \frac{\nu \circ \Xi (s) \cdot \nu(\xi)}{\vert X (s) -x \vert} \pm \frac{\nu \circ \Xi (s) \cdot \bigl( X(s)-x \bigr)}{\vert X (s) -
x \vert^2} \Bigr) \\ 
\ & \qquad \qquad \qquad \qquad \qquad \quad \ \times \, {\rm f} (r,x,\xi) \ (i \rho) \ e^{i \rho \, [ \pm \vert X (s) - x \vert + s -r]} \ 
ds \, dr \, d \rho \, dx \, d \xi .
\end{array}  \]
In the second sum, there is still $ \rho $ in factor (which corresponds to the lost of one derivative). The next idea is to eliminate 
this weight $ \rho $ by exploiting the underlying presence of oscillations. This is the principle of non-stationary phase. To this end, 
the strategy is to perform an integration by parts with respect to $ r $ and $ x $. As usual, this operation costs time and spatial 
derivatives of the symbol $ {\rm f} $. It must be done without introducing unmanageable derivatives of $ {\rm f} $. In practice, 
taking into account the Vlasov equation, we can convert the derivative $ \part_r + \nu(\xi) \cdot \nabla_x $ into derivatives with
respect to $ \xi $. With this in mind, we look at
\begin{equation} \label{lookatderivative}
\begin{array}{rl}
\displaystyle \bigl( \part_r + \nu(\xi) \cdot \nabla_x \bigr) \! \! \! & \displaystyle \Bigl \lbrace e^{i \rho \, [ \pm \vert X (s) - x \vert 
+ s -r]} \Bigr \rbrace \\
 & \displaystyle = \Bigl( -1 \pm \frac{x - X (s)}{\vert x - X (s) \vert} \cdot \nu (\xi) \Bigr) \ (i \rho) \ e^{i \rho \, [ \pm \vert X (s) - x 
 \vert + s -r]} .
 \end{array}
\end{equation}
Since $ \vert \nu (\xi) \vert < 1 $, the multiplicative factor in the right hand side is negative, and therefore it can be inverted. 
This argument exploits a microlocal ellipticity property. As mentioned in {\it b)} of Paragraph \ref{An alternative path}, it is 
referred to as a {\it non-resonant smoothing} property \cite{MR2124491}.

\noindent The identity (\ref{lookatderivative}) connects $ \part_r + \nu(\xi) \cdot \nabla_x $ to the multiplication by $ \rho $. In view of 
(\ref{linkrhat}), this amounts to apply the derivative $ \partial_p $. Now, recall that $ R (\partial_{x_i} {\rm f}) = \omega_i \, 
\partial_p (R {\rm f}) $, and thereby $ \omega_i \, \partial_p $ may be viewed on the Radon side as a condensed version 
of the spatial derivatives $ \partial_{x_i} $. Thus,  a link is established between $ \part_r + \nu(\xi) \cdot \nabla_x $ and
$ \partial_{x_i} $. 

\noindent Historically \cite{Bouchut,MR1877669}, this was done by converting the derivative $ \partial_{x_i} {\rm u} $ 
inside (\ref{deefdeQ}) into $ \part_r {\rm f} + \nu(\xi) \cdot \nabla_x {\rm f} $. Here, there is no need of such division lemma,
alinea {\it c)} of Paragraph \ref{An alternative path}. We just observe that the light cones (related to the wave equation for the 
electromagnetic fields), which may be viewed as the level surfaces in $ \RR \times \RR^3 $ of the phase functions $ \pm \vert 
X (s) - x \vert + s -r $, are transversal to the derivative $ \part_r + \nu(\xi) \cdot \nabla_x $ (which is related to the time-spatial 
transport in the Vlasov equation). With the help of (\ref{lookatderivative}), in the second sum defining $ \breve {\rm D} $, we 
can interpret the coefficient which is in factor of $ {\rm f} $ according to
\[ \begin{array}{rl}
\displaystyle \Bigl( \frac{\nu \circ \Xi (s) \cdot \nu(\xi)}{\vert X (s) -x \vert} \! \! \! & \displaystyle \pm \frac{\nu \circ \Xi (s) \cdot 
\bigl( X(s)-x \bigr)}{\vert X (s) - x \vert^2} \Bigr) \ (i \rho) \ e^{i \rho \, [ \pm \vert X (s) - x \vert + s -r]} \\
= \! \! \! & \displaystyle - k_\mp (s,x,\xi) \ \bigl( \part_r + \nu(\xi) \cdot \nabla_x \bigr) \Bigl \lbrace e^{i \rho \, [ \pm \vert X (s) - x \vert 
+ s -r]} \Bigr \rbrace ,
\end{array} \]
with
\[ k_\mp (s,x,\xi) := \frac{1}{\vert X (s) -x \vert} \ K_\mp \Bigl( s, \frac{x-X(s)} {\vert x -X (s) \vert} ,\xi \Bigr) \, , \qquad 
K_\mp(s,\omega,\xi) := \frac{\nu \circ \Xi (s) \cdot \bigl(\nu(\xi) \mp \omega \bigr)}{1\mp \omega \cdot \nu(\xi)} .\] 
Let $ T_k $ be the multiplicative operator by the function $ k $, that is 
\begin{equation} \label{multiplicativeoperator}
T_k({\rm f}) (r,s,x,\xi) := k(s,x,\xi) \, {\rm f} (r,x,\xi) .
\end{equation}
Now, we can integrate by parts with respect to $ r $ and $ x $ to obtain 
\[ \begin{array}{rl}
\displaystyle \breve {\rm D} = \! \! \! & \displaystyle \sum_\pm  \mp \, \frac{1}{8 \pi^2} \int_0^t \! \int_0^s \! \int_\RR \! \int_{\RR^3} \! 
\int_{\RR^3} \frac{\nu \circ \Xi (s) \cdot \bigl( X(s)-x \bigr)}{\vert X (s) -x \vert^3} \ {\rm f} (r,x,\xi) \ e^{i \rho \, [ \pm \vert X (s) 
- x \vert + s -r]} \ ds \, dr \, d \rho \, dx \, d \xi \\
+ \! \! \! & \displaystyle \sum_\pm  \pm \, \frac{1}{8 \pi^2} \int_0^t \int_\RR \int_{\RR^3} \Bigl \lbrace \int_0^s \int_{\RR^3} \bigl( 
\part_r + \nu(\xi) \cdot \nabla_x \bigr) \bigl( T_{k_\mp} ({\rm f}) \bigr) \ e^{i \rho \, [ \pm \vert X (s) - x \vert + s -r]} \ dx \, dr \Bigr 
\rbrace \ d\xi \, d \rho \, ds \\
- \! \! \! & \displaystyle \sum_\pm  \pm \, \frac{1}{8 \pi^2} \int_0^t \int_\RR \int_{\RR^3} \Bigl \lbrack \int_{\RR^3} T_{k_\mp} 
({\rm f}) \ e^{i \rho \, [ \pm \vert X (s) - x \vert + s -r]} \ dx \Bigr \rbrack_0^s \ d\xi \, d \rho \, ds .
\end{array}  \]
This yields two boundary terms ($ {\rm D}_{b0} $ and $ {\rm D}_{bs} $ issued respectively from $ r = 0 $  and $ r = s $),
a term $ {\rm D}_l $ including the integrands which have $ {\rm f} $ in factor, as well as a contribution $ {\rm D}_n $ which
implies derivatives of $ {\rm f} $. More precisely, we have $ \breve {\rm D} = {\rm D}_{b0} + {\rm D}_{bs} + {\rm D}_l + {\rm D}_n $ 
with
\[ {\rm D}_{b0} = \sum_\pm \pm {\rm D}_{b0}^\pm \, , \qquad {\rm D}_{bs} = \sum_\pm \mp {\rm D}_{bs}^\pm \, , 
\qquad {\rm D}_l = \sum_\pm \pm {\rm D}_l^{\pm} \, , \qquad {\rm D}_n = \sum_\pm \pm {\rm D}_n^{\pm} \, , \]
where
\renewcommand\arraystretch{1.4} 
\begin{align*}
& \displaystyle {\rm D}_{b0}^\pm := \frac{1}{8 \pi^2} \int_0^t \int_\RR \int_{\RR^3}  \int_{\RR^3} T_{k_\mp} ({\rm f}) 
(0,s,x,\xi) \ e^{i \rho \, [ \pm \vert X (s) - x \vert + s ]} \ ds \, d \rho \, dx \, d \xi , \\
& \displaystyle {\rm D}_{bs}^\pm := \frac{1}{8 \pi^2} \int_0^t \int_\RR \int_{\RR^3}  \int_{\RR^3} T_{k_\mp}({\rm f}) 
(s,s,x,\xi) \ e^{\pm  i \rho \, \vert X (s) - x \vert} \ ds \, d \rho \, dx \, d \xi  , \\
& \displaystyle {\rm D}_l^{\pm} := \frac{1}{8 \pi^2}  \int_0^t \int_0^s \int_\RR \int_{\RR^3}  \int_{\RR^3}T_{\tilde k_\mp} 
({\rm f}) \ e^{i \rho \, [ \pm \vert X (s) - x \vert + s -r]} \ ds \, dr \, d \rho \, dx \, d \xi , \\
& \displaystyle {\rm D}_n^{\pm} :=  \frac{1}{8 \pi^2} \int_0^t \int_0^s \int_\RR \int_{\RR^3}  \int_{\RR^3} T_{k_\mp} \bigl(
\part_r {\rm f}  + \nu(\xi) \cdot \nabla_x {\rm f} \bigr) \ e^{i \rho \, [ \pm \vert X (s) - x \vert + s -r]} \ ds \, dr \, d \rho \, 
dx \, d \xi .
\end{align*} 
 
 \noindent By construction, we find 
 \[ \tilde k_\mp := - \frac{\nu \circ \Xi (s) \cdot \bigl( X(s)-x \bigr)}{\vert X (s) -x \vert^3} + \nu(\xi) \cdot \nabla_x k_\mp 
 = \frac{1}{\vert X (s) -x \vert^2} \ \tilde K_\mp \Bigl( s, \frac{x-X(s)} {\vert x -X (s) \vert} ,\xi \Bigr)  , \]
 with
 \[ \tilde K_\mp ( s, \omega ,\xi ) := \frac{1}{\langle \xi \rangle^2} \ \frac{\nu \circ \Xi (s)  \cdot \bigl( \omega \mp \nu(\xi) 
 \bigr)}{\bigl( 1\mp \omega \cdot \nu(\xi) \bigr)^2} . \] 
Now, the time spatial derivative $ \part_r {\rm f} + \nu(\xi) \cdot \nabla_x {\rm f} $ can be exchanged with velocity 
derivatives. Indeed, using the Vlasov equation, it can be converted into derivatives with respect to $ \xi $ (which are harmless 
because the coefficients are smooth in $ \xi $ and because the phase does not depend on $ \xi $). We find that
\[ {\rm D}_n^{\pm} =  \frac{1}{8 \pi^2}  \int_0^t \int_0^s \int_\RR \int_{\RR^3}  \int_{\RR^3} \nabla_\xi k_\mp \cdot {\rm F} \ 
{\rm f} \ e^{i \rho \, [ \pm \vert X (s) - x \vert + s -r]} \ ds \, dr \, d \rho \, dx \, d \xi . \]
On the other hand, in $ {\rm D}^\pm_l $ and $ {\rm D}^\pm_n $, we can switch the order of integrations according to
\[ \int_0^t \Bigl( \int_0^s \cdots \, dr \Bigr) \, ds = \int_0^t \Bigl( \int_r^t \cdots \, ds \Bigr) \, dr . \] 


\subsubsection{Epilog} \label{Outcoming} In general, the Fourier analysis as well as the Radon analysis of nonlinear 
partial differential equations lead to complex formulas. And indeed, the above oscillatory integrals seem complicated.
Surprisingly, they can be significantly simplified by applying Lemma \ref{splitDhgdzajh}. As a matter of fact, all the 
computations can be made explicit:
 \begin{itemize}
\item [-] Study of $ {\rm D}_{b0} $. We take $ {\rm K} = K_\mp $, $ {\rm X} = X (s) $, $ \alpha = 1 $ and $ \tau = s \geq 0 $.
We find that $ {\rm D}_{b0}^+ = 0  $ and $ {\rm D}_{b0} = - {\rm D}_{b0}^-  $. With $ {\rm W}_0 = - s \, K_+ / (4 \pi) $, we 
can recognize $ {\rm D}_0 $ as in (\ref{cD0}).
\item [-] Study of $ {\rm D}_{bs} $. We take $ {\rm K} = K_\mp $, $ {\rm X} = X (s) $, $ \alpha = 1 $ and $ \tau = 0 $ to see that 
$ {\rm D}_{bs}^\pm = 0 $.
\item [-] Study of $ {\rm D}_l $. We take $ {\rm K} = \tilde K_\mp $,  $ {\rm X} = X (s) $, $ \alpha = 2 $ and $ \tau = s - r \geq 0 $ 
to obtain $ \mathcal D_l^{+} = 0 $ so that $ {\rm D}_l = - {\rm D}_l^{-} $. This is coherent with (\ref{cDi1})
where $ {\rm W}_l = - \tilde K_+ / (4 \pi) $ as in (\ref{weightsi}).
\item [-] Study of $ {\rm D}_n $. We deal with a vector valued version of (\ref{hlv za}) where $ {\rm f} $ is replaced by 
 $ {\rm F} \, {\rm f} $. We take $ {\rm K} = \nabla_\xi K_\mp $,  
$ {\rm X} = X (s) $, $ \alpha = 1 $ and $ \tau = s - r \geq 0 $ to see that $ \mathcal D_n^+ = 0 $ and therefore that $ {\rm D}_n 
= - {\rm D}_n^{-} $. By this way, with $ {\rm W}_n = - (s-r) \, \nabla_\xi K_+ / (4 \pi) $, we find that $ {\rm D}_n $ is given by (\ref{cDi2}).
\end{itemize}

\noindent The proof of Proposition \ref{theorem1}  is now complete. 


\section{The control of the momentum spread} \label{RadonFourier} 
In this section, we consider a solution $ {\rm U} $ to the Cauchy problem (\ref{VlasoveqRVMD})-(\ref{maxwelleqRVMD})-(\ref{stationarystate}), 
assumed to be smooth, compactly supported (in $ \xi $), and defined on $ [0,T[ $ where $ T \equiv T({\rm U}_0) \in \RR_+^* \cup \{+\infty \} $ 
is the maximum life-span of this smooth solution. From (\ref{characteristicsmeaning}), we know that the support of $ {\rm f} (t,\cdot) $ 
is the image by the map $ (X,\Xi) (t,\cdot) $ of the support of $ {\rm f}_0 $. Given $ t \in [0,T[ $, we can define the 
maximal size $ {\rm Q} (t) $ of the spatial support of $ {\rm f} (t,\cdot) $, which is 
\begin{equation}\label{maximalsizess} 
\qquad {\rm Q} (t) := \inf \ \bigl \lbrace R \in \RR_+ \, ; \, {\rm f} (t,x,\xi) = 0 \ \text{for all} \  \xi \in \RR^3 \ \text{and for all} \ x 
\in \RR^3 \ \text{with} \ R \leq \vert x \vert \bigr \rbrace  ,
\end{equation}
as well as the maximal size $ {\rm P} (t) $ of the momentum support of $ {\rm f} (t,\cdot) $, which is
\begin{equation}\label{maximalsizesm} 
\qquad {\rm P} (t) := \inf \ \bigl \lbrace R \in \RR_+ \, ; \, {\rm f}(t,x,\xi) = 0 \ \text{for all} \  x \in \RR^3 \ \text{and for all} \ \xi 
\in \RR^3 \ \text{with} \ R \leq \vert \xi \vert \bigr \rbrace .
\end{equation}
By construction, we have 
\[ \text{supp} \ {\rm f} \subset \bigl \lbrace (t,x,\xi) \in [0,T[ \times \RR^3 \times \RR^3 \, ; \, \vert x \vert \leq {\rm Q} (t) \, , \, 
\vert \xi \vert \leq {\rm P} (t) \, \bigr \rbrace . \]
The quantity $ {\rm Q} $ yields a control on the size of the {\it spatial domain of influence}. In view of the first equation of (\ref{diffeoom}), 
the spatial speed of propagation is bounded by one, so that
\begin{equation}\label{controlofcq} 
{\rm Q} (t) \leq {\rm Q}_0 + t \, , \qquad {\rm Q}_0 := {\rm Q} (0) . 
\end{equation}
Both Vlasov and Maxwell's equations have a finite spatial  speed of propagation (bounded by $ 1 $). Thus, exploiting the notion 
of region of influence, to prove the local smooth solvability, it suffices to work with solutions that are compactly supported with 
respect to the space variable $ x $. To simplify, we can localize the spatial and momentum support in the same ball (say of 
size $ {\rm P}_0 $). With this in mind, we replace the condition inside (\ref{mathfrak N}) on $ \text{\rm supp} \ {\rm f}_0 $ by
\begin{equation}\label{repetcompact space} 
\text{\rm supp} \ {\rm f}_0 \subset B(0,{\rm P}_0]  \times B(0,{\rm P}_0] .
\end{equation}
On the other hand, the quantity $ {\rm P} $ gives a bound on the extent of the {\it momentum domain of influence}. It provides insight 
into the {\it momentum spread}. Without a control involving the sup-norm of $ {\rm F} $ (that is equivalently of $ {\rm E} $ and 
$ {\rm B}) $, the second equation of (\ref{diffeoom}) does not provide with a bound for the momentum speed of propagation. This is 
usually resolved by looking at Lipschitz bounds on $ {\rm U} = ({\rm f} , {\rm E} , {\rm B} ) $. From the pioneering works 
\cite{MR969207,MR816621,MR1039487,zbMATH03951266} on smooth solutions, we know that there exists a time $ \cT_s \in \RR_+^* $ 
and a continuous function $ \cF_s : [0,\cT_s] \rightarrow \RR_+ $ (the subscript $ s $ stands for \underline{s}mooth) depending both 
on the Lipschitz norm of $ {\rm U}_0 $ such that 
\begin{equation}\label{conofcps} 
{\rm P} (t) \leq \cF_s (t) \, , \qquad \forall \, t \in [0,\cT_s] \, , \qquad 0 < \cT_s \leq T . 
 \end{equation}
A key issue is whether (\ref{conofcps}) remains true under less restrictive criteria on the initial data $ {\rm U}_0 $. Our purpose 
here is to remove the Lipschitz condition on $ {\rm f}_0 $ and to relax the cost of two derivatives concerning $  ({\rm E}_0 , {\rm B}_0 ) $.

\begin{prop}[Control on the momentum spread by mild information] \label{theorempms} Fix $ {\rm P}_0 \in \RR_+ $, and 
consider the corresponding subspace $ \mathfrak N $ adjusted as in (\ref{mathfrak N}). Given any $  {\rm S}_0 \in \RR_+^* $,
select initial data $ {\rm U}_0 \in \mathfrak N $  satisfying $ \cN ({\rm U}_0) \leq {\rm S}_0 $ with $ \cN $ as in (\ref{mathcal N}). 
Then, there exist a time $ \cT \in \RR_+^* $ and a continuous increasing function $ \cF : [0,\cT] \rightarrow \RR_+ $, both 
depending only on $ {\rm S}_0 $ so that $ \cT \equiv \cT({\rm S}_0) $ and $ \cF \equiv \cF({\rm S}_0;\cdot) $, such that
\begin{equation}\label{conofcpbis} 
\qquad {\rm P} (t) \leq \cF ( {\rm S}_0 ; t ) , \qquad \forall \, t \in \bigl[ 0, \min \bigl( T({\rm U}_0) ; \cT({\rm S}_0) \bigr) \bigr[ . 
\end{equation}
Moreover, under (\ref{repetcompact space}), the behavior of $ \cT $ for small values of $ {\rm S}_0 $ is bounded below 
according to
\begin{equation}\label{behaviorcTsmallvalues} 
\exists \, c \in \RR_+^* ; \qquad c \ {\rm S}_0^{-1} \leq \cT({\rm S}_0) \, , \qquad \forall \, {\rm S}_0 \in \RR_+^* . 
\end{equation}
\end{prop}

\noindent From Glassey-Strauss continuation criterion \cite{MR816621}, assuming that $ T ({\rm U}_0) \in \RR_+^* $, 
the size of $ {\rm P} (t) $ must explode when $ t \rightarrow T ({\rm U}_0)- $. Looking at (\ref{conofcpbis}), since $ \cF $ 
is continuous on $ \bigl \lbrack 0, \cT({\rm S}_0) \bigr \rbrack $ and therefore bounded, there is a contradiction if 
$ T({\rm U}_0) < \cT ({\rm S}_0) $. Expressed in terms of  (\ref{lower bound fo T}), this means that
\begin{equation}\label{lowerboundforlifespan} 
0 < \cT({\rm S}_0) \leq  T(\mathfrak N, \cN;{\rm S}_0) \leq T ({\rm U}_0) . 
\end{equation}
Theorem \ref{inimaintheo} is a direct consequence of (\ref{lowerboundforlifespan}). Given $ {\rm P}_0$, the line 
(\ref{behaviorcTsmallvalues}) specifies how fast $ T(\mathfrak N, \cN;{\rm S}_0) $ tends to $ + \infty $ when 
$ {\rm S}_0 $ goes  to zero. Without (\ref{repetcompact space}), a version of (\ref{behaviorcTsmallvalues}) is 
still available under adaptations. To this end, the spatial support of $ {\rm f}_0 $ must be truncated on a larger ball 
of size $ C \, t $ with $ C > 1 $, and the smallness parameters must be revisited.

\noindent Now, let $ {\rm G} : \RR^3 \rightarrow \RR^q $ with $ q \in \NN^* $. The usual fluid description of plasmas 
(MHD) involves macroscopic quantities like
\[ {\rm M}_{\rm G} (t,x) := \int_{\RR^3} \, {\rm G} (\xi) \ {\rm f} (t,x,\xi) \ d\xi . \]
For $ {\rm G} \equiv 1 $, we deal with the number density. For $ {\rm G} \equiv \nu(\xi) $, we recover the current density. 
For $ {\rm G}_n (\xi) = \xi \otimes \xi \otimes \cdots \otimes \xi $, where $ \xi $ is multiplied $ n $ times with $ n \in \NN^* $, 
we find the n-th moment, with in particular the momentum density (for $ n=1 $).

\begin{cor}[Sup-norm controls on all fluid quantities under mild information] \label{corrempms} In the context of 
Proposition \ref{theorempms}, for all $ {\rm G} \in L^\infty (\RR^3;\RR^q) $, we have
\begin{equation}\label{Sup-normcontrolquantities} 
\qquad \forall \, (t,x) \in [0,\cT({\rm S}_0)] \times \RR^3 \, , \qquad \parallel {\rm M}_{\rm G} (t,x) \parallel \leq 4 \pi \, {\rm S}_0 
\sup_{\vert \xi \vert \leq \cF({\rm S}_0;t)} \, \parallel {\rm G} (\xi) \parallel^3 < + \infty .
\end{equation}
\end{cor}

\noindent This furnishes a range of {\it a priori sup-norm estimates} which, contrary to (\ref{conofcps}), do not require 
any regularity on $ {\rm f}_0 $, and which implement (relatively) weak estimates on $ ({\rm E}_0,{\rm B}_0) $. This may 
seem surprising in the quasilinear context (\ref{VlasoveqRVMD})-(\ref{maxwelleqRVMD}) under study. Recall however 
that such bounds are basically inherited from the transport part (the Vlasov equation) after its (complicated) interaction 
with Maxwell's equations. 

\noindent Our construction is also a gateway to a notion of solutions which is at the interface between the strong and weak versions 
of respectively \cite{MR816621,MR1039487} and \cite{MR1003433}.

\begin{cor}[Strong-weak solutions] \label{strong-weak} Fix $ ({\rm P}_0 , {\rm S}_0) \in \RR_+ \times \RR_+^* $. Select 
$ {\rm f}_0 \in L^\infty (\RR^3 \times \RR^3) $ satisfying (\ref{cestpourf0s}), (\ref{cestpourf0infty}) and (\ref{compatibility conditionsrmf0}).
Choose $ {\rm E}_0 \in H^1(\RR^3) $ and $ {\rm B}_0 \in H^1(\RR^3) $. Assume that $ \cN ({\rm U}_0) \leq {\rm S}_0 $
where $ \cN $ is as in (\ref{mathcal N}) with $ \bar p = 2 $. Then, we can find a time $ \cT \equiv \cT({\rm S}_0) \in \RR_+^* $ 
depending only on $ {\rm S}_0 $ such that the Cauchy problem (\ref{VlasoveqRVMD})-(\ref{maxwelleqRVMD})-(\ref{stationarystate}) 
has a unique solution satisfying 
\begin{equation}\label{swuniquesol} 
\text{\rm supp} \ {\rm f} (t,\cdot) \subset \RR^3 \times B \bigl(0, \cF ( {\rm S}_0 ; t ) \bigr)] \, , \qquad \forall t \in [0,\cT] 
\end{equation}
as well as $ {\rm E} \in H^1([0,\cT] \times \RR^3) $ and $ {\rm B} \in H^1([0,\cT] \times \RR^3) $.
\end{cor}

\begin{prof} Any bounded function $ {\rm f}_0 $ satisfying (\ref{cestpourf0s}), (\ref{cestpourf0infty}) and (\ref{compatibility conditionsrmf0}) 
 can be approximated by a sequence $ ( {\rm f}^n_0)_n $ subject to (\ref{cestpourf0})-(\ref{compatibility conditionsrmf0}) 
 uniformly with respect to $ n \in \NN $. Similarly, $ {\rm E}_0 $ and $ {\rm B}_0 $ can be approximated by  
 $ ({\rm E}_0^n)_n $ and $ ({\rm B}_0^n)_n $ with $ {\rm E}_0^n \in C^2_c $ and $ {\rm B}_0^n \in C^2_c $. Moreover, it 
 can be ensured that $ {\rm U}_0^n := ({\rm f}^n_0,{\rm E}_0^n,{\rm B}_0^n) $ is such that $ \cN ({\rm U}_0^n) \leq {\rm S}_0 $ 
 for all $ n \in \NN $. Theorem \ref{inimaintheo} gives access to solutions $ {\rm U}^n $ on $ [0,\cT] $ associated with the initial 
 data $ {\rm U}_0^n $, and satisfying (\ref{conofcpbis}). By compactness  arguments 
 (based on averaging lemmas \cite{MR1003433,Rhein}), passing to the limit ($ n \rightarrow + \infty $), we can extract a 
 corresponding weak solution $ {\rm U} = ({\rm f},{\rm E},{\rm B}) \in L^\infty \times L^2 \times L^2 $ which still satisfies 
 (\ref{conofcpbis}). But from \cite{MR2124491}-Theorem 1 (or alternatively \cite{MR3900881}) together with (\ref{conofcpbis}), 
 we can deduce\footnote{This argument has been reported to us by Nicolas Besse. Observe that the information (\ref{conofcpbis}) 
 is crucial to recover the $ H^1 $-regularity. It is missing 
 in the case of the weak solutions provided by DiPerna-Lions \cite{MR1003433}} that $ \xi $-averages of $ u $ (without 
 the need of momentum cutoff) are in $ H^2 ([0,\cT] \times \RR^3) $, and therefore that  $ {\rm E} \in H^1 $ and $ {\rm B} 
 \in H^1 $. From there, applying \cite{MR1022305}, we can recover the uniqueness.
\end{prof} 

\noindent The information (\ref{conofcpbis}) is also adapted to {\it concrete applications} and to further stability results 
which, as in \cite{Bouchut}, may be inherited from (\ref{conofcpbis}). Since the study of strongly magnetized plasmas were our point of entry \cite{MR4084146,MR4357273}, 
this axis of research is prioritized in the present article. For the moment, we can only say that $ T $ is all the greater given that 
$ {\rm S}_0 $ is small. And we can already guess that $ T $ could remain fixed for a large magnetic field whose size is 
adequately compensated by the smallness of $ {\rm S}_0 $. This will be confirmed in Section \ref{Application}. 

\noindent In Subsection \ref{Preparatory work}, we lay the background for a phase space analysis. These preliminaries lead 
in Subsection \ref{Weighted integrals} to the study of weighted integrals. This results in Subsection \ref{Proof of Propositiontheorempms} 
in the proof of Proposition \ref{theorempms}. From now on, we will work implicitly with $ t < T $.


\subsection{Preparatory work} \label{Preparatory work}  In Paragraph \ref{Momentum spread}, we make the connection 
between the size $ {\rm P} $ of the momentum spread and quantities $ {\rm D}^{\rm a}_\star $ (derived from the 
$ {\rm D}_\star $). In Paragraph \ref{Comparison of measures}, we take the pushforward of the measure $ ds \, d \omega \, d \xi $ 
to recover (modulo a Jacobian) the Liouville measure $ dx \, d\xi $. In Paragraph \ref{Three useful tools}, we show three lemmas that are helpful 
in Paragraph \ref{Radon microlocal analysis of the weight functions} to estimate the weight functions $ \vert {\rm W}_\star \vert $. 


\subsubsection{Control of the momentum spread through the representation formula} \label{Momentum spread} 
By construction, for all $ (s,x) \in [0,t] \times \RR^3 $, the momentum support of $ {\rm f} (s,x,\cdot) $ is contained 
in the ball of radius
\begin{equation}\label{ball of radius} 
 \langle {\rm P} \rangle_\infty (t) := \sup \ \bigl \lbrace \bigl( 1 + {\rm P}(s)^2 \bigr)^{1/2} \, ; \, s \in [0,t] \bigr \rbrace .
 \end{equation}
By compactness arguments using the continuity of $ {\rm f }$, we can find some $ (t_0,x_0,\xi_0) $ such that
\[ (t_0,x_0,\xi_0) \in [0,t] \times \text{supp} \ {\rm f} (t_0,\cdot)  \, , \qquad \langle \xi_0 \rangle = \langle {\rm P} 
\rangle (t_0) = \langle {\rm P} \rangle_\infty (t) . \]
In view of (\ref{characteristicsmeaning}), the position $ (x_0,\xi_0) $ is necessarily the image by $ (X,\Xi)(t_0,\cdot) $
of some $ (y_0,\eta_0) $ in the support of $ {\rm f}_0 $. In other words, we have $ \xi_0 = \Xi(t_0,y_0,\eta_0) $
for some $ \eta_0 $ satisfying $ \vert \eta_0 \vert \leq {\rm P}_0 $. It follows that
\begin{equation}\label{compactness arguments} 
\langle {\rm P} \rangle_\infty (t) = \langle \Xi (t_0,y_0,\eta_0) \rangle = \langle\eta_0 \rangle + {\rm D} (t_0,y_0,\eta_0) 
\leq 1 + {\rm P}_0 + \vert {\rm D} (t_0,y_0,\eta_0) \vert . 
 \end{equation}
By this way, the control of $ \langle {\rm P} \rangle_\infty (t) $ boils down to the study of $ \vert {\rm D} (t_0,y_0,\eta_0) \vert $. 
Moreover, for $ s \not = t_0 $, we find that 
$ (X,\Xi)(s,y_0,\eta_0) $ must be in the support of $ {\rm f} (s,\cdot)$ so that
 \begin{equation}\label{supplball of radius} 
\langle \Xi(s,y_0,\eta_0) \rangle \leq \langle {\rm P} \rangle_\infty (s) \leq \langle {\rm P} \rangle_\infty (t) \, , 
\qquad \forall s \in [0,t] .
 \end{equation}
Unless necessary, we will no more mention the selection of $ (y_0,\eta_0) $. For $ \star \in \{0,h,l,ne,nb \} $, we 
denote by $ {\rm D}^{\rm a}_\star $ (with superscript \og $ {\rm a} $ \fg{} for {\underline a}bsolute value) the expression 
$ {\rm D}_\star $ where $ {\rm f} $, $ {\rm W}_\star $, $ {\rm E} $ and $ {\rm B} $ are replaced by $ \vert {\rm f} \vert $ 
($ = {\rm f} $), $ \vert {\rm W}_\star \vert $, $ \vert {\rm E} \vert $ and $ \vert {\rm B} \vert $, and where the scalar products 
(like $ {\rm W}_n \cdot {\rm E} $) are just substituted for products of norms (like $ \vert {\rm W}_n \vert \, \vert {\rm E} \vert $). 
By definition, the $ {\rm D}^{\rm a}_\star $ are increasing functions of the time $ t $, and therefore $ {\rm D}_\star (t_0) \leq 
{\rm D}^{\rm a}_\star (t_0) \leq {\rm D}^{\rm a}_\star (t) $. From Proposition \ref{theorem1} together with 
(\ref{compactness arguments}), we can assert that
\begin{equation} \label{control of the momentum spread}  
0 \leq \langle {\rm P} \rangle_\infty (t) \leq 1 + {\rm P}_0 + {\rm D}^{\rm a}_0 (t) + {\rm D}^{\rm a}_h (t) + {\rm D}^{\rm a}_l (t) 
+ {\rm D}^{\rm a}_{ne} (t) + {\rm D}^{\rm a}_{nb} (t)  . 
\end{equation}
Retain also that, from (\ref{cestpourf0infty}) and (\ref{characteristicsmeaning}), we get 
\begin{equation} \label{estimateforrmf}  
0 \leq {\rm f} (r,x,\xi) \leq {\rm S}_0 \, , \qquad \forall \, (r,x,\xi) \in [0,T[ \times \RR^3  \times \RR^3 .
\end{equation}
Looking at (\ref{control of the momentum spread}), we can already assert that $ \langle {\rm P} \rangle_\infty $ is 
qualitatively controlled by zero-order information on $ {\rm U} $, with no derivatives of $ {\rm f} $, 
$ {\rm E}_h $, $ {\rm E} $ or $ {\rm B} $. In fact, the situation is even better since the influence of $ {\rm f} $, 
$ {\rm E} $ and $ {\rm B} $ is expressed through integrals with respect to $ d s \, d\omega \, d \xi $. This aspect 
is examined in the next paragraph.


\subsubsection{Comparison of measures, microlocal pictures and related difficulties} \label{Comparison of measures} 
All the integrals $ {\rm D}^{\rm a}_\star $ (except for $ \star = h $) involve the differential element $ d s \, d\omega $.
It is worth noting that this $ d s \, d\omega $ appears in our analysis after special procedures, by mixing different aspects:
\begin{itemize}
\item [-] the part $ d s $ comes from (\ref{deefdeQFF}). We do not look at the speed of propagation (in momentum)
which in our context may be large. Instead, we consider some integral version of it. When dealing with (\ref{diffeoom}), 
we exploit Duhamel's principle.
 \item [-] the part $ d \omega $ is issued from (\ref{inverseformrexpl}). The Radon transform involves  integrals on 
 families of (two-dimensional) hyperplanes. Its inverse implies  integrals on (two-dimensional) spheres with 
 respect to $ d \omega $.
\end{itemize}

\noindent By assembling $ ds $ and $ d \omega $, we find $ d s \, d\omega $. This combination of the variables $ s $
and $ \omega $  emerges also at the level of (\ref{cD}). Indeed, in the $ {\rm D}^{\rm a}_\star $, all expressions, like  
$ {\rm f} (r,\cdot,\xi) $, are evaluated at the specific position $ X(s) + (s-r) \omega $ (note that $ r = 0 $ in case of 
$ {\rm D}^{\rm a}_0 $). This furnishes a link to the Lebesgue measure $dx $. To this end, we can (for instance) 
exploit Lemma 2.2 of \cite{MR3291372} which is recalled below for the sake of completeness.

\begin{lem}\label{by C. Pallard} [by C. Pallard] Let $ r \in \RR $. The map 
\[ \begin{array}{rcl}
\cX \, : \, [r,t] \times {\mathbb S^2} & \longrightarrow & \RR^3 \\
(s,\omega) & \longmapsto & \cX (s,\omega) := X(s) + (s-r) \, \omega
\end{array} \]
is a $ C^1 $-diffeomorphism onto some region of $ \RR^3 $. Its Jacobian $ J $ is given by
\begin{equation} \label{jacobiadit}  
J(r,s,\omega) := (s-r)^2 \ \bigl( 1 + \omega \cdot \nu \circ \Xi(s) \bigr) \, , \qquad J \ d s \, d\omega = dx .
\end{equation}
\end{lem} 

\noindent Thus, the pushforward (by $ \cX $) of the measure $ d s \, d\omega $ is absolutely continuous with respect 
to the Lebesgue measure: it is just $ J^{-1} \, dx $. Now, we can extend this argument at a microlocal level. The Radon 
transform is an integral (hence non local) operator and, as such, it is not suitable for transferring localizations. However, 
what counts inside the $ {\rm D}^{\rm a}_\star $ is $ \cX $, which can effectively translates (microlocal) localizations in terms 
of $ (s,\omega, \xi) $ into (microlocal) localizations in terms of $ (x,\xi) $ - and conversely. In doing so, the landscape is 
changing. The advantage is that working with $ (s,\omega, \xi) $ instead of $ (x,\xi) $ is much easier.

\noindent This transfer of microlocal localizations is achieved by the map $ \cX $ which is built (through $ X $) on a 
complete nonlinear evolution (on a special solution to the RVM system) and on the choice of a specific characteristic 
(namely $ X $ for a special selection of $ y_0 $ and $ \eta_0 $). Passing to the phase space, we have to deal with the 
pushforward by $ \cX \otimes Id_\xi $ of the measure $ d s \, d\omega \, d \xi $. Modulo the Jacobian $ J $, this enables 
a connection with the Liouville measure $ dx \, d \xi $. 

\noindent One of the main interests of Proposition \ref{theorem1} associated with (\ref{jacobiadit}) is to produce explicit 
weights which indicate where $ {\rm f} $ , $ {\rm E} $ and $ {\rm B} $ contribute most to the $ {\rm D}^{\rm a}_\star $. 
Depending upon the relative positions of $ \Xi(s) $, $ \omega $ and $\xi $, these weights can be larger or smaller.
To perform estimates by taking care of the respective localizations of $ \Xi(s) $, $ \omega $ and $\xi $ is what we call 
here {\it Radon Fourier analysis}. We use the expression \og  {\it Radon} \fg{} because the angle $ \omega $ (resp. 
the measure $ d \omega $) has appeared after a Radon transform (resp. through its inverse $ R^{-1} $). 

\noindent In the end, we are faced with three types of singularities related to the content of the weight functions 
$ \vert {\rm W}_\star \vert $ or to negative powers of $ J $ thus introduced:
 \begin{itemize}
 \item [a)] Complication at $ \omega = - \xi / \vert \xi \vert $ coming from negative powers of $ 1 + \omega 
\cdot \nu (\xi) $ inside the $ \vert {\rm W}_\star \vert $. This leads to unbounded coefficients in the perspective
of a global analysis ($ \xi \in \RR^3 $).
\item [b)] Problem at $ s = r $ due to the introduction of $ J^{-1} $ (or $ J^{-1/p} $ in case of H\" older estimates).
This induces a singularity in the time variable.
\item [c)] Difficulty at $ \omega = - \Xi(s) / \vert  \Xi(s) \vert $ issued from negative powers of $ 1 + \omega 
\cdot \nu \circ \Xi(s) $ which are provided by $ J^{-1} $ (or $ J^{-1/p} $). This reflects some interplay between 
the phase space and the momentum component $ \Xi(s) $ of the characteristics.
\end{itemize}


\subsubsection{Three useful tools} \label{Three useful tools} The origin of negative powers of $ 1 + \omega 
\cdot \nu (\eta) $ inside (\ref{cD}) is manifold. For $ \eta = \xi $, they go with the gain of one derivative (in 
Lemma \ref{splitD1}). For $ \eta = \xi $ again, they are also due to the rate at which the transversality condition 
between the derivative $ \partial_t + \nu(\xi) \cdot \nabla_x $ and the light cones can degenerate for large values 
of $ \vert \xi \vert $. Or for $ \eta = \Xi(s) $, they could be issued from the inverse $ J^{-1/p} $ of the Jacobian. 
From a quantitative perspective, they may furnish large weights. We examine below what happens in diverse 
situations {\bf A}, {\bf B} and {\bf C}.

\smallskip

\noindent $ \bullet $ {\bf A}. This case deals with sup-norm estimates on the coefficients. It is particularly 
relevant for the study of $ {\rm D}^{\rm a}_{ne} $ and $ {\rm D}^{\rm a}_{nb} $ because only $ L^2 $- estimates are available 
when dealing with $ {\rm E} $ and $ {\rm B} $.

\begin{lem}\label{Maximalloss} [Maximal loss due to the proximity of $ \nu (\eta) $ to the light cone] We have 
\begin{equation} \label{simpleboundgen}  
0 \leq \bigl( 1 + \omega \cdot \nu (\eta) \bigr)^{-1} \leq 2 \, \langle \eta \rangle^2 \, , \qquad \forall \, (\omega,\eta) \in \mathbb S^2 
\times \RR^3 . 
\end{equation}
\end{lem} 

\begin{prof} It suffices to remark that
$ \bigl( 1 + \omega \cdot \nu (\eta) \bigr)^{-1} \leq \bigl( 1 - (\vert \eta \vert / \langle \eta \rangle) \bigr)^{-1} = \langle \eta \rangle \,
(\langle \eta \rangle + \vert \eta \vert ) . $
\end{prof}

\noindent The bound (\ref{simpleboundgen}) is (almost) optimal when $|\eta|$ is large, and  $\omega $ becomes close to 
$ - \eta / \vert \eta \vert $. However, for other values of $ \omega $, the upper bound (\ref{simpleboundgen}) furnishes only 
a rough control. The question is therefore to evaluate the impact of this singular factor in terms of the measure $ d \omega $. 

\smallskip

\noindent $ \bullet $ {\bf B}. Large weights may have a limited impact when they focus on a domain of small measure. This
is helpful at the level of $ {\rm D}^{\rm a}_0 $ and $ {\rm D}^{\rm a}_l $ because we know that $ {\rm f}_0 $ and $ {\rm f} $ are bounded
functions. Below, this effect appears after integration with respect to $ \omega $.
\begin{lem}\label{gainaveraging} [Gain after averaging along the sphere] For all $ (\delta,\eta) \in \RR_+ \times \RR^3 $, 
we have 
\begin{equation} \label{simpleboundgenaveraging}  
0 \leq \int_{\mathbb S^2} \bigl( 1 + \omega \cdot \nu (\eta) \bigr)^{-\delta} \ d \omega \lesssim \left \lbrace \begin{array}{lcl}
 \langle \eta \rangle^{2 \, (\delta-1)} & \text{when} & 1 < \delta , \\
1 + \ln \, \langle \eta \rangle & \text{when} & \delta = 1 , \\
1 & \text{when} & \delta < 1 .
\end{array} \right.  
\end{equation}
\end{lem} 

\begin{prof} By performing a rotation in $ \omega $, we can always assume that $ \nu(\eta) = (\vert \eta \vert / \langle 
\eta \rangle) \, {}^t (0,0,1) $. Then, we can work in spherical coordinates as in (\ref{sphericalcoordinates}) to see that
\[ \begin{array}{rl}
\displaystyle \int_{\mathbb S^2} \frac{d \upomega}{\bigl( 1 + \upomega \cdot \nu (\eta) \bigr)^\delta} \! \! \! & \displaystyle 
= \int_0^\pi \int_0^{2 \pi} \frac{\sin \vartheta \ d \vartheta \, d \varphi}{\bigl( 1 + \vert \eta \vert \, \cos \vartheta / \langle 
\eta \rangle \bigr)^\delta} = (2 \pi) \ \frac{\langle \eta \rangle}{\vert \eta \vert} \int_{- \vert \eta \vert / \langle \eta \rangle
}^{+\vert \eta \vert / \langle \eta \rangle} \frac{ds}{(1-s)^\delta} \\
\  & \displaystyle =  (2 \pi) \ \frac{\langle \eta \rangle}{\vert \eta \vert} \times  \left \lbrace \begin{array}{lcl}
 \displaystyle  \frac{1}{1-\delta} \, \Bigl \lbrack - \Bigl( 1 - \frac{\vert \eta \vert}{\langle \eta \rangle} \Bigr)^{1-\delta}
 + \Bigl( 1 + \frac{\vert \eta \vert}{\langle \eta \rangle} \Bigr)^{1-\delta} \Bigr \rbrack & \text{when} & \delta \not = 1 , \\
 \displaystyle - \ln \, \Bigl( 1 - \frac{\vert \eta \vert}{\langle \eta \rangle} \Bigr) + \ln \, \Bigl( 1 + \frac{\vert \eta \vert}{\langle 
 \eta \rangle} \Bigr)& \text{when} & \delta = 1 .
\end{array} \right.  
\end{array} \]
By this way, we can easily deduce (\ref{simpleboundgenaveraging}).\end{prof}

\smallskip

\noindent $ \bullet $ {\bf C}. On different occasions, we will have to evaluate the distance between $ \nu (\eta) $ and $ \omega $.
This will serve for instance to control the vector valued functions $ {\rm W}_{ne} $ and $ {\rm W}_{nb} $.
\begin{lem}\label{Comparisonsingularweights} [Comparison between $ \nu (\eta) $ and $ \omega $] We have 
\begin{equation} \label{comparisonboundgen}  
\vert \nu (\eta) + \omega \vert \leq \sqrt 2 \ \bigl( 1 + \omega \cdot \nu (\eta) \bigr)^{1/2} \, , \qquad \forall \, (\omega,\eta) \in 
\mathbb S^2 \times \RR^3 . 
\end{equation}
\end{lem} 

\begin{prof} This is just because $ \vert \nu (\eta) + \omega \vert^2 = \vert \nu (\eta) \vert^2 + 2 \ \omega \cdot \nu (\eta)  
+ 1 \leq 2 \ \bigl( 1 + \omega \cdot \nu (\eta) \bigr) $.
\end{prof}


\subsubsection{Bounds for the weight functions} \label{Radon microlocal analysis of the weight functions}  
The purpose of this paragraph is to evaluate carefully the amplitudes of $ \vert {\rm W}_0 \vert $, $ \vert {\rm W}_l \vert $, $ \vert 
{\rm W}_{ne} \vert $ and $ \vert {\rm W}_{nb} \vert $. 

\smallskip

\noindent $ \bullet $ {\it Study of $ \vert {\rm W}_0 \vert $.} From (\ref{weights0}), remark that
\[ {\rm W}_0 (s,\omega,\xi) := \frac{s}{4 \pi} \ \Bigl( 1 - \frac{\bigl( \nu \circ \Xi (s) + \omega \bigr) \cdot  \bigl(\nu(\xi) + 
\omega \bigr)}{1 + \omega \cdot \nu(\xi)} \Bigr) . \]
With (\ref{comparisonboundgen}), it is obvious that 
\begin{equation} \label{comparisonboundgensuitefor0} 
\vert {\rm W}_0 (s,\omega,\xi) \vert \leq \frac{s}{4 \pi} \ \Bigl( 1 + 2 \ \frac{\bigl( 1 + \omega \cdot \nu \circ \Xi (s) \bigr)^{1/2}}
{\bigl( 1 + \omega \cdot \nu(\xi) \bigr)^{1/2}} \Bigr) .
\end{equation}
 
\smallskip

\noindent $ \bullet $ {\it Study of $ \vert {\rm W}_l \vert $} \label{Study ofWi} From (\ref{weightsi}), observe that
\[ {\rm W}_l (s,\omega,\xi) := \frac{1}{4 \pi} \ \frac{1}{\langle \xi \rangle^2} \, \Bigl( \frac{1}{\bigl( 1+ \omega \cdot \nu(\xi) \bigr)} 
- \frac{\bigl( \nu \circ \Xi (s) + \omega \bigr) \cdot \bigl( \omega + \nu(\xi) \bigr)}{\bigl( 1+ \omega \cdot \nu(\xi) \bigr)^2} \Bigr) . \]
Then, as a corollary of Lemmas \ref{Maximalloss} and 
 \ref{Comparisonsingularweights}, we can assert that
 \begin{equation} \label{cestpourcdiameprem}  
\vert {\rm W}_l (s,\omega,\xi) \vert \leq \frac{1}{\pi} \ \Bigl( 1 + \frac{\bigl( 1 + \omega \cdot \nu \circ \Xi (s) \bigr)^{1/2}}
{\bigl( 1 + \omega \cdot \nu(\xi) \bigr)^{1/2}} \Bigr) . 
\end{equation}

\smallskip

\noindent $ \bullet $ {\it Study of $ \vert {\rm W}_{ne} \vert $ and $ \vert {\rm W}_{nb} \vert $} \label{Study ofWne} We start 
by computing
\[ \begin{array}{rl}
\displaystyle \nabla_\xi \Bigl \lbrace \frac{\nu \circ \Xi (s) \cdot \bigl( \nu (\xi) + \omega \bigr)}{1 + \omega \cdot \nu (\xi)}  
\Bigr \rbrace = \! \! \! & \displaystyle \frac{1}{\langle \xi \rangle \, \bigl(1 + \omega \cdot \nu (\xi) \bigr)} 
\ \bigl \lbrack \nu \circ \Xi(s) - \bigl( \nu (\xi) \cdot \nu \circ \Xi(s) \bigr) \, \nu (\xi) \bigr \rbrack \\
\displaystyle \ \! \! \! & 
\displaystyle - \, \frac{\nu \circ \Xi (s) \cdot \bigl( \nu (\xi) + \omega \bigr)}{\langle \xi \rangle \, \bigl(1 + 
\omega \cdot \nu (\xi) \bigr)^2} \ \bigl \lbrack \omega- \bigl( \nu (\xi) \cdot \omega \bigr) \, \nu (\xi) \bigr \rbrack .
 \end{array} \]
This is a vector valued function which can be decomposed with respect to the moving ``frame" made of the three directions  
$ \omega $, $ \nu \circ \Xi(s) + \omega$ and $ \nu (\xi) + \omega $. This gives rise to
\[ \begin{array}{rl}
\displaystyle \nabla_\xi \Bigl \lbrace \frac{\nu \circ \Xi (s) \cdot \bigl( \nu (\xi) + \omega \bigr)}{1 + \omega \cdot \nu (\xi)}  
\Bigr \rbrace \! \! \! & \displaystyle = - \, \frac{1}{\langle \xi \rangle} \ \frac{1+ \omega \cdot \nu \circ \Xi(s)}{1 + \omega \cdot 
\nu (\xi)} \ \omega \\
\ & \displaystyle \! \! \!  \! \! \! + \, \frac{1}{\langle \xi \rangle} \ \frac{1}{1 + \omega \cdot \nu (\xi)} \ \bigl( \nu \circ \Xi(s) + \omega \bigr) \\
\ & \displaystyle \! \! \!  \! \! \! + \, \frac{1}{\langle \xi \rangle} \ \Bigl( + \frac{1+ \omega \cdot \nu \circ \Xi(s)}{1 + \omega \cdot 
\nu (\xi)} - \frac{\bigl( \nu (\xi) + \omega \bigr) \cdot \bigl(\omega + \nu \circ \Xi(s) \bigr)}{\bigl( 1 + \omega \cdot \nu (\xi) \bigr)^2} 
\Bigr) \ \bigl( \nu (\xi) + \omega \bigr) . 
 \end{array} \]
 The three vectors $ \omega $, $ \nu \circ \Xi(s) + \omega$ and $ \nu (\xi) + \omega $ are clearly uniformly bounded (by $ 2 $) 
 as functions of $ (s,\omega,\xi) $. When doing the above decomposition, we can observe that the coefficient in factor of $ \omega $ 
 is small (at least smaller than what appears at first sight ) due to various cancellations that are revealed during its decomposition. 
 On the other hand, the sizes of $ \nu \circ \Xi(s) + \omega$ and $ \nu (\xi) + \omega $ can be estimated through Lemma 
 \ref{Comparisonsingularweights}. Briefly, from (\ref{cDndecomponbinftyssb}), we can deduce that
\begin{equation} \label{cestpoureeeee}  
\vert {\rm W}_{n \star} (r,s,\omega,\xi) \vert \leq \frac{3 \, \sqrt 2}{2 \, \pi} \ \frac{(s-r)}{\langle \xi \rangle} \ \frac{\bigl(1+ \omega \cdot 
\nu \circ \Xi(s) \bigr)^{1/2}}{\bigl(1 + \omega \cdot \nu (\xi) \bigr)} \, , \qquad \forall \star \in \{ e,b \} . 
\end{equation}


\subsection{Weighted integrals} \label{Weighted integrals} The goal of this subsection is to estimate the contributions 
provided by the $ {\rm D}^{\rm a}_\star $ with $ \star \in \{ 0,h,l,ne,nb \} $. This will be done in separate paragraphs, one 
for each $ {\rm D}^{\rm a}_\star $. Before starting, we would like to accurately define the scope of our discussion. Indeed,
different courses of action are possible when studying the $ {\rm D}^{\rm a}_\star $.

\noindent In the perspective of continuation criteria, one might attempt to minimize the powers of $ \langle {\rm P} \rangle_\infty $ 
needed to control the $ {\rm D}^{\rm a}_\star $. As in \cite{MR3437855,MR3291372,MR3721415}, it seems that additional 
conditions (to be identified) are needed to recover the global existence. This interesting option is not pursued here.

\noindent  On the other hand, to fit in with Proposition \ref{theorempms}, one can insist on the role of $ {\rm S}_0 $ when 
looking at the momentum increment. This is the path that we follow below. However, when doing this, a cautionary note 
is in order. This is because the functional $ \cN $ (and thereby $ {\rm S}_0 $) involves different types of norms: first, the 
sup-norm concerning $ {\rm f}_0 $; and secondly, the Sobolev norm $ W^{1,\bar p} $ in the case of $ ({\rm E}_0,{\rm B}_0) $. 
Let us explain the origin of this distinction:
 \begin{itemize}
 \item [-] The handling of the $ {\rm D}^{\rm a}_\star $ with $ \star \in \{ 0,l,ne,nb \} $ does not prove to be demanding  in terms 
 of regularity.  It will only require the use of $ \parallel f_0 \parallel_0 $ and $ \text{\small $ \pmb{\mathscr{E}} $}_{\! 0} $.
\item [-] The manipulation of $ {\rm D}^{\rm a}_h $ could be based just on a sup-norm estimate concerning $ {\rm E}_h $. It is 
the transcription of such uniform bound in terms of the initial data $ ({\rm E}_0,{\rm B}_0) $ that generates the implementation 
of $ W^{1,\bar p} $. By the way, note that the use of $ \bar p \in ]3/2,3[ $ instead of $ \bar p = + \infty $ (or even higher levels of 
regularity) is a subtle refinement that will be clarified in Paragraph \ref{Control ofDh}.
\end{itemize}

\noindent In other words, the focus is on the {\it minimal regularity} required on the complete {\it initial data} $ {\rm U}_0 $ 
in order to control the quantities $ {\rm D}^{\rm a}_\star $. In so doing, for the sake of simplicity, we have highlighted in 
Proposition \ref{theorempms} the role of the sole parameter $ {\rm S}_0 $, which serves in fact to cover different aspects. 
For instance, from (\ref{cestpourf0}), we get that
 \begin{equation}\label{consertotenRVMgencotef0} 
\int_{\RR^3} \int_{\RR^3} \langle \xi \rangle \ {\rm f}_0 (x,\xi) \ dx \, d\xi  \lesssim \langle {\rm P} \rangle_0^7 \  {\rm S}_0 \, , 
\qquad \langle {\rm P} \rangle_0 := (1+ {\rm P}_0^2)^{1/2} ,
\end{equation}
where the symbol $ \lesssim $ is for $ \leq C $ with some (universal) constant $ C $ not depending on $ {\rm P}_0 $ 
or $ {\rm S}_0 $. On the other hand, the Sobolev embedding theorem (which holds since $ \bar p < 3 $) ensures that 
$ {\rm E}_0 \in L^{\tilde p} $ for some $ \tilde p > 3 $. By interpolation, we get that $ W^{1,\bar p}  \hookrightarrow L^2 $. 
Now, since $ \cN ({\rm U}_0) \leq {\rm S}_0 $, we have
 \begin{equation}\label{incluconse} 
\int_{\RR^3} |{\rm E}_0(x)|^2 \ dx + \int_{\RR^3} |{\rm B}_0(x)|^2 \ dx \lesssim {\rm S}_0^2 .
\end{equation}
Coming back to (\ref{consertotenRVMgen}) and using (\ref{repetcompact space}), it follows that (for fixed $ {\rm P}_0 $ 
and small $ {\rm S}_0 $)
\begin{equation}\label{estimonslee0} 
\text{\small $ \pmb{\mathscr{E}} $}_{\! 0} \lesssim \langle {\rm P} \rangle_0^7 \ {\rm S}_0 +  {\rm S}_0^2 \lesssim {\rm S}_0 .
\end{equation}
It is clear that the $ {\rm S}_0 $ of (\ref{consertotenRVMgencotef0}) comes from the sup-norm of $ {\rm f}_0 $, while 
the $ {\rm S}_0 $ of (\ref{incluconse}) is issued from the $ W^{1,\bar p} $-norm of $ ({\rm E}_0,{\rm B}_0) $. These 
contributions are mixed (with different powers) at the level of (\ref{estimonslee0}). To avoid having to introduce too 
much material, our decision is to not make the various origins of $ {\rm S}_0 $ apparent in the final statement.
But the interested reader can easily trace $ {\rm S}_0 $ in the forthcoming analysis.


\subsubsection{Study of $ {\rm D}^{\rm a}_0 $} \label{Control ofD0} By definition, $ {\rm D}^{\rm a}_0 $ is built on $ {\rm U}_0 $.
It is therefore a known quantity. Still, it is interesting to estimate $ {\rm D}^{\rm a}_0 $ to see how this works and to have 
access to its time behavior.

\begin{lem} \label{controllem ofD0} [Control of $ {\rm D}^{\rm a}_0 $] For all $ p \in ]2,+\infty] $, we have
\renewcommand\arraystretch{1.6}
\begin{equation}\label{thisisfor0} 
\begin{array}{rl}
\displaystyle {\rm D}_0^{\rm a} (t) \! \! \! & \displaystyle = \int_0^t \int_{\mathbb S^2} \int_{\RR^3} \vert {\rm W}_0 
(s,\omega,\xi) \vert \ {\rm f}_0 \bigl( X (s) +  s \omega,\xi \bigr) \ ds \, d \omega \, d \xi  \\
\ & \displaystyle \lesssim \, {\rm S}_0^{1-(1/p)} \ \text{\small $ \pmb{\mathscr{E}} $}_{\! 0}^{1/p} \ t^{2-(3/p)} \ (1 + {\rm P}_0^2
)^{(3/2)-(2/p)} .
\end{array}
  \end{equation}
  \renewcommand\arraystretch{1}
 \end{lem}

\begin{prof} Let $ q $ be the conjugate index of $ p $. Remark that
\[ {\rm D}_0^{\rm a} (t) = \int_0^t \int_{\mathbb S^2} \int_{\RR^3} \frac{s^{-(2/p)} \, \vert {\rm W}_0 (s,\omega,\xi) \vert}{\langle 
\xi \rangle^{1/p} \, \bigl(1+ \omega \cdot \nu \circ \Xi(s) \bigr)^{1/p}} \times J(0,s,\omega)^{1/p} \, \langle \xi \rangle^{1/p} \ 
{\rm f}_0 \bigl(X (s) + s \omega,\xi \bigr) \ ds \, d \omega \, d \xi . \]
As explained in Paragraph \ref{Comparison of measures}, see (\ref{jacobiadit}), we can assert that
\[ \parallel J^{1/p} \, \langle \xi \rangle^{1/p} \, {\rm f}_0 \bigl(X (s) + s \omega,\xi \bigr) \parallel_{L^p ([0,t] \times 
\mathbb S^2 \times\RR^3)} \leq \Bigl( \int_{\RR^3} \int_{\RR^3}  \langle \xi \rangle \, {\rm f}_0(x,\xi)^p \ dx \, d\xi \Bigr)^{1/p} \! \leq 
{\rm S}_0^{1-(1/p)} \ \text{\small $ \pmb{\mathscr{E}} $}_{\! 0}^{1/p} . \]
By H\"older's inequality, exploiting (\ref{comparisonboundgensuitefor0}) and the condition $ p > 2 $ for the second term in the 
right hand side of (\ref{comparisonboundgensuitefor0}), we find that
\[ \begin{array}{rl}
\displaystyle {\rm D}_0^{\rm a} (t) \lesssim \! \! \! \! & \displaystyle {\rm S}_0^{1-(1/p)} \ \text{\small $ \pmb{\mathscr{E}} $}_{\! 0}^{1/p} \
\Bigl( \int_0^t \int_{\mathbb S^2} \int_{\vert \xi \vert \leq {\rm P}_0} \frac{s^{q-(2q/p)} \, \langle \xi \rangle^{-(q/p)} }{\bigl(1+ 
\omega \cdot \nu \circ \Xi(s) \bigr)^{q/p}} \ ds \, d \omega \, d \xi \Bigr)^{1/q} \\
+ \! \! \! \! & \displaystyle {\rm S}_0^{1-(1/p)} \ \text{\small $ \pmb{\mathscr{E}} $}_{\! 0}^{1/p} \ \Bigl( \int_0^t \int_{\mathbb S^2} 
\int_{\vert \xi \vert \leq {\rm P}_0} \frac{s^{q-(2q/p)} \, \langle \xi \rangle^{-(q/p)} }{\bigl(1+ \omega \cdot \nu (\xi) \bigr)^{q/p}} \ 
ds \, d \omega \, d \xi \Bigr)^{1/q} .
\end{array} \]
We first integrate with respect to $ \omega $. Since $ q/p < 1 $ (since again $ p > 2 $), from Lemma \ref{gainaveraging}, we have
\[ {\rm D}_0^{\rm a} (t) \lesssim {\rm S}_0^{1-(1/p)} \ \text{\small $ \pmb{\mathscr{E}} $}_{\! 0}^{1/p}  \ \Bigl( \int_0^t s^{q-(2q /p)} \ ds 
\Bigr)^{1/q} \ \Bigl(\int_{\vert \xi \vert \leq {\rm P}_0} \langle \xi \rangle^{-(q/p)} \ d \xi \Bigr)^{1/q} . \]
Since $ q-(2q/p) > 0 $, we end up with (\ref{thisisfor0}). 
 \end{prof} 


\subsubsection{Study of $ {\rm D}_h^{\rm a} $} \label{Control ofDh} Remark that
\[ 0 \leq {\rm D}_h^{\rm a} := \int_0^t \vert \nu \circ \Xi (s) \vert \ \vert {\rm E}_h \bigl( s,X(s)\bigr) \vert \ ds \leq \int_0^t \vert {\rm E}_h 
\bigl(s,X(s) \bigr) \vert \ ds \, , \qquad t \leq T . \]
As already noted, the field $ ({\rm E}_0,{\rm B}_0) $ has an impact on all the $ {\rm D}_\star $ as well as $ X $ because the Lorentz 
force $ {\rm F} $ is built with (\ref{defdeEetB}), where $ {\rm A}_h $ (and therefore $ {\rm E}_0 $ and $ {\rm B}_0 $)  is activated. It is 
particularly interesting to further examine its influence on $ {\rm D}_h^{\rm a} $. There are two ways of thinking.

\smallskip

\noindent $ \bullet $ The {\it linear} viewpoint. That is concentrating on the only role of $ {\rm E}_h $. This method could be based 
on the following observations marked by $ a) $ and $ b) $. 

\noindent $ a) $ In the local (in time) version of Proposition \ref{theorempms}, it turns out that the restriction on $ ({\rm E}_0,{\rm B}_0) $
may be exchanged with the mild assumption
\begin{equation}\label{mildmildeh} 
{\rm E}_h \in L^1_{loc} \bigl( \RR_+ ; L^\infty (\RR^3) \bigr) .
  \end{equation}
\noindent $ b) $ The field $ {\rm E}_h $ allows to absorb the main contribution brought by $ ({\rm E}_0,{\rm B}_0) $. The wave equation
inside (\ref{homogeneous version}) is {\it linear} and completely {\it decoupled} from (\ref{VlasoveqRVMD})-(\ref{maxwelleqRVMD}). 
Moreover, the information (\ref{mildmildeh}) is available for a whole range of bounded initial data $ ({\rm E}_0,{\rm B}_0) $. The 
$ W^{1,\bar p} $-condition is not necessarily (and also not sufficient) for that. \hfill $ \circ $

\smallskip

\noindent It would be enough to deal with (\ref{mildmildeh}), but this would be disappointing in terms of the Cauchy problem for the 
RVM system. Moreover, that would ignore a subtle nuance arising between the time integration of $ {\rm E}_h \bigl( s,X(s)\bigr) $
and the one of $  {\rm E}_h (s,x) $. First, recall that (\ref{mildmildeh}) is not easy to find \cite{MR1240537}. In particular, the endpoint 
Strichartz estimate
\begin{equation}\label{Strichartz estimate}
 \parallel {\rm E}_h \parallel_{L^2_t L^\infty_x} \lesssim  \parallel {\rm E}_0 \parallel_{\dot H^1} +  \parallel \nabla_x \times 
{\rm B}_0 \parallel_{L^2}  
\end{equation}
is known to be false \cite{MR1646048}. Furthermore, for data $ {\rm E}_0 $ and $ {\rm B}_0 $ in $ W^{1,\bar p} (\RR^3) $, it is not clear 
that the integral of $ {\rm E}_h $ along {\it any} space-time curve makes sense. At this low level of regularity, the meaning of $ {\rm D}_h $ 
cannot be based solely on the properties of $ {\rm E}_h $. We have to change the perspective.

\smallskip

\noindent $ \bullet $ The {\it nonlinear} viewpoint. This means to look at the expression $ {\rm D}_h^{\rm a} $ as a nonlinear functional, 
pursuant to the influence of $ X $. This approach has a clear advantage. The counter-examples to the inequality (\ref{Strichartz estimate}) 
are exhibited by concentrating solutions along the light cone. But the special curve $ \bigl \lbrace \bigl(s,X(s) \bigr) ; s \in [0,t] \bigr \rbrace $ 
intersects the light cone transversally, and therefore the time integration of $ {\rm E}_h $ when computing $ {\rm D}_h^{\rm a} $ reduces 
this alignment effects. This (relativistic) feature is a key ingredient because it allows to exploit (as in \cite{Bouchut,MR816621,MR3356994,MR1231427}) 
the different speeds of propagation between the slow particles and the fields (which propagate at the speed of light). It is crucial here to 
make sense of $ {\rm D}_h $ in the context of (\ref{mathcal N}). \hfill $ \circ $

\begin{lem} \label{controllem ofDah} [Control of $ {\rm D}_h^{\rm a} $] We have
\begin{equation}\label{dhtodofinal}
\qquad {\rm D}_h^{\rm a} \lesssim {\rm S}_0 \ (t+1) \ \langle {\rm P} \rangle_\infty (t)^{2/3} . 
  \end{equation}
 \end{lem}

\begin{prof} Let us turn to another interpretation of $ {\rm D}_h $. The field $ {\rm E}_h $ is given
by Kirchhoff's formula
\[ {\rm E}_h (t,x) = t \, \cM_t (\nabla_x \times {\rm B}_0) + \part_t \bigl( t \, \cM_t ({\rm E}_0) \bigr) = t \, \cM_t (\nabla_x 
\times {\rm B}_0) + t \, \cM_t (\omega \cdot \nabla_x {\rm E}_0) + \cM_t ({\rm E}_0) , \]
where, for $ t \in \RR_+ $, we have introduced the mean operator $ \cM_t : L^\infty (\RR^3) \rightarrow L^\infty (\RR^3) $
defined by
\[ \cM_t (k) (x) := \frac{1}{4 \pi} \int_{\mathbb S^2} k(x+ t \omega) \ d \omega \, , \qquad \vert \! \vert \! \vert \cM_t 
\vert \! \vert \! \vert_{\cL(L^\infty)} \leq 1 . \]
After substitution, this means that
\begin{equation}\label{dhtodo}
\begin{array}{rl}
 {\rm D}_h^{\rm a} \lesssim \! \! \! & \displaystyle \sum_{i=1}^3 \int_0^t \int_{\mathbb S^2} s \ \bigl( \vert \partial_{x_i} 
 {\rm E}_0 \vert + \vert \partial_{x_i} {\rm B}_0 \vert \bigr) \bigl(s,X(s) + s \omega \bigr) \ ds \, d \omega \\
& \displaystyle + \int_0^t \int_{\mathbb S^2} \vert {\rm E}_0 \vert \bigl(s,X(s) + s \omega \bigr) \ ds \, d \omega . 
\end{array}
\end{equation}
First, consider the contribution coming from $ \partial_{x_i} {\rm E}_0 $ (do the same with $ \partial_{x_i} {\rm B}_0 $). Since
$ {\rm E}_0 $ is selected in $ W^{1,\bar p} $, we know that $ \partial_{x_i} {\rm E}_0 \in L^{\bar p} $. Let $ \bar q := \bar p / 
(\bar p -1) $ be the H\"older conjugate of $ \bar p $. By H\" older's inequality, we can assert that
\[ \begin{array}{l}
\displaystyle \int_0^t \int_{\mathbb S^2} \! s \, \vert \partial_{x_i} {\rm E}_0 \vert \bigl(s,X(s) + s \omega \bigr) \, ds \, d \omega \\
\qquad \qquad \qquad \quad \displaystyle \lesssim \Bigl( \int_0^t \int_{\mathbb S^2} \! J \, \vert \partial_{x_i} {\rm E}_0 \vert^{\bar p}  
\bigl(s,X(s) + s \omega \bigr) \, ds \, d \omega \Bigr)^{1/\bar p} \, \Bigl( \int_0^t \int_{\mathbb S^2} \frac{s^{\bar q}}{J^{\bar q/\bar p}} \, 
ds \, d \omega \Bigr)^{1/\bar q} \\
\qquad \qquad \quad \qquad \displaystyle \lesssim {\rm S}_0 \ \parallel {\rm E}_0 \parallel_{W^{1,\bar p}} \, \Bigl( \int_0^t 
\int_{\mathbb S^2} \frac{s^{\bar q - (2 \bar q / \bar p)}}{\bigl( 1 + \omega \cdot \nu \circ \Xi(s) \bigr)^{\bar q / \bar p} } \ ds \, 
d \omega \Bigr)^{1/\bar q} .
\end{array} \]
To estimate the right hand side, we start by applying Lemma \ref{gainaveraging} with $ 1 \leq \delta = \bar q / \bar p \leq 2 $; then, 
we use (\ref{supplball of radius}) and finally, to obtain the time integrability near $ s= 0 $, it suffices to remark that we have 
$ \bar q - (2 \bar q / \bar p) > -1 $ (since $ 3/2 < \bar p $). This furnishes
\[ \Bigl( \int_0^t \int_{\mathbb S^2} \frac{s^{\bar q - (2 \bar q / \bar p)}}{\bigl( 1 + \omega \cdot \nu \circ \Xi(s) \bigr)^{\bar q / \bar p} } \ 
ds \, d \omega \Bigr)^{1/\bar q} \lesssim t^{2-(3/\bar p)} \ \langle {\rm P} \rangle_\infty (t)^{(4/\bar p)-2} . \]
Since  $ 2-(3/\bar p) < 1 $ and $ (4/\bar p)-2 \leq 2/3 $ when $ 3/2 <\bar p \leq 2 $, we find the right hand side of (\ref{dhtodofinal}).
Secondly, we look at the contribution brought inside (\ref{dhtodo}) by $ {\rm E}_0 $. The difficulty is that \og $ s $ \fg{} 
is no more in factor. But this lost may be compensated by extra integrability concerning $ {\rm E}_0 $.

\noindent Indeed, the Sobolev embedding theorem gives $ {\rm E}_0 \in L^{\tilde p} $ for some $ \tilde p > 3 $. Let $ \tilde q := \tilde p / (\tilde p -1) $ be 
the H\"older conjugate of $ \tilde q $ so that $ 2 \tilde q / \tilde p < 1 $. By H\" older's inequality, we have
\[ \begin{array}{l}
\displaystyle \int_0^t \int_{\mathbb S^2} \vert {\rm E}_0 \vert \bigl(s,X(s) + s \omega \bigr) \ ds \, d \omega \\
\displaystyle \qquad \qquad \qquad \quad  \lesssim \Bigl( \int_0^t \int_{\mathbb S^2} J \, \vert {\rm E}_0 \vert^{\tilde p}  
\bigl(s,X(s) + s \omega \bigr) \ ds \, d \omega \Bigr)^{1/\tilde p} \, \Bigl( \int_0^t \int_{\mathbb S^2} J^{- \tilde q / \tilde p} \ 
ds \, d \omega \Bigr)^{1/\tilde q} \\
\qquad \qquad \qquad \quad \displaystyle \lesssim {\rm S}_0 \ \parallel {\rm E}_0 \parallel_{L^{\tilde p}} \ \Bigl( 
\int_0^t \int_{\mathbb S^2} \frac{s^{- 2 \tilde q / \tilde p}}{\bigl( 1 + \omega \cdot \nu \circ \Xi(s) \bigr)^{\tilde q / 
\tilde p} } \ ds \, d \omega \Bigr)^{1/\tilde q} \\
\qquad \qquad \qquad \quad  \displaystyle \lesssim {\rm S}_0 \ \parallel {\rm E}_0 \parallel_{W^{1,\bar p}} \ 
t^{-(2/ \tilde p)+(1/\tilde q)} \lesssim {\rm S}_0 \ (t+1) . 
\end{array} \]
Again, this is consistent with (\ref{dhtodofinal}).
\end{prof} 


\subsubsection{Study of $ {\rm D}_l^a $} \label{Control ofDb} Observe that $ {\rm W}_l (s,\cdot) $ is an odd function. Thus, there 
is no benefit from the sign condition on $ {\rm f} $ when studying $ {\rm D}_l^a $.

\begin{lem} \label{controllem ofDgfaejk} [Control of $ {\rm D}^{\rm a}_l $] For all $ p \in ]3,+\infty] $, we have
\renewcommand\arraystretch{1.6}
\begin{equation}\label{thisisforthel} 
\begin{array}{rl}
\displaystyle {\rm D}_l^a (t) \! \! \! & \displaystyle = \int_0^t \Bigl( \int_r^t \int_{\mathbb S^2} \int_{\RR^3} {\rm W}_l (s,\omega,\xi) \  
{\rm f} \bigl(r,X(s) + (s-r) \omega,\xi \bigr) \ ds \, d \omega \, d \xi \Bigr) \, dr \\
\ & \displaystyle \lesssim \, {\rm S}_0^{1-(1/p)} \ \text{\small $ \pmb{\mathscr{E}} $}_{\! 0}^{1/p} \ t^{1-(3/p)} \ \int_0^t \langle 
{\rm P} \rangle_\infty (r)^{3-(4/p)} \ dr . 
\end{array}
  \end{equation}
  \renewcommand\arraystretch{1}
 \end{lem}

\begin{prof} The inequalities (\ref{comparisonboundgensuitefor0}) and 
(\ref{cestpourcdiameprem}) are similar, except that $ (s-r) $ does not appear in factor in the right hand side of (\ref{cestpourcdiameprem}).
Due to this additional difficulty, there are some nuances in comparison to what has been done in Paragraph \ref{Control ofD0}. From 
(\ref{cestpourcdiameprem}), as soon as $ p> 2 $, we have
\[ \begin{array}{rl}
\displaystyle {\rm D}_l^a (t) \lesssim \! \! \! \! & \displaystyle \int_0^t \Bigl( \int_r^t \int_{\mathbb S^2} \int_{\RR^3}  \frac{(J^{1/p} \ \langle 
\xi \rangle^{1/p} \ {\rm f}) \bigl(r,X(s) +  (s-r) \omega,\xi \bigr)}{(s-r)^{2/p} \ \bigl( 1 + \omega \cdot \nu \circ \Xi(s) \bigr)^{1/p} \ \langle \xi 
\rangle^{1/p}} \ ds \, d \omega \, d \xi \Bigr) \, dr \\
\ & \displaystyle + \int_0^t \Bigl( \int_r^t \int_{\mathbb S^2} \int_{\RR^3}  \frac{(J^{1/p} \ \langle 
\xi \rangle^{1/p} \ {\rm f}) \bigl(r,X(s) +  (s-r) \omega,\xi \bigr)}{(s-r)^{2/p} \ \bigl( 1 + \omega \cdot \nu (\xi) \bigr)^{1/2} \ \langle \xi 
\rangle^{1/p}} \ ds \, d \omega \, d \xi \Bigr) \, dr .
\end{array}\]
We apply H\"older's inequality 
\[ \begin{array}{rl}
\displaystyle {\rm D}_l^a (t) \lesssim \! \! \! \! & \displaystyle {\rm S}_0^{1-(1/p)} \ \text{\small $ \pmb{\mathscr{E}} $}_{\! 0}^{1/p} \ 
\int_0^t \Bigl( \int_r^t \int_{\mathbb S^2} \int_{\vert \xi \vert \leq {\rm P} (r)}  \frac{ds \, d \omega \, d \xi }{(s-r)^{2q/p} \ \bigl( 1 + 
\omega \cdot \nu \circ \Xi(s) \bigr)^{q/p} \ \langle \xi \rangle^{q/p}} \Bigr)^{1/q} \, dr \\
+ \! \! \! \! & \displaystyle {\rm S}_0^{1-(1/p)} \ \text{\small $ \pmb{\mathscr{E}} $}_{\! 0}^{1/p} \ \int_0^t \Bigl( \int_r^t \int_{\mathbb S^2} 
\int_{\vert \xi \vert \leq {\rm P} (r)}  \frac{ds \, d \omega \, d \xi }{(s-r)^{2q/p} \ \bigl( 1 + \omega \cdot \nu (\xi) \bigr)^{q/2} \ \langle 
\xi \rangle^{q/p}} \Bigr)^{1/q} \, dr .
\end{array}\]
Knowing that $ q/p < 1 $ and $ q/2 < 1 $, we can integrate with respect to $ \omega $ through Lemma \ref{gainaveraging}. Then, we 
have to take $ p > 3 $ (so that $ 2q/p < 1 $) in order to be sure that  the integral with respect to $ ds $ is convergent, giving 
rise to (\ref{thisisforthel}).
 \end{prof} 


\subsubsection{Study of $ {\rm D}^{\rm a}_{ne} $ and $ {\rm D}^{\rm a}_{nb} $} \label{Control ofDn} The expressions 
$ {\rm D}^{\rm a}_{ne} $ and $ {\rm D}^{\rm a}_{nb} $ with $ \star \in \{ e,b\} $ can be viewed as bilinear forms. Indeed, 
with $ {\rm W}_{n\star}^{\rm a} := \vert {\rm W}_{n\star} \vert $ where $ {\rm W}_{n\star} $ is as in (\ref{cDndecomponbinftyssb}) 
together with the conventions $ {\rm G}_e := \vert {\rm E} \vert  $ and $ {\rm G}_b := \vert {\rm B} \vert $, we have to deal with
\[ {\rm D}_{n\star}^{\rm a} (t) = \int_0^t \Bigl( \int_r^t \int_{\mathbb S^2} \int_{\vert \xi \vert \leq \langle {\rm P} \rangle_\infty (r)} 
\bigl \lbrack (J^{-1/2} \, {\rm W}_{n\star}^{\rm a} \, {\rm f}) \cdot (J^{1/2} \, {\rm G}_\star) \bigr \rbrack \bigl(r,X(s) + (s-r) \omega,
\xi \bigr) \ ds \, d \omega \, d \xi \Bigr) \, dr . \]

\begin{lem} \label{controllem ofDaneee} [Control of $ {\rm D}_{n\star}^{\rm a} $] For all $ p \in ]3,+\infty] $, we have
\begin{equation} \label{atthisstage}   
{\rm D}_{n\star}^{\rm a} (t) \lesssim {\rm S}_0^{1-(1/2p)} \ \text{\small $ \pmb{\mathscr{E}} $}_{\! 0}^{(1/2)+(1/2p)} \ t^{(1/2) - (3/2 p)} 
 \int_0^t \langle {\rm P} \rangle_\infty (r)^{4-(3/p)} \ dr  . 
\end{equation}
 \end{lem}

\begin{prof} By Cauchy-Schwarz inequality, since $ {\rm G}_\star (r,\cdot) \in L^2(\RR^3) $ with a bound which can be 
viewed as coming from (\ref{consertotenRVM}), we have 
\[ \begin{array}{rl}
\displaystyle {\rm D}_{n\star}^{\rm a} (t) \! \! \! & \displaystyle \lesssim {\rm S}_0^{1/2} \ \text{\small $ \pmb{\mathscr{E}} $}_{\! 0}^{1/2} \\ 
\ & \displaystyle \times \int_0^t \langle {\rm P} \rangle_\infty (r)^{3/2} \ \Bigl( \int_r^t \int_{\mathbb S^2} \int_{\vert \xi \vert \leq 
\langle {\rm P} \rangle_\infty (r)} \! \! (J^{-1} \, {\rm W}_{n\star}^2 \, {\rm f}) \bigl(r,X(s) + (s-r) \omega,\xi \bigr) \ ds \, d \omega \, 
d \xi \Bigr)^{1/2} \, dr . \end{array}  \]
From (\ref{cestpoureeeee}), for $ \star \in \{ e,b \} $, we obtain that
\[ J^{-1} \, {\rm W}_{n\star}^2 \lesssim \langle \xi \rangle^{-2} \ \bigl(1 + \omega \cdot \nu (\xi) \bigr)^{-2} . \]
Again, we select some $ p > 3 $. By H\"older's inequality, we get
\[ \begin{array}{l}
\displaystyle \int_r^t \int_{\mathbb S^2} \int_{\vert \xi \vert \leq \langle {\rm P} \rangle_\infty (r)} (J^{-1} \, {\rm W}_{n\star}^2 \, 
{\rm f}) \bigl(r,X(s) + (s-r) \omega,\xi \bigr) \ ds \, d \omega \, d \xi \\
 \displaystyle \qquad \quad \lesssim \int_r^t \int_{\mathbb S^2} \int_{\vert \xi \vert \leq \langle {\rm P} \rangle_\infty (r)} 
\frac{J^{1/p} \ \langle \xi \rangle^{1/p} \ {\rm f} \bigl(r,X(s) + (s-r) \omega,\xi \bigr)}{(s-r)^{2/p} \, \langle \xi \rangle^{2+ (1/p)} \, 
\bigl(1 + \omega \cdot \nu (\xi) \bigr)^2 \, \bigl(1 + \omega \cdot \nu \circ \Xi(s) \bigr)^{1/p}} \ ds \, d \omega \, d \xi \\
\displaystyle \qquad \quad \lesssim {\rm S}_0^{1-(1/p)} \ \text{\small $ \pmb{\mathscr{E}} $}_{\! 0}^{1/p} \  \Bigl( \int_r^t \int_{\mathbb S^2} 
\int_{\vert \xi \vert \leq \langle {\rm P} \rangle_\infty (r)} \frac{\langle \xi \rangle^{-2q-(q/p)}}{(s-r)^{2q/p}} \ \frac{\bigl(1 + 
\omega \cdot \nu \circ \Xi(s) \bigr)^{-q/p}}{\bigl(1 + \omega \cdot \nu (\xi) \bigr)^{2q}} \ ds \, d \omega \, d \xi \Bigr)^{1/q} .
 \end{array} \]
On the one hand, from Lemma \ref{Maximalloss}, we have
\[ \langle \xi \rangle^{-2q-(q/p)} \ \bigl(1 + \omega \cdot \nu (\xi) \bigr)^{-2q} \lesssim \langle \xi \rangle^{2q-(3 q/p)} \ 
\bigl(1 + \omega \cdot \nu (\xi) \bigr)^{-q/p} . \]
On the other hand, from Lemma \ref{gainaveraging} together with the condition $ 2q/p < 1 $, we can assert that
\[ \begin{array}{rl} 
\displaystyle \int_{\mathbb S^2} \bigl(1 + \omega \cdot \nu (\xi) \bigr)^{-q/p} \ \bigl(1 + \omega \cdot \nu \circ \Xi(s) 
\bigr)^{-q/p} \ d \omega \lesssim \! \! \! & \displaystyle \int_{\mathbb S^2} \bigl(1 + \omega \cdot \nu (\xi) \bigr)^{-2 q/p} 
\ d \omega \\
\ \! \! \! & \displaystyle + \int_{\mathbb S^2} \bigl(1 + \omega \cdot \nu \circ \Xi(s) \bigr)^{-2 q/p} \ d \omega \lesssim 1 . 
 \end{array} \]
This implies that 
\[ \Bigl( \int_r^t \int_{\mathbb S^2} \int_{\vert \xi \vert \leq \langle {\rm P} \rangle_\infty (r)} \frac{\langle \xi 
\rangle^{-2q-(q/p)}}{(s-r)^{2q/p}} \ \frac{\bigl(1 + \omega \cdot \nu \circ \Xi(s) \bigr)^{-q/p}}{\bigl(1 + \omega \cdot 
\nu (\xi) \bigr)^{2q}} \ ds \, d \omega \, d \xi \Bigr)^{1/q} \lesssim t^{1 - (3/p)} \ \langle {\rm P} \rangle_\infty (r)^{5-(6/p)} . \]
 From there, it is easy to deduce (\ref{atthisstage}).
 \end{prof} 
 

\subsection{Proof of Proposition \ref{theorempms}} \label{Proof of Propositiontheorempms} In this section, we work 
with bounded  values of $ {\rm S}_0 \in \RR_+^* $, which may go to $ 0 $. Let $ \iota \in ]0,1[ $. We perform the 
change of time variable 
\[ \tau := {\rm S}_0^{1-\iota} \, t \, , \qquad \langle \tilde {\rm P} \rangle_\infty ( \tau ) :=  \langle {\rm P} \rangle_\infty 
( {\rm S}_0^{\iota-1} \, \tau ) . \]
Let $ p > 3 $. For small values of $ {\rm S}_0 $, we can exploit (\ref{estimonslee0}) to replace $ \text{\small $ \pmb{\mathscr{E}} $}_{\! 0} $
by $ {\rm S}_0 $. With this in mind, we combine (\ref{control of the momentum spread}) with (\ref{thisisfor0}), (\ref{dhtodofinal}), 
(\ref{thisisforthel}) and (\ref{atthisstage}) to see that
\[ \begin{array}{rl}
\displaystyle \langle \tilde {\rm P} \rangle_\infty (\tau) \lesssim 1 + {\rm P}_0 + \! \! \! & \displaystyle {\rm S}_0^{(3/p)-1 + 
\iota \, (2-3/p)} \ \tau^{2-(3/p)} \ (1 +{\rm P}_0^2 )^{(3/2)-(2/p)} + ({\rm S}_0^\iota \, \tau + {\rm S}_0) \ \langle \tilde {\rm P} 
\rangle_\infty (\tau)^{2/3} \\
\ + \! \! \! & \displaystyle  {\rm S}_0^{(3/p)-1 + \iota \, (2-3/p)} \ \tau^{(p-3)/p} \int_0^\tau  \langle \tilde {\rm P} \rangle_\infty (r)^{3-(4/p)} \ dr \\
\ + \! \! \! & \displaystyle {\rm S}_0^{(3/2p)+ \iota \, (p-3)/(2p)} \ \tau^{(p-3)/(2p)} \ \int_0^\tau  \langle \tilde {\rm P} \rangle_\infty (r)^{4-(3/p)} \ dr .
 \end{array} \]
We can adjust $ p >3 $ sufficiently close to $ 3 $ to ensure that 
 \[ 0 < (3/p)-1 + \iota \, (2-3/p) . \]
We select $ \tau \in [0, \tau_0 ] $ with arbitrary $ \tau_0>0 $, and using Young's inequality, one can absorb the power $ 2/3 $ of 
$ \langle \tilde {\rm P} \rangle_\infty (\tau) $. There remains 
\[ \langle \tilde {\rm P} \rangle_\infty (\tau) \lesssim 1 +  \int_0^\tau  \langle \tilde {\rm P} \rangle_\infty (r)^4 \ dr . \]
This estimate is compatible with \href{https://en.wikipedia.org/wiki/Bihari-LaSalle_inequality}{Bihari-LaSalle inequality}. It furnishes a 
time $ \tilde \tau \leq \tau_0 $ and a continuous increasing function $ \tilde \cF : [0,\tilde \tau ] \rightarrow \RR_+ $ (not depending on 
$ \iota $) such that
\[ \langle \tilde {\rm P} \rangle_\infty (\tau) \leq \tilde \cF (\tau) \, , \qquad \forall \, \tau \in [0, \tilde \tau] , \]
or equivalently
\begin{equation} \label{dcujzefz}  
{\rm P} (t) \leq \langle {\rm P} \rangle_\infty (t) \leq \tilde \cF ( {\rm S}_0^{1-\iota} \, t) \, , \qquad \forall \, t \in [0,T[ \cap 
[0, {\rm S}_0^{\iota-1} \, \tilde \tau] . 
\end{equation}
Passing to the limit ($ \iota \rightarrow 0 + $), with $ c = \tilde \tau $ and $ \cF ({\rm S}_0;t) := \tilde \cF ( {\rm S}_0 \, t)  $, 
we recover (\ref{conofcpbis}) .
\hfill $ \square $


\section{Application} \label{Application}  

This section is still devoted to the study of (\ref{VlasoveqRVMD})-(\ref{maxwelleqRVMD})-(\ref{stationarystate}). 
But from now on, we assume that the initial data  $ {\rm U}_0 \equiv {\rm U}_0^\eps $ takes the form of 
(\ref{stationarystateperturb}) where $ {\rm U}_a^\eps(0,\cdot) $ is issued from a well-prepared approximate 
solution $ {\rm U}_a^\eps $ (in the sense of Definition \ref{Local well-prepared approximate solution} below). 
Thus, we consider a family of Cauchy problems indexed by $ \eps \in ]0,1] $. From \cite{zbMATH03951266}, 
for all $ \eps \in ]0,1] $, there exists locally in time, say on the interval $ \lbrack 0,T^\eps \rbrack $ with 
$ T^\eps \in \RR_+^* $, a unique smooth solution $ {\rm U} \equiv{\rm U}^\eps $. The main goal is to show 
that we can extract a lower bound $ \cT \in \RR_+^* $ such that $ \cT \leq T^\eps $ for all $ \eps \in ]0,1] $. 
To this end, Proposition \ref{theorempms} is of no use. Indeed, as soon as $ {\rm B}_a^\eps \sim 
\eps^{-1} $, the condition $ \cN ({\rm U}_0^\eps) \lesssim 1 $ is not uniformly satisfied when $ \eps \rightarrow 0 $. 
However, the proof of Proposition \ref{theorempms} does not exploit a number of specificities that can be 
detected by working  in the vicinity of $ {\rm U}_a^\eps $. Far beyond \cite{MR4084146,MR4357273,MR4244270}, 
this will allow us to incorporate the dense, hot and strongly magnetized framework.

\noindent In Subsection \ref{Preliminary material}, we furnish some preparatory material. In Subsection 
\ref{discussenergy}, we discuss the issue of energy estimates. This is an opportunity to  identify the 
challenges posed by the introduction of $ {\rm U}_a^\eps $. In Subsection \ref{prove Theorem mainmain}, 
we prove Theorem \ref{mainmaintheo}.


\subsection{Preliminary background} \label{Preliminary material}   
We start in Paragraph \ref{Well-prepared approximate solutions} by generalizing the choice (\ref{prototype}). In Paragraph 
\ref{equations for the perturbation}, we write the equations which are satisfied by the perturbation $ U = (f,E,B) $. In Paragraph 
\ref{compilitdazh}, we specify the exact meaning of \guillemotleft compatible initial data\guillemotright.

 
\subsubsection{Well-prepared approximate solutions} \label{Well-prepared approximate solutions} The system built with 
(\ref{VlasoveqRVMD}) and the two equations at the left of (\ref{maxwelleqRVMD}) is denoted by $ \cL ({\rm U} , \part) {\rm U} = 0 $.
This quasilinear system (with integral source term $ {\rm J} $) is obtained by ignoring the compatibility conditions.

\begin{defi}[Well-prepared approximate solutions] \label{Local well-prepared approximate solution} We say 
that the family $ \{ {\rm U}^\eps_a \}_\eps $ is a \underline{well}-p\underline{re}p\underline{ared} approximate solution to 
(\ref{VlasoveqRVMD})-(\ref{maxwelleqRVMD}) if:
\begin{itemize}
\item [$ \bullet $] There exists a time $ \cT_a \in \RR_+^* $ such that, for all $ \eps \in ]0,1]  $, we have
\begin{equation}\label{cTaa} 
{\rm U}^\eps_a = ({\rm f}^\eps_a,{\rm E}^\eps_a,{\rm B}^\eps_a) \in C^1 ([0,\cT_a] \times \RR^3 \times 
\RR^3;\RR^7 ) .
\end{equation}
\item [$ \bullet $] There exists a constant $ {\rm P}_a \in \RR_+^* $ such that, for all $ \eps \in ]0,1]  $, we have
 \begin{equation}\label{supppourUepsa} 
 \text{\rm supp} \ {\rm f}^\eps_a \subset [0,\cT_a] \times \RR^3 \times B (0,{\rm P}_a] . 
 \end{equation}
 \item [$ \bullet $] There exist positive constants $ {\rm S}_a $ and $ {\rm H}_a $ such that (in the sup-norm with  respect 
 to the domain $ [0,\cT_a] \times \RR^3 \times \RR^3 $), for all $ \eps \in ]0,1]  $, we have the sup-norm estimates
 \begin{subequations}\label{Lipesteforuepsa} 
 \begin{align}
 & \displaystyle \parallel {\rm f}^\eps_a \parallel_\infty \leq {\rm S}_a , \label{Linftyteforuepsaf} \\
 & \displaystyle  \parallel {\rm E}^\eps_a \parallel_\infty \leq {\rm S}_a , , \label{LinftyteforuepsaE} \\
 & \displaystyle \parallel \eps \, {\rm B}^\eps_a \parallel_{\infty} \leq {\rm H}_a . \label{LipesteforuepsaB}
 \end{align}
 \end{subequations}  
 \item [$ \bullet $] There exists a constant $ {\rm L}_a $ such that, for all $ \eps \in ]0,1]  $, we have the Lipschitz estimate
  \begin{equation}\label{Lipesteforuepsaf} 
\parallel \nabla_{t,x,\xi} {\rm f}^\eps_a \parallel_{\infty}\leq {\rm L}_a .
 \end{equation}
 \item [$ \bullet $] Let $ {\rm E}^\eps_h $ be the solution to the linear wave equation (\ref{homogeneous version})
 with initial data
 \[ {{\rm E}_h^\eps}_{\mid t=0} = {\rm E}_a^\eps (0,\cdot) \, , \qquad \part_t {{\rm E}_h^\eps}_{\mid t=0} = \nabla_x 
 \times {\rm B}_a^\eps (0,\cdot) . \]
There exists a positive constant $ {\rm S}_h \in \RR_+^* $ such that, for all $ \eps \in ]0,1]  $, we have
\begin{equation} \label{Linftyteforuepsah}
\parallel {\rm E}^\eps_h \parallel_{L^1 ([0,\cT_a] ; L^\infty_x)} \leq {\rm S}_h .
\end{equation}
\item [$ \bullet $] Introduce the remainder $ {\rm R}^\eps_a := \cL ({\rm U}^\eps_a , \part) {\rm U}^\eps_a $. There 
exist positive constants $ {\rm S}_r $ and $ {\rm N}_r $ such that, for all $ \eps \in ]0,1] $, we have 
\begin{subequations}\label{remainderupesa} 
\begin{align}
 & \displaystyle 
\parallel {\rm R}^\eps_a \parallel_{\infty} \leq {\rm S}_r , \label{remainderupesa1} \\
 & \displaystyle \parallel {\rm R}^\eps_a \parallel_{L^2} \leq {\rm N}_r \, \eps , \label{remainderupesa2} 
  \end{align}
 \end{subequations} 
 where $ \parallel {\rm R}^\eps_a \parallel_{L^2} $ stands for the following $ L^2 $-norm:
\[ \parallel {\rm R}^\eps_a \parallel_{L^2} \, :=  \bigl( \parallel {\rm R}_{{\rm f} a}^\eps \parallel^2_{L^2
( [0,\cT_a] \times \RR^3 \times \RR^3)} + \parallel {\rm R}_{{\rm E} a}^\eps \parallel^2_{L^2 ( [0,\cT_a] 
\times \RR^3)} + \parallel {\rm R}_{{\rm B} a}^\eps \parallel^2_{L^2 ( [0,\cT_a] \times \RR^3)} \bigr)^{1/2} . \]
   \item [$ \bullet $] For all $ (\eps,x) \in ]0,1] \times \RR^3 $, at the time $ t = 0 $, we impose the compatibility conditions
   \begin{subequations}\label{compatibility conditionforUepsa2} 
\begin{align}
 & \displaystyle  \nabla_x \cdot {\rm E}^\eps_a (0,x) =  \int_{\RR^3} {\rm f}^\eps_a (0,x,\xi) \, d \xi - \rho^\eps (x) , 
 \label{compatibility conditionforUepsa21} \\
 & \displaystyle \nabla_x \cdot {\rm B}^\eps_a (0,x) = 0 . \label{compatibility conditionforUepsa22} 
  \end{align}
 \end{subequations} 
  \end{itemize}
\end{defi}

\noindent Let us come back to the situation (\ref{prototype}), where $ {\rm B}_e $ and $ {\rm M}^\eps $ are adjusted 
as in (\ref{conditionforthebe}) and (\ref{conditionforthedenM}). This means to deal with $ \tilde {\rm U}^\eps_a $, which 
gives rise to $ \tilde {\rm E}^\eps_h $ and $ \tilde {\rm R}^\eps_a $. The four conditions (\ref{cTaa}), (\ref{supppourUepsa}), 
(\ref{Lipesteforuepsa}) and (\ref{Lipesteforuepsaf}) are clearly satisfied for all $ \cT_a \in \RR_+^* $. Due to the last condition 
inside (\ref{conditionforthebe}), we simply find that $ \tilde {\rm E}^\eps_h \equiv 0 $ so that (\ref{Linftyteforuepsah}) is evident. 
In the same vein, we have $ \tilde {\rm R}^\eps_a := \cL (\tilde {\rm U}^\eps_a , \part) \tilde {\rm U}^\eps_a \equiv 0 $ so that 
(\ref{remainderupesa}) is achieved. The condition (\ref{compatibility conditionforUepsa21}) for $ \tilde {\rm E}^\eps_a \equiv 0 $ 
can be guaranteed by adjusting $ \rho^\eps $ as indicated in (\ref{basiccompatibility conditionsrm}), while the condition 
(\ref{compatibility conditionforUepsa22}) for $ \tilde {\rm B}^\eps_a = \eps^{-1} \, {\rm B}_e $ is a consequence of the second 
condition inside (\ref{conditionforthebe}). In brief, {\it the family $ \{ \tilde {\rm U}^\eps_a \}_\eps $ is a well-prepared approximate 
solution.}

\smallskip

\noindent For a better understanding, Definition \ref{Local well-prepared approximate solution} must be supplemented 
with a number of remarks. Indeed, the lines (\ref{cTaa}), (\ref{supppourUepsa}) and (\ref{Linftyteforuepsaf}) are just 
extensions of (\ref{cestpourf0}). But the other constraints  serve other purposes which must be clarified. We discuss 
the rest of (\ref{LinftyteforuepsaE})-(\ref{Lipesteforuepsaf}), (\ref{LipesteforuepsaB}), (\ref{Linftyteforuepsah}), 
(\ref{remainderupesa}) and (\ref{compatibility conditionforUepsa2}) in separate Paragraphs, respectively in {\it ra)}, 
{\it rb)}, {\it rc)}, {\it rd)} and {\it re)} below. 

\smallskip

\noindent $ \circ $ {\it ra) Link between the notion of well-prepared data and the Lipschitz estimate on 
$ {\rm f}^\eps_a $}. By definition, we have the decomposition $ {\rm R}^\eps_a = ({\rm R}_{{\rm f} a}^\eps,
{\rm R}_{{\rm E} a}^\eps, {\rm R}_{{\rm B} a}^\eps) \in \RR \times \RR^3 \times \RR^3 $, with in particular
\begin{equation} \label{aidercomprehension}
{\rm R}_{{\rm f} a}^\eps := \partial_t {\rm f}^\eps_a + \nu (\xi) \cdot \nabla_x {\rm f}^\eps_a  + {\rm E}^\eps_a  \cdot 
\nabla_\xi {\rm f}^\eps_a + \bigl( \nu(\xi) \times {\rm B}^\eps_a \bigr) \cdot \nabla_\xi {\rm f}^\eps_a .
\end{equation}
Then, from (\ref{supppourUepsa}), it is easy to see that
 \begin{equation} \label{remainderupesasupp}
 \text{\rm supp} \ {\rm R}_{{\rm f} a}^\eps \subset [0,\cT_a] \times \RR^3 \times B (0,{\rm P}_a] . 
 \end{equation}
Then, from (\ref{LinftyteforuepsaE}), (\ref{Lipesteforuepsaf}) and (\ref{remainderupesa1}), we can also infer that
\begin{equation} \label{wellwellprepa}
\parallel \bigl( \nu(\xi) \times {\rm B}^\eps_a \bigr) \cdot \nabla_\xi {\rm f}^\eps_a \parallel_\infty \leq {\rm S}_r
+ (2+ {\rm S}_a) \ {\rm L}_a.
\end{equation}
This is not a consequence of (\ref{LipesteforuepsaB}) and (\ref{Lipesteforuepsaf}). This is possible only if 
$ {\rm B}^\eps_a $ and $ {\rm f}^\eps_a $ are adjusted accordingly. In particular, this should hold true at the 
time $ t = 0 $.  This is the common notion of {\it well-prepared initial data}, see for instance \cite{MR4084146,MR4357273} 
or the books \cite{MR2562165,Rauch}. In coherence with (\ref{Lipesteforuepsaf}), where the $ {\rm L} $ inside 
$ {\rm L}_a $ is for \underline{L}ipschitz, this prevents the emergence of rapid variations inside $ {\rm f}^\eps_a $. 
Instead, oscillations can be (or are intended to be) introduced through $ f^\eps_0 $. \hfill $ \circ $

\smallskip

\noindent $ \circ $ {\it rb) Large magnetic fields}.  Strongly magnetized plasmas are available due to the  weight 
$ \eps $ which is put in factor of $ {\rm B}^\eps_a $ inside (\ref{LipesteforuepsaB}), where the letter $ {\rm H} $ 
inside $ {\rm H}_a $ is for \underline{H}igh. Again, the prototype of such behavior is $ \eps^{-1} \, {\rm B}_e $
inside $ \tilde {\rm U}^\eps_a $. From there, the construction of more elaborate approximate solutions is a subject 
in its own right. This is a way of revealing various physical phenomena. In the absence of coupling, if we concentrate 
only on  the Vlasov part, this involves a WKB analysis which ties in with the recent advances 
\cite{MR3582249,gyrokineticCF}.  \hfill $ \circ $

\smallskip

\noindent $ \circ $ {\it rc) The impact of initial electromagnetic fields}.  The condition (\ref{stationarystateperturb}) may 
seem restrictive since we impose $ (E^\eps_0, B^\eps_0) \equiv 0 $. This is to forget that the presence at the time $ t = 0 $ 
of a non-zero electromagnetic component is considered as an integral part of the construction of $ {\rm U}^\eps_a $. 
Indeed, the access to $ {\rm U}^\eps_a $ may be achieved in two steps. The first is to solve the Cauchy problem 
(\ref{homogeneous version}) with initial data $ ({\rm E}_0, {\rm B}_0) \equiv ( {\rm E}^\eps_a, {\rm B}^\eps_a) 
(0,\cdot) $ in order to get $ ({\rm E}_h^\eps , {\rm B}_h^\eps) $ ; the second is to seek approximate solutions in the 
form $  {\rm U}^\eps_a = ( {\rm f}^\eps_a, {\rm E}_h^\eps +  \check {\rm E}^\eps_a, {\rm B}_h^\eps + \check 
{\rm B}^\eps_a) $ with $ (\check {\rm E}^\eps_a, \check {\rm B}^\eps_a) (0,\cdot) \equiv 0 $. When doing this,
there is not total freedom concerning the choice of $ ( {\rm E}^\eps_a, {\rm B}^\eps_a) (0,\cdot) $. Indeed, it is 
necessary to ensure the property (\ref{Linftyteforuepsah}). Now, for reasons explained in Paragraph \ref{Control ofDh}, 
this cannot be taken for granted under only $ L^2 $ or even $ H^1 $-controls on $ ( {\rm E}^\eps_a, {\rm B}^\eps_a) 
(0,\cdot) $.   \hfill $ \circ $

\smallskip

\noindent $ \circ $ {\it rd) About the remainder}. The condition (\ref{remainderupesa}) gives meaning to the word 
\guillemotleft approximate\guillemotright. The smallness of the remainder is measured by the relatively mild
constraint (\ref{remainderupesa2}). \hfill $ \circ $

\smallskip

\noindent $ \circ $ {\it re) About the compatibility conditions}. In view of (\ref{compatibility conditionsrm}), the restrictions 
(\ref{compatibility conditionforUepsa21}) and (\ref{compatibility conditionforUepsa22}) seem unavoidable. But it is difficult
and pointless to preserve the compatibility conditions when constructing approximate solutions $ {\rm U}^\eps_a $. It is 
preferable to rather leave some freedom to the remainder $ {\rm R}^\eps_a $ and therefore to the choice of $ {\rm U}^\eps_a $ 
while, as will be seen, the exact solution $ {\rm U}^\eps $ is sure to propagate the initial compatibility conditions. \hfill $ \circ $

\smallskip

\noindent Changing $ T^\eps $ if necessary, we can always assume that we work with $ T^\eps \leq \cT_a $.

  
\subsubsection{The equations for the perturbation} \label{equations for the perturbation} Recall the decomposition (\ref{typeofsolution})
where $ {\rm U} $ has been put in the form $ {\rm U} = {\rm U}^\eps_a + U $ with $ U = (f,E,B) $. The Lorentz force $ {\rm F} $ acts on 
the charged particles through a contribution coming from $ {\rm U}^\eps_a $ and a part issued from the self-consistent electromagnetic 
field $ (E,B) $. We have $ {\rm F} = {\rm F}^\eps_a + F $ with
\begin{equation}\label{decompofexpa} 
{\rm F}^\eps_a := {\rm E}^\eps_a (t,x) + \nu(\xi) \times {\rm B}^\eps_a (t,x) \, , \qquad F^\eps \equiv F :=  E (t,x) + \nu(\xi) \times B (t,x) . 
\end{equation}
After substitution of $ {\rm U} $ as in (\ref{typeofsolution}) inside (\ref{VlasoveqRVMD})-(\ref{maxwelleqRVMD}), the Vlasov equation becomes
\begin{equation}\label{Vlasoveq} 
\partial_t f + \, \nu (\xi) \cdot \nabla_x f + ({\rm F}^\eps_a + F) \cdot \nabla_\xi f  + F \cdot \nabla_\xi 
f^\eps_a  = {\rm R}_{{\rm f} a}^\eps ,
\end{equation}
and the Maxwell's equations reduce to
 \begin{subequations}\label{maxwelleq} 
 \begin{align}
 & \displaystyle \partial_t {E} - \nabla_x \times {B} + \int_{\RR^3} \nu(\xi) \ f(t,x,\xi) \, d \xi = {\rm R}_{{\rm E} a}^\eps  , \label{maxwelleq1} \\
 & \displaystyle \partial_t {B} + \nabla_x \times {E} = {\rm R}_{{\rm B} a}^\eps . \label{maxwelleq2} 
 \end{align}
 \end{subequations}  
 
\noindent The system (\ref{Vlasoveq})-(\ref{maxwelleq}) is denoted by $ \cL^\eps_a (U , \part) U = 0 $. 

  
\subsubsection{Compatible initial data}\label{compilitdazh}  In view of (\ref{stationarystateperturb}), the initial data are adjusted 
according to
 \begin{equation}\label{initial data} 
f_{\mid t=0} = f^\eps_0 , \qquad E_{\mid t=0} = 0 , \qquad B_{\mid t=0} = 0 .
\end{equation}
Constraints inherited from (\ref{cestpourf0}), (\ref{compatibility conditionsrm})-(\ref{compatibility conditionforUepsa2}) 
and (\ref{LipesteforuepsaB}) must be imposed on $ f^\eps_0 $. 

\begin{defi}\label{assumptiononofoeps}[Compatible family of initial data] We say that the family $ \{ f_0^\eps \}_\eps $ 
is \underline {com}p\underline {atible} when the following five conditions are satisfied:
 \begin{itemize}
\item [a)] Regularity:
\begin{equation}\label{Regularity} 
f_0^\eps \in C^1_c (\RR^3 \times \RR^3;\RR) \, , \qquad \forall \, \eps \in ]0,1] .
\end{equation}
\item [b)] Uniform control on the the size of the momentum support:
\begin{equation}\label{minimalradiusforf} 
\text{supp} \ f_0^\eps (\cdot) \subset \RR^3 \times B(0,{\rm P}_0] \, , \qquad \forall \, \eps \in ]0,1] . 
\end{equation}  
 \item [c)] Uniform sup-norm estimate:
\begin{equation}\label{Lebesgue bounds}
\qquad \parallel f_0^\eps \parallel_\infty := \! \! \sup_{(x,\xi) \in \RR^3 \times \RR^3} \vert f_0^\eps \vert 
(x,\xi) \leq {\rm S}_0 , \qquad \forall \, \eps \in ]0,1].
\end{equation}
 \item [d)] $ L^2 $-smallness:
 \begin{equation}\label{absorptioforsm} 
\parallel f_0^\eps \parallel_{L^2} \lesssim {\rm N}_0 \, \eps  \, , \qquad \forall \, \eps \in ]0,1] .
\end{equation}
 \item [e)] Compatibility conditions (issued from $ {\rm U}^\eps_a $):
 \begin{equation}\label{compatibility conditions} 
\int_{\RR^3} {\rm f}^\eps_0 (0,x,\xi) \ d \xi = 0 \, , \qquad \forall \, (x,\eps) \in \RR^3 \times ]0,1].
 \end{equation} 
  \end{itemize}
\end{defi}

\noindent In view of (\ref{stationarystateperturb}), (\ref{compatibility conditionforUepsa2}) and (\ref{compatibility conditions}),
the compatibility conditions (\ref{compatibility conditionsrm}) are satisfied at time $ t = 0 $. They are propagated so that
 \begin{subequations}\label{compaeqt} 
 \begin{align}
 & \displaystyle  \nabla_x \cdot {\rm E}^\eps_a (t,x) + \nabla_x \cdot E (t,x) = \int_{\RR^3} {\rm f}^\eps_a (t,x,\xi) \ d \xi + 
 \int_{\RR^3} f (t,x,\xi) \ d \xi - \rho^\eps (x) , \label{compaeqt1} \\
 & \displaystyle \nabla_x \cdot {\rm B}^\eps_a (t,x) +\nabla_x \cdot B (t,x) = 0 . \label{compaeqt2}  
  \end{align}
 \end{subequations}

\noindent We can test Definition \ref{assumptiononofoeps} with $ f_0^\eps \equiv \eps \, f_0 $ where $ f_0 \in C^1_c (\RR^3 
\times \RR^3) $ does not depend on $ \eps $ and is as indicated in (\ref{cestpourf0}) and (\ref{compatibility conditionsrmf0})
with $ \rho \equiv 0 $. We have obviously a), b) and c) as well as e). Moreover, we have d) with $ {\rm N}_0 = \parallel f_0 
\parallel_{L^2} $. Thus, under these assumptions on $ f_0 $, we can assert that:

\smallskip

\centerline{\it The family $ \{ \eps \, f_0 \}_\eps $ is compatible.}

\smallskip

\noindent Now, we can apply Theorem \ref{mainmaintheo} with $ {\rm U}^\eps_a \equiv 0 $ and $ f^\eps_0 = \eps \, f_0 $.
Then, just fix $ \eps = 1 $ to recover the conclusions of Theorem \ref{inimaintheo}. Theorem \ref{mainmaintheo} can 
therefore be viewed as an extension of Theorem \ref{inimaintheo}. Its proof implies a few subtleties in comparison 
to what has been done in Section \ref{RadonFourier}.

 
 \subsection{Estimating the pertubation} \label{discussenergy} In the perturbative context (\ref{Vlasoveq})-(\ref{maxwelleq}), it 
 becomes more complicated to control $ f $ and $ (E,B) $. In Subsection \ref{oscillating characteristic flow}, we comment on the structure of
 the characteristic flow and, from there, we deduce sup-norm estimates on $ f (t,\cdot) $. In Subsection \ref{totalenergyofRVMFsystem}, 
 we explain why $ \text{\small $ \pmb{\mathscr{E}} $} $ is of no use any more, and we follow an alternative path in order to manage the 
 $ L^2 $-norm of $ (E,B) (t,\cdot) $.
 
 
\subsubsection{The oscillating characteristic flow} \label{oscillating characteristic flow} In the context of (\ref{typeofsolution}), 
the differential system (\ref{diffeoom}) breaks down into 
\renewcommand\arraystretch{2}
\begin{equation} \label{diffeoomeps}
\left \lbrace \begin{array} {ll}
\displaystyle \frac{d X}{dt} (t,y,\eta) = \nu (\Xi) , \quad & X(0,y,\eta) = y , \\
\displaystyle \frac{d \Xi}{dt} (t,y,\eta) =  {\rm F}^\eps_a (t,X,\Xi) + F(t,X,\Xi) \ , \quad & \Xi(0,y,\eta) 
= \eta .
\end{array} \right. 
\end{equation}
\renewcommand\arraystretch{1}

\noindent In a first attempt, the impact of the self-consistent electromagnetic field $ (E,B) $ may be neglected. In 
particular, in the case of strongly magnetized plasmas (when $ {\rm B}^\eps_a \sim \eps^{-1} $), the dominant part 
inside (\ref{diffeoomeps}) may be thought as 
\renewcommand\arraystretch{2}
\begin{equation} \label{diffeoomepsa}
\left \lbrace \begin{array} {ll}
\displaystyle \frac{d X^\eps_a}{dt} (t,y,\eta) = \nu (\Xi^\eps_a) , \quad & X^\eps_a (0,y,\eta) = y , \\
\displaystyle \frac{d \Xi^\eps_a}{dt} (t,y,\eta) =  {\rm F}^\eps_a (t,X^\eps_a,\Xi^\eps_a) \ , \quad & \Xi^\eps_a (0,y,\eta) 
= \eta .
\end{array} \right. 
\end{equation}
\renewcommand\arraystretch{1}

\noindent The study of (\ref{diffeoomepsa}), when for instance $ {\rm E}^\eps_a = {\rm E}_a $ and $ {\rm B}^\eps_a 
= {\rm B}_a + \eps^{-1} \, {\rm B}_e $ with a given field $ ({\rm E}_a,{\rm B}_a) $ not depending on $ \eps $, is very 
informative. This provides WKB expansions for $ (X^\eps_a,\Xi^\eps_a) $ which reveal  the high complexity of the 
underlying motions \cite{MR3582249,gyrokineticCF}. The flow associated with (\ref{diffeoomepsa}) is 
strongly oscillating in both time, space and momentum; the directions $ \Xi^\eps_a / \vert \Xi^\eps_a \vert $ are 
rapidly oscillating, but not $ \langle \Xi^\eps_a \rangle $. Let us consider the frozen version $ {\rm F} {\rm D}^\eps_a $
of $ {\rm D} $, which is
\[ {\rm F} {\rm D}^\eps_a (t,y,\eta) := \int_0^t \nu \circ \Xi^\eps_a (s,y,\eta) \cdot {\rm E}^\eps_a \bigl( s,X^\eps_a
(s,y,\eta)\bigr) \ ds . \]
From (\ref{LinftyteforuepsaE}), we can see that $ \vert {\rm F} {\rm D}^\eps_a \vert \leq t \, {\rm S}_a $. This is 
already an indication that $ {\rm D} $ should remain under control  only where the impact of $ U $ can be taken 
into account. To this end, by analogy with (\ref{maximalsizesm}) and (\ref{ball of radius}), we can define
\[ P (t) := \inf \ \bigl \lbrace R \in \RR_+ \, ; \, f (t,x,\xi) = 0 \ \text{for all} \  x \in \RR^3 \ \text{and for all} \ \xi \in \RR^3 \ \text{with} \ 
R \leq \vert \xi \vert \bigr \rbrace , \]
as well as
\begin{equation}\label{cestpourPinfty}
\langle P \rangle_\infty (t) := \sup_{s \in [0,t]} \ \bigl( 1 + P(s)^2 \bigr)^{1/2} . 
\end{equation}
The solutions to the Vlasov equation (\ref{VlasoveqRVMD})  are constant along the characteristics, so that
\[ 0 \leq {\rm f}^\eps_a (t,x,\xi) + f^\eps (t,x,\xi) = {\rm f}^\eps_a \bigl( 0, X(-t,x,\xi) , \Xi(-t,x,\xi) \bigr)  + f^\eps_0 
\bigl( X(-t,x,\xi) , \Xi(-t,x,\xi) \bigr) . \]
As a consequence, we can assert that
\begin{equation} \label{supnormonf}
\parallel f^\eps (t,\cdot) \parallel_\infty \lesssim {\rm S}_a + {\rm S}_0 \, , \qquad \forall (t,\eps) \in [0,T^\eps] \times ]0,1] .
\end{equation}

 
\subsubsection{Propagation of the $ L^2 $-norm} \label{totalenergyofRVMFsystem} The system (\ref{VlasoveqRVMD})-(\ref{maxwelleqRVMD})
is endowed with the conserved total energy $ \text{\small $ \pmb{\mathscr{E}} $} (t) $, see (\ref{consertotenRVMgen})-(\ref{consertotenRVM}). 
In the strongly magnetized case, we find that  $ \text{\small $ \pmb{\mathscr{E}} $} (t) = \text{\small $ \pmb{\mathscr{E}} $} (0) \sim \eps^{-1} $
and the use of $ \text{\small $ \pmb{\mathscr{E}} $} (t) $ does not help in the perspective of uniform estimates. Instead, we could consider the 
{\it relative energy} $ \text{\small $ \mathscr{E} $} (t) $ defined by
\[ \text{\small $ \mathscr{E} $} (t) := \int_{\RR^3} \int_{\RR^3} \langle\xi\rangle \ {f(t,x,\xi)} \ dx \, d\xi + \frac12 \int_{\RR^3} 
(|E|^2 + |B|^2)(t,x) \ dx . \]
Let us first examine what happens in the case of $ \tilde {\rm U}_a^\eps $, with $ \tilde {\rm U}_a^\eps $ as in (\ref{prototype}).

\begin{lem} \label{consebutspeci} Assume that $ {\rm U}_a^\eps \equiv \tilde {\rm U}_a^\eps $. Then, for all $ t \in [0,T^\eps] $, 
we have $ \text{\small $ \mathscr{E} $} (t) = \text{\small $ {\mathscr{E}} $}(0) $.
 \end{lem}

\begin{prof} The proof of this conservation of $ \text{\small $ \mathscr{E} $} $ may be achieved through direct computations 
based on (\ref{Vlasoveq})-(\ref{maxwelleq}). But it is more instructive to deduce it from the well-known conservation \cite{Rhein} 
of the total energy $ \text{\small $ \pmb{\mathscr{E}} $} (t) $. When $ {\rm U}_a^\eps \equiv \tilde {\rm U}_a^\eps $, this gives rise to
\[ \begin{array}{l}
\displaystyle \int_{\RR^3} \int_{\RR^3} \langle \xi \rangle \ ({\rm M} + f) (t,x,\xi) \ dx \, d\xi + \frac12 \int_{\RR^3} |E(t,x)|^2 \ dx + 
\frac12 \int_{\RR^3} |\eps^{-1} \, {\rm B}_e (x) + B(t,x)|^2 \ dx \\
\displaystyle \qquad = \int_{\RR^3} \int_{\RR^3} \langle \xi \rangle \ ({\rm M} + f_0) (x,\xi) \ dx \, d\xi + \frac12 \int_{\RR^3} 
|\eps^{-1} \, {\rm B}_e (x) |^2 \ dx .
\end{array} \]
This is the same as
\[ \text{\small $ \mathscr{E} $} (t) + \eps^{-1} \, \int_{\RR^3} {\rm B}_e (x) \cdot B (t,x) = \text{\small $ \mathscr{E} $} (0) . \]
Since $ {\rm B}_e $ is an irrotational field, it can be written as the gradient of a magnetic potential $ {\rm A}_e $. After 
integration by parts, this yields
\begin{equation}\label{cancelllllll}  
\int_{\RR^3} {\rm B}_e (x) \cdot B(t,x) \ dx = - \int_{\RR^3} {\rm A}_e (x) \ \nabla_x \cdot B(t,x) \ dx . 
 \end{equation}
But, for $ {\rm B}_a^\eps \equiv \eps^{-1} \, {\rm B}_e $ with $ {\rm B}_e $ as in (\ref{conditionforthebe}), the Gauss's law 
(\ref{compaeqt2}) reduces to $ \nabla_x \cdot B \equiv 0 $, which leads to the expected result.
\end{prof}

\noindent To prove that $ 0 < \cT \leq T^\eps $, it seems essential to compensate the large factor $ \eps^{-1} $ which may 
occur inside (\ref{Vlasoveq}) in front of $ \nabla_\xi f $. To this end, the idea is to obtain some smallness information on $ U $. 
As already mentioned, this cannot be achieved by passing through $ \text{\small $ \pmb{\mathscr{E}} $} $ (as in Section 
\ref{RadonFourier}). This can also not be obtained by using $ \text{\small $ \mathscr{E} $} $. Let us briefly explain why.
In view of (\ref{initial data}), assuming instead of (\ref{absorptioforsm}) that $ \parallel f_0^\eps \parallel_{L^1} \lesssim 
\eps $, Lemma \ref{consebutspeci} would imply that $ \text{\small $ \mathscr{E} $} (t) \lesssim \eps $. But resorting 
to $ \text{\small $ \mathscr{E} $} $ is inadequate. There are two reasons for this:
\begin{itemize}
\item [$ - $] the quantity $ \text{\small $ \mathscr{E} $} $ is not conserved when dealing in a wider context than $ \tilde {\rm U}_a^\eps $,
with $ {\rm U}_a^\eps $ as in Definition \ref{Local well-prepared approximate solution}. This is due (in particular) to the influence of the 
remainder $ {\rm R}_a^\eps $. In fact, the situation is even more problematic: under the relaxed condition (\ref{compaeqt2}), the error 
term may be of size $ \eps^{-1} $ because (\ref{cancelllllll}) does not hold (with $ {\rm B}_e $ replaced by $ {\rm B}_a^\eps $);
\item [$ - $] the expression $ \text{\small $ \mathscr{E} $} $  is not exploitable because $ f $ represents a perturbation and, 
as such (in contrast to $ {\rm f} $), it is not necessarily sign definite. 
\end{itemize}

\noindent Instead of looking at $ \text{\small $ \mathscr{E} $} $, we perform usual energy estimates at the level of 
the system (\ref{Vlasoveq})-(\ref{maxwelleq}). With this in mind, we introduce the square of the $ L^2 $-energy of $ U(t) $, that is
\[ \EE (t) := \parallel U(t,\cdot) \parallel^2_{L^2} = \int_{\RR^3} \int_{\RR^3} f(t,x,\xi)^2 \ dx \, d\xi + \int_{\RR^3}|E(t,x)|^2 \ dx + 
\int_{\RR^3} |B(t,x)|^2 \ dx . \]

\begin{lem}[$ L^2 $-estimate on the perturbation $ U $]\label{propagationL2energy} For all $ t \leq T^\eps $, we have
\begin{equation}\label{propaL2enercons}  
 \EE (t) \leq \eps^2 \ ( {\rm N}_0^2 + {\rm N}_r^2 ) \ \exp \ \Bigl( t + 4 \, ({\rm L}_a + 1) \, \int_0^t \langle P \rangle_\infty (s)^{3/2} \, ds \Bigr) .
 \end{equation}
 \end{lem}

\begin{prof} Multiply (\ref{Vlasoveq}) by $ 2 f $ and integrate on the phase space; multiply (\ref{maxwelleq1}) and (\ref{maxwelleq2}) 
by respectively $ 2 E $ and $ 2 B $ and integrate with respect to $ dx $. By this way, we find that
\[ \partial_t \EE (t) = - 2 \, \int_{\RR^3} \int_{\RR^3} ( f \ F \cdot \nabla_\xi {\rm f}^\eps_a + E \cdot \nu (\xi) \ f) \ dx \, d \xi + 2 \, 
\langle U , {\rm R}_a^\eps \rangle_{L^2 \times L^2} . \]
From (\ref{Lipesteforuepsaf}) and Cauchy-Schwarz inequality, we have
\renewcommand\arraystretch{2}
\[ \begin{array}{l}
\displaystyle \vert \int_{\RR^3} \int_{\RR^3} ( f \ F \cdot \nabla_\xi {\rm f}^\eps_a + E \cdot \nu (\xi) \ f) \ dx \, d \xi \vert \\
\displaystyle \qquad \qquad \leq ({\rm L}_a+1) \, \int_{\RR^3} \int_{\RR^3} ( \vert E \vert + \vert B \vert ) (s,x) \ \vert f(s,x,\xi) \vert \ dx \, d \xi \\
\displaystyle \qquad \qquad \leq \sqrt 2 \ ({\rm L}_a+1) \, \Bigl( \int_{\vert \xi \vert \leq P(s)} d \xi \Bigr)^{1/2} \ \Bigl( \int_{\RR^3} ( \vert E \vert^2 + \vert B \vert^2 ) 
(s,x) \ dx \Bigr)^{1/2} \ \EE (s)^{1/2} \\
\displaystyle \qquad \qquad \leq 4 \ ({\rm L}_a+1) \ P(s)^{3/2} \ \EE (s) .
\end{array} \]
\renewcommand\arraystretch{1}

\noindent After integration in time, there remains
\[ \EE (t) \leq \EE (0) + \int_0^t \parallel {\rm R}^\eps_a (s,\cdot) \parallel^2_{L^2} \ ds + \int_0^t \bigl( 1 + 4 \, ({\rm L}_a +1) \,  
\langle P \rangle_\infty (s)^{3/2} \bigr) \ \EE (s) \ ds . \]
Now, it suffices to apply (\ref{remainderupesa2}) and (\ref{absorptioforsm}) together with Gr\"onwall's inequality.
\end{prof}

\noindent Starting with $ f_0 $ as in (\ref{minimalradiusforf}), we have 
\begin{equation}\label{cestpourI0}  
\langle P \rangle_\infty (0) \leq 1 + {\rm P}_0 \leq {\rm I}_0 \, , \qquad {\rm I}_0 := 1 +  \max \ ({\rm P}_0 , {\rm P}_a ) 
+ {\rm S}_h .
 \end{equation}
 Recall that $ T^\eps \leq \cT_a $. For $ T^\eps $ small enough, by continuity, we can always assume that 
\begin{equation}\label{startingfrombeg}  
\langle P \rangle_\infty (t) \leq 4 \, {\rm I}_0 \, , \qquad \forall \, t \in [0,T^\eps] .
  \end{equation}
Then, applying (\ref{propaL2enercons}) with $ {\rm N} := ({\rm N}_0^2 + {\rm N}_r^2 )^{1/2} $, restricting $ T^\eps \leq \cT_a $ 
again if necessary, we find 
\[ \EE (t) \lesssim \eps^2 \ {\rm N}^2 \ \exp \ \bigl( 1 + 32 \, ({\rm L}_a +1) \, {\rm I}_0^{3/2}\, \bigr) \ t  \, , 
 \qquad \forall \, t \in [0,T^\eps] . \]
  Adjust $ \tilde \cT \in ]0,1] $ with $ 0 < \tilde \cT \leq \cT_a $ small enough to be sure that
  \[ \exp \ \bigl( 1 + 32 \,  ({\rm L}_a +1) \, {\rm I}_0^{3/2}\, \bigr) \ \tilde \cT \leq 2 .\]
  We define $ T^\eps \in  ]0,1] $ as the maximal time less than $ \tilde \cT $ leading to (\ref{startingfrombeg}). 
  By construction
  \begin{equation}\label{propaL2final}  
 \EE (t) \lesssim 2 \ \eps^2 \ {\rm N}^2 \, , \qquad \forall \, t \in [0,T^\eps] \, , \qquad T^\eps \leq \tilde \cT .
   \end{equation}
The quantity $ \EE $ remains under control (and small) as long as $ \langle P \rangle_\infty $ is bounded. Remark that 
the proof of Lemma \ref{propagationL2energy} is simple and that its conclusion (\ref{propaL2final}) is not surprising. What 
is remarkable is the fact that (\ref{propaL2final}) is sufficient. Indeed, we do not need to involve costly $ H^s_\eps $-estimates 
with $ s > 0 $ large, as is typical in nonlinear geometric optics (see for instance Chapter 4 in 
\cite{MR2562165}). 

 
\subsection{Proof of Theorem \ref{mainmaintheo}} \label{prove Theorem mainmain} We start by interpreting the representation 
formula in terms of the two parts $ {\rm U}^\eps_a $ and $ U $ of $ {\rm U} $, see (\ref{typeofsolution}). We still have $ {\rm D} 
=  {\rm D}_0 + {\rm D}_h + {\rm D}_l +  {\rm D}_n $ together with (\ref{cD}) where $ {\rm f}_0 = {\rm f}^\eps_a (0,\cdot) + f^\eps_0 $, 
where the electric field $ {\rm E}_h $ is given by (\ref{homogeneous version}) with $ {\rm E}_0 $ and $ {\rm B}_0 $ as in Definition
\ref{Local well-prepared approximate solution}, whereas $ {\rm f} = {\rm f}^\eps_a + f $ and $ {\rm F} = {\rm F}^\eps_a + F $. 
Now, we can split the $ {\rm D}_\star $ into parts $ {\rm D}_{\star a}^\eps $ (with subscript \og $ a $ \fg{} for \underline{\it a}{\it pproximate})
coming from $ {\rm U}^\eps_a $, parts $ D_\star $ 
issued from $ U^\eps $, as well as bilinear terms. With this in mind, we look at $ {\rm D}_n $ as a bilinear product $ \cB ({\rm F} , 
{\rm f}) $ in terms of $ {\rm F} $ (or $ {\rm E} $ and $ {\rm B} $) and $ {\rm f}$. We denote by:
\begin{itemize}
\item [$ - $] $ {\rm D}^\eps_{0a} $ and $ D_0 $ the expression $ {\rm D}_0 $ where $ {\rm f}_0 $ is replaced respectively by
$ {\rm f}^\eps_a (0,\cdot) $ and $ f^\eps_0 $;
\item [$ - $] $ {\rm D}^\eps_{ha} $ the expression $ {\rm D}_h $ where $ {\rm E}_h $ stands for $ {\rm E}_h^\eps $ as in Definition
\ref{Local well-prepared approximate solution};
\item [$ - $] $ {\rm D}^\eps_{la} $ and $ D_l $ the expression $ {\rm D}_l $ where $ {\rm f}_0 $ is replaced respectively by
$ {\rm f}^\eps_a  $ and $ f^\eps $;
\item [$ - $] $ {\rm D}^\eps_{na} := \cB ({\rm F}^\eps_a,{\rm f}^\eps_a) $ and $ D_n := \cB (F^\eps,f^\eps) $ the expression
$ {\rm D}_n $ where the couple $ ({\rm F},{\rm f})  $ is replaced respectively by $ ({\rm F}^\eps_a , {\rm f}^\eps_a) $ and 
$ (F^\eps ,f^\eps) $.
\end{itemize}

\noindent With the above conventions, the contributions to the momentum increment brought by the approximate solution 
$ {\rm U}^\eps_{0a} $ and the perturbation $ U^\eps $ are 
\begin{equation}\label{broughtby}  
{\rm D}^\eps_a := {\rm D}^\eps_{0a} + {\rm D}^\eps_{ha} + {\rm D}^\eps_{la} +  {\rm D}^\eps_{na} \, , \qquad D :=  D_0 + D_l + D_n . 
\end{equation}
But we have also to take into account the effect of cross terms, so that
\begin{equation}\label{broughtbypoursuit} 
 {\rm D} =  {\rm D}^\eps_a + \cB(F,{\rm f}^\eps_a) + \cB({\rm F}^\eps_a,f) + D . 
 \end{equation}
We will estimate separately $ {\rm D}^\eps_a $
(in Paragraph \ref{approximateincrement}), $ \cB(F,{\rm f}^\eps_a) $ and $ \cB({\rm F}^\eps_a,f) $ (in Paragraph \ref{Influencebilinear}), 
as well as $ D $ (in Paragraph \ref{Impactperturbation}). To this end, we exploit the tools of Section \ref{RadonFourier}. New difficulties 
arise due to the presence of additional terms (especially those that have $ \eps^{-1} $ in factor), the absence of sign condition on 
$ f $, and the need to deal with the energy $ \EE $ (instead of $ \text{\small $ \pmb{\mathscr{E}} $} $ or $ \text{\small $ \mathscr{E} $ } $). 
The estimates below are not meant to be optimal (in terms of powers of $ \eps $ or $ P $), except for the crucial (singular) contribution 
$ \cB(\nu \times {\rm B}^\eps_a,f) $, where a compensation between the (possible) large size of $ {\rm B}^\eps_a $ and the $ L^2 $-smallness 
of $ f $ must be implemented. The proof of Theorem \ref{mainmaintheo} is completed in last Paragraph \ref{Conclusion}. Recall that $ 0 \leq t 
\leq T^\eps \leq \tilde \cT \leq \min (1,\cT_a) $.


\subsubsection{Control of the approximate momentum increment $ {\rm D}^\eps_a $} \label{approximateincrement} 
For illustrative purposes, we first examine the case of $ \tilde {\rm D}^\eps_a $, which is $ {\rm D}^\eps_a $ with 
$ {\rm U}^\eps_a \equiv \tilde {\rm U}^\eps_a $. Then, we turn to the general situation.

\smallskip

\noindent $ \bullet $ {\it Study of} $ \tilde {\rm D}^\eps_a $. When $ {\rm U}^\eps_a \equiv \tilde {\rm U}^\eps_a $ with 
$ \tilde {\rm U}^\eps_a $ as in (\ref{prototype}), as already seen, we have $ {\rm E}^\eps_h \equiv 0 $ and therefore 
$ \tilde {\rm D}^\eps_{ha} \equiv 0 $. Coming back to (\ref{cD}), we have to deal with
\[ \begin{array}{l}
 \displaystyle \tilde {\rm D}^\eps_{0a} :=  \int_0^t \int_{\mathbb S^2} \int_{\RR^3} {\rm W}_0 (s,\omega,\xi) \ {\rm M} (\eps, \langle 
 \xi \rangle) \ ds \, d \omega \, d \xi , \\
\displaystyle \tilde {\rm D}^\eps_{la} := \int_0^t \Bigl( \int_r^t \int_{\mathbb S^2} \int_{\RR^3} {\rm W}_l (s,\omega,\xi) \  
{\rm M} (\eps, \langle \xi \rangle) \ ds \, d \omega \, d \xi \Bigr) \, dr , \\
\displaystyle \tilde  {\rm D}^\eps_{na} := \frac{1}{\eps} \, \int_0^t \Bigl( \int_r^t \int_{\mathbb S^2} \int_{\RR^3} {\rm W}_n 
(r,s,\omega,\xi) \cdot \bigl( \nu(\xi) \times {\rm B}_e (x) \bigr) \ {\rm M} (\eps, \langle \xi \rangle) \ ds \, d \omega \, d \xi \Bigr) \, dr .
 \end{array} \]
In view of (\ref{the weights }), the functions $ {\rm W}_0 (s,\cdot) $ and $ {\rm W}_l (s,\cdot) $ are odd (with respect to both 
variables $ \omega $ and $ \xi $). It follows that $ \tilde {\rm D}^\eps_{0a} = 0 $ and $ \tilde {\rm D}^\eps_{la} = 0 $. On the 
other hand, we can replace $ {\rm W}_n  $ as indicated in (\ref{weightsn}), and then integrate by parts with respect to $ \xi $. Since
\begin{equation}\label{annulationaimposer}  
\nabla_\xi \cdot \bigl \lbrack \bigl( \nu (\xi) \times {\rm B_e} (x) \bigr) \ {\rm M}(\eps, \langle \xi \rangle)\bigr \rbrack = 
\bigl( \nu (\xi) \times {\rm B_e} (x) \bigr) \cdot \nabla_\xi \bigl \lbrack {\rm M}(\eps, \langle \xi \rangle) \bigr \rbrack = 0 , 
\end{equation}
we can assert that $ \tilde  {\rm D}^\eps_{na} = 0 $. Briefly, the stationary solution $ \tilde {\rm U}^\eps_a = 0 $ does not 
contribute to the momentum increment. We find that  $ \tilde  {\rm D}^\eps_a = 0 $.

\smallskip

\noindent $ \bullet $ {\it Study of} $ {\rm D}^\eps_a $. Let $ {\rm U}^\eps_a $ be an approximate solution in the sense 
of Definition \ref{Local well-prepared approximate solution}. From (\ref{comparisonboundgensuitefor0}) and 
(\ref{cestpourcdiameprem}) together with Lemma \ref{gainaveraging}, using (\ref{supppourUepsa}) and 
(\ref{Linftyteforuepsaf}), we can assert that
\[ \vert {\rm D}^\eps_{0a} \vert \lesssim \int_0^t \int_{\vert \xi \vert \leq {\rm P}_a} s \ {\rm S}_a \ ds \, d\xi \lesssim 
 {\rm P}_a^3 \ {\rm S}_a \ t^2 \, , \qquad \vert {\rm D}^\eps_{la} \vert \lesssim \int_0^t \Bigl( \int_r^t \int_{\vert \xi \vert 
 \leq {\rm P}_a} {\rm S}_a \ ds \, d\xi \Bigr) \, dr \lesssim  {\rm P}_a^3 \ {\rm S}_a \ t^2. \]
From (\ref{Linftyteforuepsah}), we easily get that $ \vert {\rm D}^\eps_{ha}  \vert \leq {\rm S}_h $. The most problematic 
term is $ {\rm D}^\eps_{na} $. Decompose $ {\rm F}^\eps_a $ as in (\ref{decompofexpa}). Exploit again the specific 
gradient form of $ {\rm W}_n $ inside (\ref{weightsn}) to perform an integration by parts with respect to the variable $ \xi $ 
in order to get 
\[ \begin{array}{rl} 
\displaystyle {\rm D}^\eps_{na} = \int_0^t \Bigl( \int_r^t \int_{\mathbb S^2} \int_{\RR^3} \! \! \! & \displaystyle \frac{s-r}{4 \pi} \ 
\frac{\nu \circ \Xi (s) \cdot  \bigl(\nu(\xi) + \omega \bigr)}{1 + \omega \cdot \nu(\xi)} \\
& \bigl( {\rm E}^\eps_a + \nu(\xi) \times {\rm B}^\eps_a \bigr) \cdot \nabla_\xi {\rm f}^\eps_a \bigl(r,X(s) + (s-r) 
\omega,\xi \bigr) \ ds \, d \omega \, d \xi \Bigr) \, dr .
 \end{array} \]
From Lemmas \ref{gainaveraging} and \ref{Comparisonsingularweights} together with (\ref{supppourUepsa}), 
(\ref{LinftyteforuepsaE}), (\ref{Lipesteforuepsaf}) and (\ref{wellwellprepa}), we can assert that
\[ \vert {\rm D}^\eps_{na} \vert \lesssim \int_0^t \Bigl( \int_r^t \int_{\vert \xi \vert \leq {\rm P}_a} (s-r) \ \bigl( {\rm S}_r + 
2 \, (1+ {\rm S}_a) \, {\rm L}_a \bigr) \ ds \, d\xi \Bigr) \, dr \lesssim {\rm P}_a^3 \ \bigl( {\rm S}_r + 2 \, (1+ {\rm S}_a) \, 
{\rm L}_a \bigr) \ t^3 . \]
In the end, there remains 
\begin{equation}\label{conconfroninverseepsa}  
 \vert {\rm D}^\eps_a \vert \lesssim {\rm S}_h + {\rm P}_a^3 \ {\rm S}_a \ t^2 +  {\rm P}_a^3 \ \bigl( {\rm S}_r + 2 \, 
 (1+ {\rm S}_a) \, {\rm L}_a \bigr) \ t^3 .
\end{equation}


\subsubsection{Influence of the frozen bilinear terms $ \cB(F,{\rm f}^\eps_a) $ and $ \cB({\rm F}^\eps_a,f) $} \label{Influencebilinear} 
We assume here that $ t  \leq T^\eps $.

\smallskip

\noindent $ \bullet $ {\it Control of} $ \cB(F,{\rm f}^\eps_a) $. We can use (\ref{decompofexpa}) to get $ \cB(F,{\rm f}^\eps_a) 
= \cB(E,{\rm f}^\eps_a) + \cB(\nu \times B,{\rm f}^\eps_a) $. From (\ref{supppourUepsa}) and (\ref{Linftyteforuepsaf}), it is 
clear that
\[ \begin{array}{rl}
\vert \cB(E,{\rm f}^\eps_a) \vert \leq \! \! \! & \displaystyle \int_0^t \Bigl( \int_r^t \int_{\mathbb S^2} 
\int_{\vert \xi \vert \leq {\rm P}_a} {\rm S}_a \ J^{-1/2} \ \vert {\rm W}_{ne} \vert \ J^{1/2} \ \vert E \vert \bigl(r,X(s) + (s-r) 
\omega \bigr) \ ds \, d \omega \, d \xi \Bigr) \, dr \\
\lesssim \! \! \! & \displaystyle {\rm S}_a \ {\rm P}_a^{3/2} \, \int_0^t \EE (r)^{1/2} \ \Bigl( \int_r^t \int_{\mathbb S^2} 
\int_{\vert \xi \vert \leq {\rm P}_a} J^{-1} \ \vert {\rm W}_{ne} \vert^2 \ ds \, d \omega \, d \xi \Bigr)^{1/2} \, dr . 
\end{array}  \]
From (\ref{jacobiadit}) and (\ref{cestpoureeeee}) together with Lemma \ref{gainaveraging}, since $ s-r $ appears in factor 
inside $ {\rm W}_n \equiv {\rm W}_{ne} $, we have
\begin{equation}\label{ahfgiaghf}   
\begin{array}{rl}
\displaystyle \int_r^t \int_{\mathbb S^2} \int_{\vert \xi \vert \leq {\rm P}_a} J^{-1} \ \vert {\rm W}_{ne} \vert^2 \ ds \, 
d \omega \, d \xi \leq \! \! \! & \displaystyle \int_r^t \int_{\vert \xi \vert \leq {\rm P}_a} \Bigl( \int_{\mathbb S^2} 
\frac{d \omega}{\langle \xi \rangle^2 \, \bigl( 1 + \omega \cdot \nu (\xi) \bigr)^2} \Bigr) \, ds \, d \xi \\
\lesssim \! \! \! & \displaystyle \int_r^t \int_{\vert \xi \vert \leq {\rm P}_a} ds \, d \xi \lesssim t \ {\rm P}_a^3 . 
\end{array} 
\end{equation}
Thus, coming back to (\ref{propaL2final}), we recover that $ \vert \cB(E,{\rm f}^\eps_a) \vert \lesssim {\rm S}_a 
\, {\rm P}_a^3 \,  t^{3/2} \, \eps $. Since $ \vert \nu \times B \vert \leq \vert B \vert $, the same argument  
applies to $ \cB(\nu \times B,{\rm f}^\eps_a) $. Thus, we can retain that
\begin{equation}\label{conconfronE}  
\vert \cB(F,{\rm f}^\eps_a) \vert \lesssim {\rm S}_a \, {\rm P}_a^3 \ {\rm N} \ t^{3/2} \, \eps \lesssim {\rm S}_a \, 
{\rm P}_a^3 \ {\rm N} \ t .
\end{equation}

\smallskip

\noindent $ \bullet $ {\it Control of} $ \cB({\rm F}^\eps_a,f) $. From (\ref{LinftyteforuepsaE}) and (\ref{LipesteforuepsaB}), 
we can infer that
\[ \begin{array}{rl}
\vert \cB({\rm F}^\eps_a,f) \vert \leq \! \! \! & \displaystyle ({\rm S}_a + \eps^{-1} \, {\rm H}_a) \, \int_0^t 
\Bigl( \int_r^t \int_{\mathbb S^2} \int_{\RR^3} \vert {\rm W}_{ne} \vert \ \vert f \vert \bigl(r,X(s) + (s-r) \omega , \xi \bigr) \ 
ds \, d \omega \, d \xi \Bigr) \, dr \\
\lesssim \! \! \! & \displaystyle ({\rm S}_a + \eps^{-1} \, {\rm H}_a) \, \int_0^t \Bigl( \int_r^t \int_{\mathbb S^2} 
\int_{\vert \xi \vert \leq P(r)} J^{-1} \ \vert {\rm W}_{ne} \vert^2 \ ds \, d \omega \, d \xi \Bigr)^{1/2} \\
\ & \displaystyle \qquad \qquad \qquad \quad \ \ \times \Bigl( \int_r^t \int_{\mathbb S^2} \int_{\RR^3} J \ f^2 \bigl(r,X(s) + 
(s-r) \omega , \xi \bigr) \ ds \, d \omega \, d \xi \Bigr)^{1/2}\, dr . 
\end{array}  \]
With $ {\rm P}_a $ replaced by $ P(r) $, we proceed as in (\ref{ahfgiaghf}) to find that
\[ \vert \cB({\rm F}^\eps_a,f) \vert \lesssim ({\rm S}_a + \eps^{-1} \, {\rm H}_a) \, \int_0^t t^{1/2} \, P(r)^{3/2} \ 
\EE (r)^{1/2} \ dr . \]
Knowing (\ref{propaL2final}),  we obtain that
\begin{equation}\label{conconfronf}  
\vert \cB({\rm F}^\eps_a,f) \vert \lesssim (\eps \, {\rm S}_a + {\rm H}_a) \ {\rm N} \ t^{1/2} \, \int_0^t P(r)^{3/2} \ dr 
\lesssim ({\rm S}_a + {\rm H}_a) \ {\rm N} \ \Bigl( t + \int_0^t P(r)^3 \, dr \Bigr) .
\end{equation}


\subsubsection{Impact of the perturbed momentum increment $ D $} \label{Impactperturbation} $ \, $ We assume 
again that $ t  \leq T^\eps $. We consider successively $ D_0 $, $ D_l $ and $ D_n$.

\smallskip

\noindent $ \bullet $ {\it Control of} $ D_0 $. From (\ref{minimalradiusforf}), (\ref{Lebesgue bounds}), (\ref{comparisonboundgensuitefor0}) 
and Lemma \ref{gainaveraging}, we have
\[ \vert D_0 \vert \lesssim \int_0^t \int_{\vert \xi \vert \leq {\rm P}_0} s \ {\rm S}_0 \ ds \, d \xi \lesssim {\rm P}_0^3 \ {\rm S}_0 \ t^2 . \]

\noindent $ \bullet $ {\it Control of} $ D_l $. Using (\ref{supnormonf}), and then (\ref{cestpourcdiameprem}) together with 
Lemma \ref{gainaveraging}, we can be satisfied with
\[ \vert D_l \vert \lesssim \int_0^t \Bigl( \int_r^t \int_{\vert \xi \vert \leq P (r)} \ ({\rm S}_0 + {\rm S}_a) \ ds \, d\xi \Bigr) \, dr 
\lesssim ({\rm S}_0 + {\rm S}_a) \ t \ \int_0^t P (r)^3 \ dr . \]

\noindent $ \bullet $ {\it Control of} $ D_n $. Recall that $ D_n = \cB(E,f)+  \cB(\nu \times B,f) $. We handle below
$ \cB(E,f) $, the case of $ \cB(\nu \times B,f) $ being completely similar. By Cauchy-Schwarz inequality, we have
\[ \begin{array}{rl}
\vert D_n \vert \leq \! \! \! & \displaystyle \int_0^t \Bigl( \int_r^t \int_{\mathbb S^2} \int_{\vert \xi \vert \leq P(r)} J \ 
\vert E \bigl(r,X(s) + (s-r) \omega \bigr) \vert^2 \ ds \, d \omega \, d \xi \Bigr)^{1/2} \\
\ & \displaystyle \qquad \qquad \qquad \times \Bigl( \int_r^t \int_{\mathbb S^2} \int_{\vert \xi \vert \leq P(r)} J^{-1} \
\vert {\rm W}_{ne} \vert^2 \ f^2 \ ds \, d \omega \, d \xi \Bigr)^{1/2} \, dr \\
\lesssim \! \! \! & \displaystyle \int_0^t \EE (r)^{1/2} \ P(r)^{3/2} \ ({\rm S}_0 + {\rm S}_a) \,  \Bigl( \int_r^t \int_{\vert 
\xi \vert \leq P(r)} \bigl( \int_{\mathbb S^2}  J^{-1} \ \vert {\rm W}_{ne} \vert^2 \ d \omega \bigr) \, ds \, d \xi \Bigr)^{1/2} \, dr .
\end{array}  \]
With $ {\rm P}_a $ replaced by $ P(r) $, we carry on with (\ref{ahfgiaghf}). Then, from (\ref{propaL2final}), we get 
\[ \vert D_n \vert \lesssim \eps \ {\rm N} \ ({\rm S}_0 + {\rm S}_a) \ t^{1/2} \ \int_0^t P (r)^3 \ dr . \]
In brief, we can retain that 
 \begin{equation}\label{----} 
  \vert D \vert \lesssim {\rm P}_0^3 \ {\rm S}_0 \ t^2 + ({\rm S}_0 + {\rm S}_a) \ (1+ t) \, \int_0^t P (r)^3 \ dr .
    \end{equation}
    

\subsubsection{Compilation of the preceding estimates} \label{Conclusion} We still have (\ref{compactness arguments})
for some $ (y_0,\eta_0) $ in the support of $ {\rm f}_0 $, that is for some $ \eta_0 $ satisfying $ \vert \eta_0 \vert \leq 
\max \, ( {\rm P}_a , {\rm P}_0 ) $. In view of (\ref{broughtbypoursuit}), this gives rise to
\[ \langle {\rm P} \rangle_\infty (t) \leq 1 + \max \ \bigl( {\rm P}_a , P_0 \bigr) + \vert {\rm D}^\eps_a \vert 
+ \vert \cB(F,{\rm f}^\eps_a) \vert + \vert \cB({\rm F}^\eps_a,f) \vert + \vert D \vert . \]
By construction, we have
 \begin{equation}\label{finalpourP} 
 {\rm P}(t) \leq \max \ \bigl( {\rm P}_a , P(t) \bigr) \, , \qquad  P(t) \leq \max \ \bigl( {\rm P}_a , {\rm P} (t) \bigr) ,
\end{equation}
and similar inequalities concerning $ \langle {\rm P} \rangle_\infty (t)  $ and $ \langle P \rangle_\infty (t) $. In particular
\[ \langle P \rangle_\infty (t) \leq \max \ \bigl( {\rm P}_a , \langle {\rm P} \rangle_\infty (t) \bigr) \leq 1 + \max \ \bigl( {\rm P}_a , 
P_0 \bigr) + \vert {\rm D}^\eps_a \vert + \vert \cB(F,{\rm f}^\eps_a) \vert + \vert \cB({\rm F}^\eps_a,f) \vert + \vert D \vert  . \] 
From (\ref{conconfroninverseepsa}), (\ref{conconfronE}), (\ref{conconfronf}) and (\ref{----}), together with 
the definition of $ {\rm I}_0 $ inside (\ref{cestpourI0}), we can easily infer that
\[ \langle P \rangle_\infty (t) \leq {\rm I}_0 + \tilde {\rm C} \ t + \check {\rm C} \, \int_0^t \langle P \rangle_\infty (r)^3 \, dr , \]
where the constants $ \tilde {\rm C} $ and $ \check {\rm C} $ depend only on $ {\rm L}_a $, $ {\rm H}_a $, the $ {\rm S}_\star $, 
the $ {\rm P}_\star $ and the $ {\rm N}_\star $. Introduce
\[ \beta := 16 \ \check {\rm C} \ {\rm I}_0^2 \, , \qquad \cT := \min \ \bigl( \tilde \cT , {\rm I}_0/ 
\tilde {\rm C} , (\ln 3 - \ln 2)/ \beta \bigr) \, , \qquad \tilde T^\eps := \min \ ( T^\eps, \cT )  . \]
Taking into account (\ref{startingfrombeg}), we can assert that
\[ \langle P \rangle_\infty (t) \leq 2 \ {\rm I}_0 + \beta \, \int_0^t \langle P \rangle_\infty (r) \, dr \, , \qquad 
\forall \, t \in [0,\tilde T^\eps] . \]
Then, by Gr\" onwall's inequality, the quantity $ \langle P \rangle_\infty (r) $ remains controlled according to
 \begin{equation}\label{contradictionourP} 
 \langle P \rangle_\infty (t) \leq \cF (t) := 2 \ {\rm I}_0 \ e^{\beta \, t} \leq 3 \ {\rm I}_0 \, , \qquad \forall \, t \in 
[0,\tilde T^\eps] . 
\end{equation}
If $ T^\eps < \cT \leq \tilde \cT $ so that $ \tilde T^\eps = T^\eps $, due to the definition of $ T^\eps $ just before the line
(\ref{propaL2final}), we must have $ \langle P \rangle_\infty (T^\eps ) = 4 \, {\rm I}_0 $. This is clearly a contradiction 
with (\ref{contradictionourP}) for $ t = \tilde T^\eps = T^\eps $. Necessarily, we must have $ \cT \leq T^\eps $, and 
(\ref{contradictionourP}) holds true on $ [0, \cT ] $ as required in Theorem \ref{mainmaintheo}. Of course, the function 
$ \cF $ inside (\ref{conofcp}) depends also on the various parameters $ \dagger_a $, $ {\rm S}_h $ and $ \dagger_r $ 
occurring in Definition \ref{Local well-prepared approximate solution}. However, since the family $ \{ {\rm U}^\eps_a 
\}_\eps $ is viewed as being fixed, this influence has not been reported. On the contrary, we can choose any perturbation 
$ \{ f^\eps_0 \}_\eps $ as long as it is controlled by $ {\rm P}_0 $, $ {\rm S}_0 $ and $ {\rm N}_0 $ as indicated in 
Definition \ref{assumptiononofoeps}. This is why we have highlighted inside (\ref{conofcp}) the impact of $ ({\rm P}_0 , 
{\rm S}_0 , {\rm N}_0 ) $.

\section*{Acknowledgments}
Part of this material is based upon work done while S.I. was supported by the lab IRMAR at Universit\'e of Rennes. 
He would like to thank all the members of the lab for their great hospitality. S. Ibrahim is supported by the NSERC grant No. 371637-2019.


\bibliographystyle{abbrv}
\bibliography{RVM_19_06_23.bib}

\end{document}